\documentclass[reqno]{amsart}
\usepackage{amsmath,amssymb}
\theoremstyle{plain}
\newtheorem{theorem}{Theorem}
\newtheorem{prop}[theorem]{Proposition}
\newtheorem{cor}[theorem]{Corollary}
\theoremstyle{definition}
\newtheorem*{definition}{Definition}
\newtheorem*{example}{Example}
\newtheorem{lemma}[theorem]{Lemma}
\theoremstyle{remark}
\newtheorem*{remark}{Remark}
\numberwithin{theorem}{section} \numberwithin{equation}{section}
\begin{document}
\author{Oleg Chalykh$^1$}
\thanks {$^1$ On leave of absence from:  Advanced Education and Science Centre, Moscow State
University, Moscow 119899, Russia}
\address{O.C.: Department of Mathematics,
Cornell University, Ithaca, NY 14853, USA}
\email{oleg@math.cornell.edu} 

\title{Macdonald Polynomials And  Algebraic Integrability}
\begin{abstract}
We construct explicitly (non-polynomial) eigenfunctions of the
difference operators by Macdonald in case $t=q^k$,\, $k\in{\mathbb
Z}$. This leads to a new, more elementary proof of several
Macdonald conjectures, first proved by Cherednik. We also
establish the algebraic integrability of Macdonald operators at
$t=q^k$\, ($k\in {\mathbb Z}$), generalizing the result of Etingof
and Styrkas. Our approach works uniformly for all root systems
including $BC_n$ case and related Koornwinder polynomials.
Moreover, we apply it for a certain deformation of $A_n$ root
system where the previously known methods do not work.
\end{abstract}

\maketitle
\section{Introduction}
One of the goals of this paper is to present a new (essentially,
"non-polynomial") approach to Macdonald polynomials. These
polynomials were introduced in late 80-s by I.G.Macdonald in
\cite{M1} as, informally speaking, a "discrete spectrum" of
certain remarkable symmetric difference operators. Since these
(Macdonald) operators are self-adjoint with respect to a specific
scalar product, Macdonald's theory leads directly to the families
of multivariable orthogonal polynomials. From that point of view,
they generalize various classical orthogonal polynomials of one
variable. In fact, there are as many families of Macdonald
polynomials, as simple complex Lie algebras, or root systems. Each
family depends, apart from a root system, on (at least) two
parameters $q,t$ and specializes to several remarkable families of
symmetric functions. Among them are Schur functions and characters
of corresponding Lie groups, Hall--Littlewood functions, Jack
polynomials or, more generally, multivariable Jacobi polynomials
due to Heckman and Opdam \cite{HO}. All this makes Macdonald
polynomials very interesting from the point of view of the
representation theory, combinatorics, special function theory and
mathematical physics. This also makes clear that Macdonald
polynomials are highly non-trivial. Thus, it is not surprising
that their various properties formulated by Macdonald as
conjectures, remained unproven for quite a long time. A remarkable
progress has been achieved by Cherednik, who proved the so-called
norm and evaluation conjectures and the symmetry identity for all
root systems \cite{Ch1,Ch2}. Cherednik's approach is based on his
theory of double affine Hecke algebras and it remains one of the
major achievements in this area. As an introduction into
Cherednik's theory we recommend to the reader a nicely written
survey by Kirillov \cite{Ki1}.

One of the results of the present paper is an independent proof of
these three Macdonald's conjectures. Our approach uses some
remarkable properties of the Macdonald operators in case $t=q^k$
with {\it integer} $k$. Ideologically, it goes back to the paper
by Veselov and the author \cite{CV}, where the quantum
Calogero--Sutherland--Moser problem was considered for some
special values of the coupling constant. Recall that the
Calogero--Sutherland--Moser problem \cite{Ca,Su} describes $N$
particles $x_1,\dots,x_N$ on the line whose pairwise interaction
is given by the potential $u(x)=m(m+1)\sin ^{-2}x$. In the quantum
case its Hamiltonian is the following Schr\"odinger operator in
$\mathbb R^N$: $$H=-\Delta+\sum_{i<j} m(m+1)\sin
^{-2}(x_i-x_j)\,.$$ It is a celebrated example of a
completely integrable $N$-body problem, and there are quite a few
exact results about it. In particular, it is completely integrable
in Liouville sense, i.e. $H$ is a member of a family of $N$
commuting partial differential operators (quantum integrals)
$H_1=H,H_2,\dots, H_N$. Moreover, as it was demonstrated in
\cite{CV,VSC}, for special values of the coupling constant $m$
this problem becomes "much more" integrable. Namely, it turns out
that the Calogero--Sutherland--Moser problem for $m \in {\bf Z}$
is {\it algebraically integrable}, i.e.  its quantum integrals
$H_1,\ldots,H_N$ are a part  of some bigger commutative ring
$\mathfrak R$ of partial differential operators (see \cite{VSC}
for precise formulations and results). Moreover, $Spec(\mathfrak
R)$ is an affine algebraic variety whose points
parameterize Bloch eigenfunctions of $H$. This is a
multidimensional analogue of a phenomenon, well known from the
finite-gap theory in dimension one \cite{DMN, Kr}. The main
difference with the one-dimensional case is that the corresponding
algebro-geometrical data are very "rigid" and really exceptional,
which makes the existence problem for related multidimensional
Baker--Akhiezer functions extremely difficult. This difficulty was
overcome in \cite{VSC} with the help of the theory of
multivariable hypergeometric  functions due to Heckman and Opdam
\cite{HO}. Recently, a direct independent proof was obtained in
\cite{C00}.

Continuing  \cite{C00}, here we demonstrate that a similar
phenomenon appears for $H$ being replaced by any of the Macdonald
difference operators, namely, they all are {\it algebraically
integrable} for special integer values of the parameter(s). Note
that in case of $R$ being of $A_n$ type the Macdonald operators
coincide (up to a certain gauge) with the trigonometric version of
the elliptic Ruijsenaars operators from \cite{R}, introduced as a
generalization of the elliptic Caloger--Moser problem. We should
mention at this point the paper by Etingof and Styrkas \cite{ES}
where the algebraic integrability has been established for the
Macdonald--Ruijsenaars operators with $t=q^k$,\,$k \in {\bf Z}$.
Their approach was based on earlier results by Etingof and
Kirillov \cite{EK1, EK2} who gave an interpretation of Macdonald
polynomials for $R=A_n$ in terms of the representation theory for
quantum groups. This delivers independent proofs of several
results in case $R=A_n$, see \cite{EK3,EK4,Ki2}. It is worth
noticing that the symmetry identity in this case has been proved
first by Koornwinder \cite{Ko1} (see chapter VI of Macdonald's
book \cite{M4}). However, his method, as well as the methods of
\cite{EK2,ES}, does not extend to other root systems.

One of the advantages of our approach is that it works equally
well for all root systems and (we believe) is simpler comparing to
\cite{Ch1,Ch2}. The main object is what is natural to call a {\it
Baker--Akhiezer function} $\psi$ for Macdonald operators.
 In case $R=A_n$ this coincides with
the $\psi$-function from \cite{ES}. In that part which goes back
to the papers \cite{CV,VSC}, our considerations have very much in
common with \cite{ES}. The main new ingredient is a fairly
elementary construction of the $\psi$-function. Our main
observation (Proposition \ref{inv}) is that the Macdonald
operators in case $t=q^k$,\ $k\in\mathbb Z_-$, act naturally in
the coordinate ring of a certain very specific affine algebraic
variety. This implies the existence of $\psi$ which is our central
result. From this we derive the duality, which reflects a certain
symmetry between the two arguments $x,z$ of $\psi$-function.
Basically, it means that, as a function of $x$, $\psi$ is an
eigenfunction of the Macdonald operators related to the root
system $R$ , while, as a function of $z$, it is an eigenfunction
of the Macdonald operators related to the {\it dual} root system.
Thus, we observe on this level the {\it bispectrality} of
Macdonald operators, if one uses the terminology of the
fundamental work by Duistermaat and Gr\"unbaum on bispectral
problem \cite{DG}. It is worth mentioning that in contrast with
the one-dimensional case, bispectral problem in higher dimensions
is much less investigated. For an interesting example related to
Knizhnik--Zamolodchikov equation see recent paper \cite{TV}.

Notice that our proof of the  existence of $\psi$ is an effective
one and gives a closed expression for it. The resulting formula
generalizes the main result of \cite{C00} and it is a discrete
version of one remarkable formula by Berest, who found in
\cite{Be} an elegant "universal" expression for the
$\psi$-function for the quantum (rational) Calogero--Moser
problem. His derivation, however, was based on a crucial
assumption that such a $\psi$ {\it does exist}. As we mentioned
above, that type of existence problems is highly nontrivial in
dimension $>1$. Remarkably enough, our approach being inspired by
the Berest's result, allows us to do these two things
simultaneously: we prove the existence of $\psi$ by a direct
derivation of a discrete analogue of the Berest's formula.

Everything extends without much difficulties to $BC_n$ case. A
proper generalization of Macdonald's theory for this case was
suggested by Koornwinder \cite{Ko1}. The resulting orthogonal
polynomials (Koornwinder polynomials) can be viewed as a
multivariable analog of the celebrated Askey--Wilson polynomials
\cite{AW}. Van Diejen \cite{vD1} showed that Koornwinder
polynomials are joint eigenfunctions of $n$ commuting difference
operators for which he gave an explicit expression (initially,
Koornwinder constructed one operator only). Further, in \cite{vD2}
the symmetry identity (or {\it self-duality}) was established for
a certain subfamily of Koornwinder polynomials. Then, finally,
Sahi \cite{S} proved duality in general case, using a proper
generalization of Cherednik's double affine Hecke algebra.
Together with earlier van Diejen's results \cite{vD2} this implied
also the evaluation identity and the norm identity, conjectured by
Macdonald. Our approach leads to an independent proof of these
results. As above, the key ingredient is the algebraic
integrability of the Koornwinder operator for special integer
values of the parameters.

 One of our primary motivations for this work was, in fact, our
attempt to find a difference version of the deformed
Calogero--Moser problem from \cite{VFC}. It is related to what can
be viewed as a one-parameter deformation of $A_n$ root system.
 In \cite{C97}, guided by
duality, we were able to find a rational difference version of
that quantum problem. Here we consider its trigonometric version,
proving its (algebraic) integrability. Similar to the usual $A_n$
case, the constructed difference operator is self-dual:
corresponding Baker--Akhiezer function $\psi(x,z)$ is invariant
under permuting $x,z$. A natural elliptic version seems to be
integrable, too. We would like to stress that while for the root
systems the approach based on affine Hecke algebras seems to be
the most adequate and powerful, in the deformed case none of the
previously known methods can be applied (at least,
straightforwardly). Thus, it would be very interesting to find an
algebraic structure which underlies that "deformed" root system.

The paper is organized as follows.  In Section 2 we recall the
definitions of difference operators and polynomials due to
Macdonald. Then we make our central observation about Macdonald
operators in case $t=q^{-m}$ ($m\in\mathbb Z_+$). The following
two sections form the core of the paper.

In Section 3, we define a Baker--Akhiezer function $\psi(x,z)$
associated to a datum which consists of a root system $R$ and some
additional integer parameter(s). $\psi$ is determined by
prescribing its analytic properties in the $z$-variable, which is
a proper modification of the approach from \cite{CV,ES}. We prove
first its uniqueness up to a normalization and then the existence,
obtaining as a by-product a discrete version of the Berest's
formula.

In Section 4, we explain how one should normalize $\psi$ to
achieve a remarkable symmetry between $x$ and $z$. This leads us
directly to the duality theorem, which is the central result of
this section.

In Section 5 we derive various corollaries from achieved results.
First, we obtain algebraic and complete integrability of the
Macdonald operators and prove the existence of the so-called shift
operators. Then we explain how our $\psi$ relates to Macdonald
polynomials, this generalizes Weyl's character formula and is
similar to a relation between symmetric and nonsymmetric Macdonald
polynomials (see \cite{M3}). As a corollary, we observe a nice
"localization" property for Macdonald polynomials in case
$t=q^{-m}$ ($m\in\mathbb Z_+$). We conclude the section by
explaining how our results lead to a proof of the norm identity,
evaluation formula and duality for Macdonald polynomials.

In Section 6 we sketch how to extend our approach to $BC_n$ case.
This is parallel to the previous sections, so we skip most of the
proofs. The main difference comes from $n=1$ case which is
technically more difficult compared to $A_1$.

Finally, in Section 7 we discuss a generalized
Macdonald--Ruijsenaars model, related to the deformed $A_n$
system.

\bigskip

\noindent{\bf Acknowledgments}. I am grateful to Yu.Berest and
A.P.Veselov for stimulating discussions. The work was supported by
EPSRC (grant GR/M69548).

\section{Difference operators and polynomials by Macdonald}
\subsection{Notations}
Let $V=\mathbb R^n$ be a Euclidean  space with  the scalar product
denoted as $(u,v)$. Consider an arbitrary root system $R\in V$
which is, by definition, a finite set of vectors (roots) $\alpha
\in V$ with the following two properties:

(1) $\forall \alpha \in R$ the orthogonal reflection $s_\alpha$
in $V$
$$
s_\alpha \,:\,x\longmapsto x-2\frac
{(\alpha,x)}{(\alpha,\alpha)}\alpha
$$

leaves $R$ invariant, $s_\alpha (R)=R$;

(2) $\forall \alpha, \beta \in R$\quad $2\frac {(\alpha, \beta)}{(\alpha, \alpha)}
\in \mathbb Z$

\noindent (see \cite{B} for the details).

The second property implies that the so-called root lattice $Q$
generated over $\mathbb Z$ by the roots $\alpha \in R$ is
invariant under all the reflections $s_\alpha$ and, therefore,
under the whole Weyl group $W$ generated by $s_\alpha$\,, $\alpha
\in R$. The vectors $\alpha ^\vee = 2\alpha / (\alpha, \alpha)$
form the dual root system $R^\vee$ and we denote by $Q^\vee$ the
lattice generated by all $\alpha ^\vee \in R^\vee$. Introduce also
weight and coweight lattices $P, P^\vee$ as
\begin{align*}
P &= \{ \pi \in V \vert \left( \pi ,\alpha^\vee\right) \in \mathbb
Z \quad \forall \alpha^\vee \in R^\vee \}&&\text{(weights)}\,,\\
P^\vee &= \{ \pi \in V \vert \left( \pi ,\alpha\right) \in \mathbb
Z \quad \forall \alpha \in R \}&&\text{(coweights)}\,.
\end{align*}
{From} the definitions one has inclusions $Q \subset P$, $Q^\vee
\subset P^\vee$. Taking ${\mathbb Z}_{\ge 0}$ instead of ${\mathbb
Z}$ in the last two definitions leads to dominant weights
(coweights) $P_+$ and $P_+^\vee$, respectively.

Let us fix some basis of simple roots $\alpha_1, \dots , \alpha_n$
in $R$, this determines a decomposition of $R$ and $R^\vee$ into
positive and negative parts:
$$
R=R_+\cup (-R_+)\,,\qquad R^\vee =R_+^\vee \cup (-R_+^\vee )\,.
$$
The elements $\omega_i$ of the basis, dual to $\alpha_1^\vee , \dots ,
\alpha_n^\vee$,\, $(\omega_i, \alpha_j^\vee)=\delta _{ij}$\,, are
called the {\it fundamental weights} for the system $R$.
Similarly, one defines the fundamental coweights $b_i$: $(b_i,
\alpha_j)=\delta_{ij}$\,.

In these terms the root lattice $Q$ is simply
$$Q=\bigoplus_{i=1}^n \mathbb Z\alpha_i\,.$$ Its {\it positive
part} $Q_+$ is obtained by replacing $\mathbb Z$ by $\mathbb Z_+$.
Similarly, the cone of the dominant weights $P_+$ is
$$P_+=\bigoplus_{i=1}^n \mathbb Z_+\omega_i\,.$$

Below we will assume that $R=A_n,\dots,G_2$ is reduced and
irreducible. The case $R=BC_n$ is considered in section $6$.

\subsection{Macdonald operators}\label{beg}
Let $R$ be a (reduced, irreducible) root system in Euclidean space
$V$ . Let us fix $q\in \mathbb C^\times$ and a set $k$ of
 $W$-invariant parameters $k_\alpha = k_{\alpha^\vee}$,
i.e. $k_\alpha = k_{w\alpha}$ for any $\alpha \in R$ and $w\in W$.
Below we will sometimes use the related parameters
$t_\alpha=t_{\alpha^\vee} :=q^{k_\alpha}$, denoting by $t=q^k$ the
whole set $\{t_\alpha\}$\footnote{To define $q^k$ for arbitrary
$k$ we fix the value of $\log q$, so $q^k:=e^{k\log q}$.}.
Throughout the paper we will suppose that $q$ is not a root of
unity. Something still can be done in case of roots of unity, though we
will not touch these issues here (see \cite{Ki2,Ch2}).

For $v \in V$ we denote by $T^v$ the operator acting on a function
$f(x)$ as a shift by $v$ in $x\in V$: $$ \left(T^v
f\right)(x)=f(x+v)\,. $$ Later, we will deal with the functions
depending on two variables $x,z \in V$, in that case we will use
subscripts to distinguish between shifts $T^v_x$ and $T^v_z$,
acting in $x$ and $z$, respectively. To introduce Macdonald
operators, we need the notion of a (quasi)minuscule coweight.

\bigskip

\begin{definition}
(1) A coweight $\pi \in P^\vee$ is called minuscule if $-1\le (\pi,\alpha)\le 1$
for all $\alpha \in R$.

(2) A coweight $\pi \in P^\vee$ is called quasiminuscule if it
belongs to $R^\vee$ and $-1\le (\pi,\alpha)\le 1$
for all $\alpha \in R \backslash \pi^\vee$.
\end{definition}

\bigskip

Using tables from \cite{B} one can check that all root systems
except $E_8, F_4, G_2$ have at least one nonzero minuscule
coweight. Meanwhile, for any root system $R$ the coroot
$\pi=\theta^\vee$ where $\theta$ is the {\it maximal root} for
$R$, will be quasiminuscule (see \cite{B}).

Now let $\pi \in P^\vee$ be a minuscule coweight for the system
$R$. Then the corresponding Macdonald operator $D^\pi$ is a
difference operator in $x\in V$ defined as follows \cite{M1,M2}:
\begin{equation}
\label{m} D^\pi = \sum_{\tau\in W\pi}a_\tau T^{\tau}, \quad a_\tau
(x) = \prod_{\genfrac{}{}{0pt}{}{\alpha \in R :}{(\alpha,
\tau)>0}}\, \frac {t_\alpha q^{(\alpha,x)}-t_\alpha^{-1}
q^{-(\alpha,x)}} {q^{(\alpha,x)}-q^{-(\alpha,x)}}\,.
\end{equation}
For a quasiminuscule coweight $\pi$ the formula is slightly more
complicated:
\begin{equation}
\label{m1} D^\pi = \sum_{\tau\in W\pi}a_\tau
\left(T^{\tau}-1\right) + \sum_{\tau\in W\pi} q^{-2(\rho,\tau)}\,,
\end{equation}
where
\begin{equation}
\label{m2} a_\tau (x) = \prod_{\genfrac{}{}{0pt}{}{\alpha \in R
:}{(\alpha,\tau)>0}} \frac {t_\alpha
q^{(\alpha,x)}-t_\alpha^{-1}q^{-(\alpha,x)}}
{q^{(\alpha,x)}-q^{-(\alpha,x)}} \prod_{\genfrac{}{}{0pt}{}{\alpha
\in R :}{(\alpha, \tau)=2}} \frac {t_\alpha
q^{1+(\alpha,x)}-t_\alpha^{-1} q^{-1-(\alpha,x)}}
{q^{1+(\alpha,x)}-q^{-1-(\alpha,x)}}
\end{equation}
and
\begin{equation}
\label{rho} \rho = \rho_k=\frac 12 \sum_{\alpha \in R_+}k_{\alpha}
\alpha\,.
\end{equation}

\medskip
\begin{remark}\label{cover}
 Note that in the
formula \eqref{m2} for $a_\tau$ the last product contains one
factor only (with $\alpha=\tau^\vee$). Written this way, it
formally makes sense for a minuscule coweight, too. Indeed, a
formal substitution of a minuscule $\pi$ into \eqref{m1} leads to
\eqref{m} because the constant term $\sum q^{-2(\rho,\tau)}-\sum
a_\tau$ will be zero in that case. In the remainder of the paper,
we will refer to the formula \eqref{m1} since it covers
both cases.
\end{remark}
\medskip
\begin{remark}
The following function $\Delta_k(x)$ plays an important role in
Macdonald's theory:
\begin{equation}\label{delta}
\Delta=\Delta_{k}(x)=q^{-2(\rho,x)}\prod_{\alpha\in
R_+}\prod_{i=0}^{\infty}\frac{1-q^{2i+2(\alpha,x)}}
{1-q^{2k_\alpha+2i+2(\alpha,x)}}\,.
\end{equation}
Using it, one can present the coefficients $a_\tau$ of the
Macdonald operator $D^\pi$ as follows:
\begin{equation}\label{coef}
a_\pi=T^\pi(\Delta)/\Delta\,,\qquad
a_{w\pi}(x)=a_{\pi}(w^{-1}x)\,.
\end{equation}
\end{remark}

\bigskip

\begin{example}
In case $R = A_{n-1} = \{ \pm\left(e_i-e_j\right) \vert i<j \}
\subset {\mathbb R}^{n}$ with $k_\alpha \equiv k$ each fundamental
coweight $\pi_s = e_1+\dots +e_s$ ($s=1,\dots, n$) is minuscule
and the corresponding operator $D_s=D^{\pi_s}$ has the form
\begin{equation}
\label{ru} D_s=\sum_{\genfrac{}{}{0pt}{}{I\subset
\{1,\dots,n\}}{|I|=s}} \prod_{\genfrac{}{}{0pt}{}{i\in I} {j\notin
I}} \frac {q^{k+x_i-x_j}-q^{-k-x_i+x_j}}
{q^{x_i-x_j}-q^{-x_i+x_j}}\, T^{I}\,,
\end{equation}
where $T^I$ stands for $\prod_{i\in I}T^{e_i}$. These operators
coincide (up to a certain gauge) with the quantum integrals of the
trigonometric Ruijsenaars model introduced in \cite{R}.
\end{example}

\subsection{Macdonald polynomials}
The starting point for Macdonald's theory \cite{M1,M2} is that the
operators $D^\pi$, introduced above, preserve the space spanned by
$W$-invariant exponents, or {\it orbitsums} $\mathfrak m_\lambda$:
\begin{equation}
\label{orb} \mathfrak m_\lambda(x) = \sum_{\tau\in W\lambda}
q^{2(\tau,x)}\,,
\end{equation}
where $\lambda\in P_+$ is a dominant weight and the summation is
taken over its $W$-orbit. Moreover, the action of $D^\pi$ is
lower-triangular:
\begin{equation}\label{action}
D^\pi \mathfrak m_\lambda=c_{\lambda \lambda}\mathfrak m_\lambda +
\sum_{\nu \prec \lambda}c_{\lambda\nu} \mathfrak m_\nu
\end{equation}
where the coefficients $c_{\lambda\nu}$ depend on $q,t$ and $\nu
\prec\lambda$ means that $\lambda-\nu$ belongs to $Q_+$.

To introduce the Macdonald polynomials, let us first agree about
terminology. Throughout the paper by a {\it polynomial} in $x$
we will always mean a function $f(x)$ of the form $$f=\sum_{\nu\in
P}f_\nu q^{2(\nu,x)}\,.$$ From algebraic point of view this
corresponds to considering a ring of Laurent polynomials in
$X_i=q^{(\omega_i,x)}$, where $\omega_i$ are the fundamental
weights. As well as a standard polynomial ring, it has a unique
factorization property with $q^{2(\nu,x)}$ being the only
invertible elements.

 Now the Macdonald polynomials $P_\lambda=P_{\lambda}(x;q,t)$
can be defined uniquely as polynomials of the form
\begin{equation}\label{macp}
P_\lambda=\mathfrak m_\lambda+\sum_{\nu \prec
\lambda}a_{\lambda\nu} \mathfrak m_\nu\,,\quad \lambda\in P_+\,,
\end{equation}
which are eigenfunctions of $D^\pi$:
\begin{equation}\label{eigen}
D^\pi P_\lambda=c_{\lambda \lambda} P_\lambda\,.
\end{equation}
The coefficients $a_{\lambda\nu}$ in \eqref{macp} are rational in
$q, t$ and the polynomials $P_{\lambda}$\,($\lambda\in P_+$) have
a number of remarkable properties. In particular, they are
orthogonal with respect to the following scalar product:
\begin{equation}\label{scpr}
\langle f,g \rangle_k={\rm CT}\left(f\overline{\mathstrut g}
\Delta_k \overline{\mathstrut\Delta}_k\,\right)\,,
\end{equation}
where CT means the constant term, $\Delta_k$ is the function
\eqref{delta} and the bar acts on a function $f(x)$ as
$\overline{\mathstrut f}(x)=f(-x)$. This scalar product can be
reinterpreted in terms of a certain integral, which makes perfect
sense for noninteger $k$, too.m  See \cite{M1} for the details.

\bigskip

\begin{remark}
It is not difficult to see that the eigenvalues $c_{\lambda
\lambda}$ for $\lambda\in P_+$ have the form:
\begin{equation}\label{eig}
c_{\lambda\lambda}=\sum_{\tau\in W\pi}
q^{2(\tau,\,\lambda+\rho)}\,,
\end{equation}
where $\rho=\rho_k$ is given by \eqref{rho}.
\end{remark}

\medskip

\begin{remark}
$P_\lambda$ is correctly defined if the diagonal terms in the
action \eqref{action} are distinct:
\begin{equation}\label{nondeg}
c_{\lambda\lambda}\ne c_{\nu\nu}\qquad{\rm for\ all}\quad
\nu\prec\lambda\,.
\end{equation}
This is true for generic $t_\alpha$ and in this case the
polynomials \eqref{macp} are uniquely determined by \eqref{eigen}.
Their coefficients, however, have singularities at certain $q,t$.
For instance, in case $k_\alpha=-m_\alpha$ with $m_\alpha\in
\mathbb Z_+$ which will be central in further considerations, some
first $P_{\lambda}(x;q,q^{-m})$ are not well-defined.
Nevertheless, even in this case the condition \eqref{nondeg} holds
for sufficiently large $\lambda$, i.e. if
$(\lambda,\alpha)>2m_\alpha$ for all $\alpha>0$ (at least,
\eqref{nondeg} will be true for a proper $D^\pi$ or their linear
combination, see \cite{M1,M2}). This means that for such $\lambda$
the whole set of equations \eqref{eigen} for all (quasi)minuscule
coweights together with \eqref{macp} determines $P_\lambda$
correctly.
\end{remark}

\subsection{Macdonald operators in case $k_\alpha \in\mathbb Z_-$}

Let us concentrate now on the case of integer multiplicities
$k_\alpha$. It is known that some results in Macdonald theory are
easier to achieve for integer $k_\alpha$, extending then to
non-integer values by a proper "analytic continuation" in
$k_\alpha$. However, instead of a common assumption $k_\alpha \in
\mathbb Z_+$, we will rather consider the case $k_\alpha
=-m_\alpha$ where $m_\alpha\in \mathbb Z_+$. As we mentioned
above, in this situation the Macdonald's theory is not complete.
Nevertheless, the structure of the eigenfunctions of the operators
$D^\pi$ can be described quite effectively. Next proposition will
be crucial for us.

For given root system $R$ and integer multiplicities
$m=\{m_\alpha\}$ introduce a ring $\mathfrak R$ which consists of
all polynomials $f(x)$ with the following properties: for each
$\alpha \in R_+$ and $j=1,\dots,m_\alpha$
\begin{equation}
\label{axx} f\left(x+ \frac 12 j\alpha^\vee\right) \equiv
f\left(x-\frac 12 j\alpha^\vee\right)\qquad\text{for }
q^{2(\alpha,x)}=1\,.
\end{equation}
\bigskip
\begin{prop}\label{inv}
Let $D^\pi$ be a  Macdonald operator , defined in accordance with
the formulas \eqref{m}--\eqref{m2}. Suppose that all $t_\alpha $
have the form $t_\alpha=q^{-m_\alpha}$ with $ m_\alpha\in\mathbb
Z_+$. Then the operator $D^\pi$ preserves the ring \eqref{axx}: \
$D^\pi(\mathfrak R) \subseteq \mathfrak R$.
\end{prop}

\bigskip
To prove this, we will look first at the rank-one case, $R\subset
V=\mathbb R^1$. In this case we will denote by $T$ the unit shift:
$(Tf)(x)=f(x+1)$. Similarly, $T^s$ will stand for the shift by a
scalar $s$. Let us consider a difference operator $L$ of the form
\begin{equation}\label{L}
L=a(x)T+b(x)T^{-1}\,.
\end{equation}
Suppose that its coefficients $a, b$ are meromorphic with simple
poles at $x=0$ and with no other poles at $x\in \mathbb Z$.
Further, let us fix an integer $m\in \mathbb Z_+$ and impose the
following conditions on $a, b$:
\begin{align}\label{c1}
&{\rm res}_{x=0}(a+b)=0\,,\\\label{c2}&a(j)=b(-j)\quad\text{for}\
j=\pm 1,\dots ,\pm m\,,\\\label{c3}&a(m)=0\,.
\end{align}
Introduce a ring $\mathfrak R_0$ which consists of all meromorphic
functions $f(x)$ with no poles at $x\in\mathbb Z$ and such that
\begin{equation}\label{ri1}
  f(j)=f(-j)\quad\text{for all}\ j=1,\dots,m\,.
\end{equation}

\medskip
\begin{lemma}\label{pres1} Under conditions \eqref{c1}--\eqref{c3}, the operator $L$
preserves the ring \eqref{ri1}: $$L(\mathfrak R_0)\subseteq
\mathfrak R_0\,.$$
\end{lemma}

\medskip
\begin{proof} First of all, for any $f\in\mathfrak R_0$ its
image $Lf$ will be nonsingular at $x\in\mathbb Z$. Indeed, the
only apparent pole is $x=0$. However, it disappears since the
residues of $a$ and $b$ are opposite and $(T-T^{-1})f|_{x=0}=0$
due to \eqref{ri1} at $j=1$.

Now let us prove that $Lf$ still satisfies the conditions
\eqref{ri1}, i.e. that $(T^j-T^{-j})Lf$ is zero at $x=0$. A simple
calculation gives us that
\begin{multline*}(T^j-T^{-j})Lf|_{x=0}= a(j)f(j+1)-b(-j)f(-j-1)+
b(j)f(j-1)-a(-j)f(-j+1)\\=
a(j)\left(f(j+1)-f(-j-1)\right)+b(j)\left(f(j-1)-f(-j+1)\right)
\end{multline*}
(here we used the conditions \eqref{c2}). The last expression must
be zero for all $j=1,\dots,m$ due to \eqref{c2} and \eqref{c3}.
\end{proof}

\medskip
Our next lemma is a modification of the previous one for a
three-term difference operator
\begin{equation}\label{L1}
L=a(x)\left( T^2-1 \right)+ b(x)\left( T^{-2}-1 \right)\,.
\end{equation}
Now $a$ and $b$ may have simple poles at $x=0,\,-1$ and at
$x=0,\,1$, respectively, with no other poles allowed at $x\in
\mathbb Z$. Further, we fix an integer $m$ as before and impose
the following conditions on $a, b$ (in case $m>1$):
\begin{gather*}
{\rm res}_{x=0}(a+b)=0\,,\quad {\rm res}_{x=-1}(a)+{\rm
res}_{x=1}(b)=0\,,\\a(j)=b(-j)\quad\text{for}\ j=1,\pm 2,\dots
,\pm m\,,\\a(m)=a(m-1)=0\,.
\end{gather*}
In case $m=1$ these conditions must be replaced by the following:
\begin{gather*}
{\rm res}_{x=0}(a)={\rm res}_{x=0}(b)=0\,,\quad {\rm
res}_{x=-1}(a)+{\rm res}_{x=1}(b)=0\,,\\a(1)=b(-1)=0\,.
\end{gather*}

\medskip
\begin{lemma}\label{pres2} Under the assumptions above, the operator \eqref{L1}
preserves the ring \eqref{ri1}: $L(\mathfrak R_0)\subseteq
\mathfrak R_0\,.$
\end{lemma}

\medskip
\begin{proof} Take any $f\in\mathfrak R_0$, then the only
possible poles of $Lf$ are $x=\pm 1$ and $x=0$ (if $m>1$). We have
: $$ {\rm res}_{x=0}(Lf) = {\rm
res}_{x=0}(a)\left[f(2)-f(0)\right]+{\rm
res}_{x=0}(b)\left[f(-2)-f(0)\right]\,.$$ This is zero for $m>1$
because $f(2)=f(-2)$ and ${\rm res}_{x=0}(a+b)=0$. Further,
\begin{gather*}
{\rm res}_{x=-1}(Lf)={\rm res}_{x=-1}(a)
\left[f(1)-f(-1)\right]=0\,,\\{\rm res}_{x=1}(Lf)={\rm
res}_{x=1}(b) \left[f(-1)-f(1)\right]=0\,.
\end{gather*}
So, $Lf$ has no singularities at $x\in\mathbb Z$.

Now let us check that $Lf$ still belongs to $\mathfrak R_0$, i.e.
$\left.\left(T^j-T^{-j}\right)Lf\right|_{x=0}=0$ for $j=1,\dots
,m$. A straightforward calculation gives us the following:
\begin{multline*}
\left(T^j-T^{-j}\right)Lf \\= a(x+j)\left[f(x+j+2)-f(x+j)\right]
+b(x+j)\left[f(x+j-2)-f(x+j)\right]\\
-a(x-j)\left[f(x-j+2)-f(x-j)\right]
-b(x-j)\left[f(x-j-2)-f(x-j)\right]\,.
\end{multline*}
For $j>1$ each term in this expression can be evaluated at $x=0$:
\begin{multline*}
\left.\left(T^j-T^{-j}\right)Lf
\right|_{x=0}=a(j)\left[f(j+2)-f(j)\right]
+b(j)\left[f(j-2)-f(j)\right]\\ -a(-j)\left[f(-j+2)-f(-j)\right]
-b(-j)\left[f(-j-2)-f(-j)\right]\\=a(j)\left[f(j+2)-f(-j-2)-f(j)+f(-j)\right]
\\+b(j)\left[f(j-2)-f(-j+2)-f(j)+f(-j)\right]\,.
\end{multline*}
Here we used the conditions $a(j)=b(-j)$ and $a(-j)=b(j)$.  The
resulting expression for all $j=1,\dots ,m$ will be zero due to
the properties \eqref{ri1} of $f$ and the condition
$a(m)=a(m-1)=0$.

Finally, for $j=1$ we have:
\begin{multline*}
\left(T-T^{-1}\right)Lf = a(x+1)\left[f(x+3)-f(x+1)\right]\\
-b(x-1)\left[f(x-3)-f(x-1)\right]\\+\left(a(x-1)+b(x+1)\right)\left[f(x-1)-f(x+1)\right]\,.
\end{multline*}
The last term is zero at $x=0$ because ${\rm res}_{x=0}
\left(a(x-1)+b(x+1)\right)$ is zero and $f(1)=f(-1)$. The first
two terms after evaluating at $x=0$ give
$$a(1)\left[f(3)-f(1)\right]-b(-1)\left[f(-3)-f(-1)\right]\,.$$
This is zero for $m=1,2$ since $a(1)=b(-1)=0$ in this case, and
for $m>2$ since $a(1)=b(-1)$ and $f(3)-f(-3)=f(1)-f(-1)=0$.
\end{proof}

\bigskip
The following two lemmas are a direct corollary of the previous
two.
\medskip
\begin{lemma}\label{sym}
Suppose that a difference operator $L$ of the form \eqref{L} is
invariant under the change of variable $x\to -x$, i.e.
$b(x)=a(-x)$. Further, let us suppose that

\noindent (1) $a,\,b$ have simple poles at $x=0$ and no other
poles at $x\in \mathbb Z$\,,

\noindent (2) $a(m)=0$\,.

\noindent Then $L(\mathfrak R_0)\subseteq \mathfrak R_0$.
\end{lemma}

\medskip
\begin{lemma}\label{sym1}
Suppose that a difference operator $L$ of the form \eqref{L1} is
invariant under the change of variable $x\to -x$, i.e.
$b(x)=a(-x)$. Further, let us suppose that

\noindent (1) $a,\,b$ have simple poles at $x=0,-1$ and $x=0,1$,
respectively, and no other poles at $x\in \mathbb Z$\,,

\noindent (2) $a(m)=a(m-1)=0$\,. \\(In case $m=1$ we replace it by
the condition $a(1)=0$ but now allow poles at $x=\pm 1$ only.)

\noindent
Then $L(\mathfrak R_0)\subseteq \mathfrak R_0$.
\end{lemma}

One can formulate the following inversions of lemmas \ref{pres1}
and \ref{pres2}.

\begin{lemma}\label{inverse}
Suppose we are in a situation described in lemma \ref{pres1}.
Moreover, let us impose extra conditions on the coefficients $a,\,
b$ of the operator \eqref{L} as follows:
\begin{align*}
&{\rm res}_{x=0}(a)\neq 0\,,\\ &a(j)\neq 0 \quad\text{for}\ \,
j=1,\dots,m-1\quad (\text{only in case }m>1)\,.
\end{align*}
Let $f$ be an analytic eigenfunction for $L$: $Lf=\lambda f$. Then
$f$ belongs to the ring $\mathfrak R_0$.
\end{lemma}

\medskip
\begin{lemma}\label{inverse1}
Suppose we are in a situation described in lemma \ref{pres2}. Let
us impose extra conditions on the coefficients $a,\, b$ of the
operator \eqref{L1} as follows:
\begin{align*}
&{\rm res}_{x=0}(a)\neq 0 \quad(\text{only in case }m>1)\,,\\
&{\rm res}_{x=-1}(a)\neq 0\,,\\ &a(j)\neq 0 \quad\text{for}\ \,
j=1,\dots,m-2\quad (\text{only in case }m>2)\,.
\end{align*}
Then each analytic eigenfunction $f$ of the operator $L$ must
belong to the ring $\mathfrak R_0$.
\end{lemma}

\medskip
\begin{proof}Both lemmas can be proven by reversing the arguments used to
prove lemmas \ref{pres1}, \ref{pres2}. Indeed, since $f$ is
non-singular, $Lf=\lambda f$ must be non-singular, too. Looking at
the residues of $Lf$, we obtain the condition $f(1)=f(-1)$. Thus,
$Lf=\lambda f$ must also satisfy this condition, which gives more
conditions on $f$, and so on.

\end{proof}

\bigskip
Let us apply all this to the difference operators
\begin{equation}\label{rank1}
D_1=\frac{q^{x-m}-q^{-x+m}}{q^x-q^{-x}}T+\frac{q^{x+m}-q^{-x-m}}{q^x-q^{-x}}T^{-1}
\end{equation}
and
\begin{multline}\label{rank11}
D_2=\frac{(q^{x-m}-q^{-x+m})(q^{x-m+1}-q^{-x+m-1})}{(q^x-q^{-x})(q^{x+1}-q^{-x-1})}
\left(T^2-1\right)\\+
\frac{(q^{-x-m}-q^{x+m})(q^{-x-m+1}-q^{x+m-1})}{(q^{-x}-q^{x})(q^{-x+1}-q^{x-1})}
\left(T^{-2}-1\right)\,.
\end{multline}
which are the Macdonald operators $D^\pi$ in cases $R=A_1=\{\pm
1\}\subset \mathbb R^1$ and $\pi=1$ and $2$, respectively. It is
obvious that $D_1$ and $D_2$ satisfy the conditions of lemma
\ref{sym} and \ref{sym1} (provided that $q$ is not a root of
unity). Hence, they preserve the properties \eqref{ri1}. Moreover,
instead of $x=0$ we may consider any point $x=x_0$ with
$q^{x_0}-q^{-x_0}=0$. Indeed, for such $x_0$ the symmetry $x\to
2x_0-x$ still does not change the operators $D_1$ and $D_2$ (this
reflects their invariance with respect to the corresponding affine
Weyl group). So, after shifting the origin to $x_0$ one gets the
operators with the same properties as in lemma \ref{sym} and
\ref{sym1}. Thus, the operators $D_1,\,D_2$ preserve, in fact, a
bigger ring, which is an affine version of the ring \eqref{ri1}.

Namely, let us consider the ring of all analytic functions $f(x)$
with the following properties:
\begin{equation}\label{axx1}
 f(x+j)=f(x-j)\qquad\text{for each } j=1,\dots,m \quad\text{and }
q^{2x}=1\,.
\end{equation}
The arguments above prove the following result.

\medskip
\begin{prop}\label{per}
The operators \eqref{rank1}, \eqref{rank11} with $m\in \mathbb
Z_+$ preserve the ring $\mathfrak R$ of analytic functions with
properties \eqref{axx1}.
\end{prop}

\bigskip \noindent This is, essentially, the rank-one case of
Proposition \ref{inv}.
Now we are ready to prove it in full generality.

\medskip \noindent \begin{proof}[Proof of Proposition \ref{inv}] We will
only consider the operator $D^\pi$ given by the formulas
\eqref{m1}--\eqref{m2}, since this covers the case \eqref{m}, see
Remark \ref{cover}. Choose any $\alpha\in R$, then we should prove
that $D^\pi$ preserves the ring $\mathfrak R_{\alpha}$ of
functions with the following properties:
\begin{equation}
\label{axxa} f\left(x+ \frac 12 j\alpha^\vee\right) =
f\left(x-\frac 12 j\alpha^\vee\right)\quad\text{for }
j=1,\dots,m_\alpha \quad\text{and } q^{2(\alpha,x)}=1\,.
\end{equation}

Let $s=s_\alpha\in W$ be the reflection with respect to $\alpha$.
Since the orbit $W\pi$ of the coweight $\pi$ is $s$-invariant, it
splits into pairs $\tau,\,\tau'$ with $\tau'=s(\tau)$ plus a
number of $s$-invariant $\tau$'s. This defines a splitting of
$D^\pi$ into a sum of difference operators of one of the following
three types:
\begin{gather}
D_0=a_\tau(x)(T^\tau-1)\,,\quad s(\tau)=\tau\,,\\
D_1=a_\tau(x)(T^\tau-1)+a_{\tau'}(x)(T^{\tau'}-1)\,,\quad
\tau'=\tau-\alpha^\vee\,,\\ D_2=a_{\alpha^\vee}(T^{\alpha^\vee}-1)
+a_{-\alpha^\vee}(T^{-\alpha^\vee}-1)\qquad (\text{ only for
quasiminuscule }\pi\ )\,.
\end{gather}
This follows directly from the fact that $\pi$ is
(quasi)minuscule. Moreover, since $D^\pi$ was obviously
$W$-invariant, each of $D_i$ will be invariant under the
reflection $s$. In particular, $a_\tau(x)=a_\tau(sx)$ in $D_0$,
$a_{\tau'}(x)=a_\tau(sx)$ in $D_1$, and
$a_{-\alpha^\vee}(x)=a_{\alpha^\vee}(sx)$ in $D_2$.

We claim that each of $D_i$ preserves the ring \eqref{axxa}.
First, note that in case of $D_0$ $a_\tau$ has no pole at
$q^{2(\alpha,x)}=1$ since it is $s$-invariant. So, $a_\tau$ itself
belongs to the ring $\mathfrak R_\alpha$. Also it is clear that
operator $T^{\tau}-1$ preserves this ring (the shift acts in
direction, orthogonal to $\alpha$). Hence, $D_0(\mathfrak
R_\alpha)\subset \mathfrak R_\alpha$.

Now let us consider $D_1$, it has the form
$$a_\tau(x)T^\tau+a_{\tau'}(x)T^{\tau'}-a_\tau(x)-a_{\tau'}(x)\,.$$
The sum $a_\tau+a_{\tau'}$ is $s$-invariant, hence, it is
nonsingular at $q^{2(\alpha,x)}=1$ and belongs to $\mathfrak
R_\alpha$. So, we may ignore it and consider the first two terms
only. Further, note that we can present $\tau,\,\tau'$ as
$$\tau=\frac12\alpha^\vee+v\,,\quad \tau'=-\frac12\alpha^\vee+v$$
for a certain $v$ such that $(\alpha,v)=0$. Hence, $T^v$ will
preserve the ring $\mathfrak R_\alpha$ and we can also ignore it,
reducing $D_1$ to
\begin{equation}\label{red}
a_\tau T^{\frac12\alpha^\vee}+a_{\tau'}
T^{-\frac12\alpha^\vee}\,.
\end{equation}
This operator is still $s$-invariant. Moreover, since $(\alpha,
\tau)=1$ in this case, we see from the formula \eqref{m2} for
$a_\tau$ that
$$a_\tau(x+\frac12m_\alpha\alpha^\vee)=0\qquad\text{for
}q^{2(\alpha,x)}=1\,.$$ Now in absolutely the same way as it was
for the operator \eqref{rank1}, we conclude that \eqref{red}
preserves the ring $\mathfrak R_\alpha$.

In the same manner the case of $D_2$ reduces to \eqref{rank11}.
\end{proof}

\bigskip\noindent

\medskip\noindent
\begin{remark}
One can show that the ring $\mathfrak R$ is finitely generated,
therefore it can be viewed as the coordinate ring of a certain quite
specific affine algebraic variety. For instance, for $A_1$ case it
is a rational curve with $m$ double points.
\end{remark}

\noindent\begin{remark} In our proof of Proposition \ref{inv}
essential ingredients were the $W$-invariance of the operator
$D^\pi$ and specific location of poles and zeros of its
coefficients. This has certain parallels with the {\it residue
construction} of Hecke algebras from \cite{GKV}, \cite{BEG}.
Moreover, using the results from \cite{BEG}, one can prove that
for $k_\alpha=-m_\alpha\in\mathbb Z_-$ all {\it
Macdonald--Cherednik operators} $D_1,\dots, D_n$ coming from
$W$-invariant part of the double affine Hecke algebra \cite{Ch1}
will preserve the ring $\mathfrak R$. All these operators will
commute with $D^\pi$. However, as we will see later, the
centralizer of $D^\pi$ in case $k_\alpha\in\mathbb Z_-$ is much
bigger, and it contains many {\it non-symmetric} difference
operators, preserving $\mathfrak R$.
\end{remark}

\bigskip
One can formulate an analogue of Proposition \ref{inv} for the case
of positive $k_\alpha\in \mathbb Z_+$, too. This is because these
two cases are related through a simple gauge transformation. Let
$\pi$ be a (quasi)minuscule coweight and $D_{m}$ denote the
corresponding Macdonald operator \eqref{m2} with $t=q^{-m}$. For
$m_\alpha\in\mathbb Z_+$ introduce a function $\delta_m(x)$ as
follows:
\begin{equation}\label{deltax}
 \delta_m(x)=\prod_{\alpha \in
R_+} \prod_{j=-m_\alpha}^{m_\alpha} [(\alpha,x)+j]\,,
\end{equation}
where $[a]$ denotes $[a]=q^a-q^{-a}$.
\medskip
The following fact is well-known and can be checked by a direct calculation.

\begin{lemma}\label{conj} For $m_\alpha\in \mathbb Z_+$ we have the
following relation between the Macdonald operators $D_{m}$ and
$D_{m'}$ with $m'=-1-m$: $$\delta_m^{-1}\circ D_{m}\circ
\delta_m=D_{m'}\,.$$
\end{lemma}

\medskip
\begin{cor}
Let $\mathfrak R=\mathfrak R_m$ be the ring \eqref{axx}. Then the
Macdonald operator $D$, given by \eqref{m1}--\eqref{m2} with
$t=q^{m+1}$, preserves the $\mathfrak R$-module
$U=\delta_m^{-1}\mathfrak R$:\quad $D(U) \subseteq U$.
\end{cor}

\section{Baker--Akhiezer functions for Macdonald operators}
\label{baf}

We keep mostly the notations of the previous section. So, we
consider an arbitrary (reduced, irreducible) root system $R$ and
fix a $W$-invariant set $m$ of {\it multiplicities} $m_\alpha\in
\mathbb Z_+$. Our purpose is to construct eigenfunctions of the
corresponding Macdonald operators \eqref{m1} with
$t_\alpha=q^{-m_\alpha}$ (so $k_\alpha=-m_\alpha$ in notations of
section \ref{beg}) . Keeping this in mind, we introduce
$\rho=\rho(m)$ instead of $\rho=\rho_k=-\rho_m$ from \eqref{rho}:
\begin{equation}
\label{rho1} \rho = \rho_m =\frac 12 \sum_{\alpha \in
R_+}m_{\alpha} \alpha\,.
\end{equation}
We will also use its counterpart for the dual root system
$R^\vee$:
\begin{equation}
\label{rho2}
\rho^\vee = \frac 12 \sum_{\alpha \in R_+}m_{\alpha}
\alpha^\vee\,.
\end{equation}
( Warning: $\rho^\vee\ne 2\rho/(\rho,\rho)$\,!)

In this section we often will deal with functions of two variables
$x,z\in V$. We will keep calling a sum $f(x)=\sum_{\nu\in
P}q^{2(\nu,x)}$ a {\it polynomial} in $x$. However, switching to
$z$, we will also switch from the root system $R$ to its dual
$R^\vee$. For instance, by a polynomial in $z$ we will mean a sum
$f(z)=\sum_{\nu\in P^\vee}q^{2(\nu,z)}$. Such (perhaps confusing)
terminology is caused by our implicit identification of the vector
space $V$ and its dual $V^*$. To distinguish between these spaces,
it would be natural to assume that $R\subset V^*$ and
$R^\vee\subset V$. In this case $x$ and $z$ would lie in $V$ and
$V^*$, respectively. However, we prefer not to do this, keeping
things simple.

We will apply the term {\it quasipolynomial} (in $x$ or in $z$) to
a function of the form $q^{2(x,z)}f$, where $f$ is polynomial in
$x$ or in $z$, respectively.

For a polynomial $f(x)=\sum_{\nu\in P} a_\nu q^{2(\nu,x)}$ by its
{\bf support} \ ${\rm supp}(f)$\  we will always mean the {\bf
convex hull} of all points $\nu$ with $a_\nu\ne 0$. In the same
way we define the support of $g(z)=\sum_{\nu\in P^\vee} a_\nu
q^{2(\nu,z)}$. For a quasipolynomial in $x$ of the form
$\phi=q^{2(x,z)}f(x)$ by its support ${\rm supp}(\phi)$ we will
simply mean the support of $f(x)$. In case when $\phi$ is
quasipolinomial in $x$ and $z$ at the same time, it usually will
be clear which support we are considering (either in $x$ or in
$z$).

\subsection{Baker--Akhiezer function: definition and uniqueness}

Let $\psi(x,z)$ be a function of two variables $x,z\in V$ of the
form
\begin{equation}\label{psi}
\psi = q^{2(x,z)}\sum_{\nu\in\mathcal N} \psi_\nu q^{2(\nu,z)}\,,
\end{equation}
where the coefficients $\psi_\nu=\psi_\nu(x)$ depend on $x$,
$(x,z)$ is the scalar product in $V$ and the summation  in
\eqref{psi} is taken over all coweights $\nu\in P^\vee$ from the
following polytope $\mathcal N$:
\begin{equation}
\label{nu} \mathcal N = \{ \nu= \rho^\vee-\sum_{\alpha \in R_+}
l_\alpha \alpha^\vee \mid 0\le l_\alpha \le m_\alpha \}\,.
\end{equation}
Using our conventions about terminology, these conditions on a
function $\psi$ can be rephrased as follows: $\psi$ is
quasipolynomial in $z$ with ${\rm supp}(\psi)\subseteq \mathcal
N$.

Suppose that $\psi$ satisfies also the following conditions: for
each $\alpha\in R$ and $s=1,\dots ,m_\alpha$
\begin{equation}\label{axzpsi}
\psi\left(x,\, z+{\frac 12 s\alpha}\right)\equiv \psi\left(x,
\,z{-\frac 12 s\alpha}\right)\quad\text{for}\ \
q^{2(\alpha^\vee,z)}=1\,.
\end{equation}

\medskip
\begin{definition}
A function $\psi(x,z)$ with the properties
\eqref{psi}--\eqref{axzpsi} is called a {\bf Baker--Akhiezer (BA)
function} associated to the data $\{R,m\}$.
\end{definition}

\bigskip
Our terminology is motivated by the fact that in case $R=A_1$  such a $\psi$
is a Krichever's Baker--Akhiezer function
\cite{Kr,Kr1} associated to a specific singular rational curve. In
contrast with the one-dimensional case, in higher dimension the
main problem is to prove the existence of such a function. We do
this in the next subsection. Let us presume now that such a $\psi$
{\it does exist}.

\medskip
\begin{prop}\label{prop1}
Properties \eqref{psi}-- \eqref{axzpsi} determine $\psi$ uniquely
up to a factor  depending on $x$.
\end{prop}

Proof is based on the following two lemmas.

\begin{lemma}\label{lemma1}
Let a quasipolynomial in $z$\  $\psi(x,z)= q^{2(x,z)}\sum_{\nu\in
P^\vee}\psi_\nu q^{2(\nu,z)}$ satisfy the conditions
\eqref{axzpsi}. Then for each $\alpha^\vee \in R^\vee$ and for any
$\nu\in P^\vee$ the set of integers $j$ such that
$\psi_{\nu+j\alpha^\vee} \ne 0$ either is empty or contains at
least two integers $j_1, j_2$ with $|j_1-j_2|\ge m_\alpha$.
\end{lemma}

\medskip
\begin{lemma}\label{lemma2}
Let $l_1,\dots, l_r$ be a set of non-parallel segments in affine
Euclidean space $V$ and $\Omega \subset V$ be a convex domain in
$V$. Suppose that for each $l_i$ and for any line $l$, which
intersects $\Omega$ and is parallel to $l_i$, the intersection
$\Omega \cap l$ has the length greater or equal than $|l_i|$. Then
$\Omega$ can be presented as $\Omega=\Omega '\# \mathcal N$ for
some convex domain $\Omega '$,where $\mathcal N=l_1\#l_2\#\dots
\#l_r$. Here $\#$ denotes the Minkowski addition in $V$.
\end{lemma}

\medskip

We recall that the Minkowski sum of two subsets $A,\,B$ of an
affine space $V$ is formed by all the points $a+b$, where $a$ and
$b$ run over $A$ and $B$, respectively. The addition of points, of
course, depends on the choice of origin, but the resulting set
will be the same up to a shift. This operation is relevant to the
multiplication of polynomials: if $f_1, f_2$ are two polynomials
in $x$, and $N_i={\rm supp}(f_i)$ then $N={\rm supp}(f_1f_2)$ is
the Minkowski sum of $N_1$ and $N_2$.

\bigskip

\noindent \begin{proof}[Proof of Lemma \ref{lemma1}] For a given
$\alpha\in R_+$, substitution of the $\psi$ into \eqref{axzpsi}
gives the following set of relations: $$ \sum_{\nu\in P^\vee}
\psi_\nu q^{2(\nu,z)}
\left(q^{(s\alpha,x+\nu)}-q^{-(s\alpha,x+\nu)}\right)=0
\quad\text{for}\quad q^{2(\alpha^\vee,z)}=1 $$
$(s=1,\dots,m_\alpha)$.

These relations split up into separate linear equations for each
"$\alpha$-string" $\nu_j=\nu_0+j\alpha^\vee$\ ($j\in\mathbb Z$):
\begin{equation}\label{lin}
 \sum_j \psi_j\left( (q_j)^s-(q_j)^{-s}\right)=0 \qquad (
s=1\dots m_\alpha )\,,
\end{equation}
 where
$\psi_j := \psi_{\nu_0+j\alpha^\vee}$ and $ q_j :=
q^{2j}q^{(\alpha, x+\nu_0)}$.

Suppose now that among the coefficients $\psi_j$ only
$\psi_1,\dots,\psi_{m_\alpha}$ do not vanish. In this situation we
would have a homogeneous linear system of $m_\alpha$ equations for
$m_\alpha$ unknowns $\psi_j$. Thus, it would be sufficient to show
that this system is non-degenerate for generic $x$. To see this,
we can look at the asymptotic behaviour of its determinant at
large $x$ (cf. \cite{ES}). More precisely, we consider the
corresponding matrix $A=(a_{ij})_{i,j=1,\dots,m_\alpha}$ with
$a_{ij}=(q_j)^i-(q_j)^{-i}$ where $q_j := q^{2j}q^{(\alpha,
x+\nu_0)}$. Then for large $x$ such that $q^{(\alpha, x+\nu_0)}\gg
1$ the determinant $\det A$ asymptotically equals the Vandermonde
determinant $\det (q_j^i)$ which is nonzero since $q$ is not a
root of unity.

 Thus, the system  is non-degenerate for generic $x$ and all
$\psi_j$ must vanish, which proves the lemma.
\end{proof}

\bigskip
\noindent \begin{proof}[Proof of Lemma \ref{lemma2}] We will prove
the lemma by induction in the number $s$ of the segments. For an
easier reference to the assumptions from the lemma, let us say in
such a situation that $\Omega$ {\bf dominates} over the segments
$l_1,\dots, l_s$.

First, suppose we have just one segment $l_1$ and a convex domain
$\Omega$ which dominates over $l_1$. Then $\Omega=l_1\# \Omega'$,
where $\Omega'$ is the intersection $\Omega\cap T(\Omega)$ of
$\Omega$ and its image under the shift $T$ for a vector
$\overrightarrow {l_1}$, associated with the segment $l_1$ (in
either of two possible directions). This proves the lemma in case
$s=1$.

Now suppose that $\Omega$ dominates over $l_1,\dots, l_s$. Take
the first segment $l_1$ and consider the convex domain $\Omega'$
constructed above, so we have $\Omega=l_1\# \Omega'$. We claim
that $\Omega'$ still dominates over $l_2,\dots, l_s$.

To prove this, let us take, for instance, $l_2$ and choose any
line $l$, parallel to it. Now we may consider a two-dimensional
section of the $\Omega$ passing through $l$ and parallel to $l_1$.
The resulting two-dimensional domain will, obviously, dominate
over $l_1$ and $l_2$. Thus, essentially, we need to check the
lemma in dimension two, for $s=2$. This is very simple, and we
leave it to the reader.

So, we have proved that $\Omega'$ dominates over $l_2,\dots, l_s$.
Now the statement of the lemma follows by an obvious induction.
\end{proof}

\bigskip
\noindent \begin{proof}[Proof of Proposition \ref{prop1}] First,
notice that the polytope $\mathcal N$ in \eqref{nu} is exactly the
Minkowski sum of the segments $m_\alpha\alpha^\vee$ with
$\alpha^\vee\in R_+^\vee$ (abusing notation, we denote a
vector and associated segment by the same symbol). Now let $\psi$
be any BA function. By definition, we have an inclusion $${\rm
supp}(\psi)\subseteq \mathcal N\,.$$ On the other hand, Lemma
\ref{lemma1} implies that the polytope ${\rm supp}(\psi)$ must
dominate over each of the segments
$m_\alpha\alpha^\vee$,\,$\alpha^\vee\in R_+^\vee$. Hence, by Lemma
\ref{lemma2} it must contain (a copy of) the polytope \eqref{nu}
which is their Minkowski sum. Altogether this proves that for each
(nonzero) BA function $\psi$ one has the equality $${\rm
supp}(\psi)= \mathcal N\,.$$

Let now $\psi'$, $\psi''$ be two Baker--Akhiezer functions.
Consider their linear combination $\psi=\psi'-c(x)\psi''$, which
still satisfies the conditions \eqref{psi}--\eqref{axzpsi}. We can
choose $c(x)$ in such way that the resulting function $\psi$ will
have zero coefficient at one of the vertices of the ${\rm
supp}(\psi')={\rm supp}(\psi'')=\mathcal N$. So, we will have a
strict inclusion ${\rm supp}(\psi)\subset \mathcal N$. Thus, the
only possibility is that such $\psi$ is zero, hence,
$\psi'=c(x)\psi''$.
\end{proof}

\medskip
\begin{cor}\label{only}
For a BA function $\psi$, the nonzero coefficients $\psi_\nu$ in
\eqref{psi} can appear only for $\nu=\rho^\vee-\sum_{\alpha \in
R_+} l_\alpha \alpha^\vee$ with {\bf integer} $l_\alpha$. In other
words, the summation in \eqref{psi} is taken effectively only over
the set $\rho^\vee+Q^\vee\subset P^\vee$.
\end{cor}

\medskip
\begin{proof}
Suppose we have other terms, then let us remove them from the sum
\eqref{psi}. This would not affect the conditions \eqref{axzpsi}.
Indeed, in the process of proving lemma \ref{lemma1} we saw that
these conditions split into separate linear equations involving
$\nu$'s from the same coset in $P^\vee/Q^\vee$. But from the
uniqueness of $\psi$ it follows that the resulting function must
remain the same. Hence, there were no other terms at all.
\end{proof}

\subsection{Existence of BA function}

Comparison of Proposition \ref{inv} and conditions \eqref{axzpsi}
suggests the idea to use a Macdonald operator acting in the
$z$-variable in order to construct a Baker--Akhiezer function
$\psi$. Let $\omega\in P$ be a (quasi)minuscule {\it weight} for
the system $R$ and $D_z^\omega$ be the Macdonald operator
\eqref{m1} corresponding to the {\it dual system} $R^\vee$ and
acting in the $z$-variable:
\begin{equation}
\label{m1z} D_z^\omega = \sum_{\tau\in W\omega}a_\tau
\left(T_z^{\tau}-1\right) + \sum_{\tau\in W\omega}
q^{-2(\rho^\vee,\tau)}\,,
\end{equation}
where  $\rho^\vee$ is given by \eqref{rho2} and
\begin{equation}
\label{m2z} a_\tau (z) = \prod_{\genfrac{}{}{0pt}{}{\alpha \in R
:}{(\alpha,\tau)>0}} \frac {[(\alpha^\vee,z)-m_\alpha]}
{[(\alpha^\vee,z)]}\prod_{\genfrac{}{}{0pt}{}{\alpha \in R
:}{(\alpha^\vee, \tau)=2}} \frac {[(\alpha^\vee,z)-m_\alpha+1]}
{[(\alpha^\vee,z)+1]} \,,
\end{equation}
where $[a]=q^a-q^{-a}$.

 Introduce the ring $\mathfrak R^\vee$
which is a counterpart of the ring \eqref{axx} and consists of all
polynomials $f(z)=\sum_{\nu\in P^\vee}f_\nu q^{2(\nu,z)}$ with the
following properties: for each $\alpha \in R_+$ and
$j=1,\dots,m_\alpha$
\begin{equation}
\label{axz} f\left(z+ \frac 12 j\alpha\right) \equiv
f\left(z-\frac 12 j\alpha\right)\qquad\text{for }
q^{2(\alpha^\vee,z)}=1\,.
\end{equation}
Then, according to Proposition \ref{inv}, the operator
$D_z^\omega$ will preserve the ring $\mathfrak R^\vee$: \
$D_z^\omega(\mathfrak R^\vee) \subseteq \mathfrak R^\vee$.

Now we need one technical lemma which shows that the action of
$D^\omega_z$ on $\mathfrak R^\vee$ is "lower-triangular".
\bigskip
\begin{lemma}\label{lead} Let $D=D_z^\omega$ be the Macdonald
operator \eqref{m1z}. Suppose that both $f$ and $\widetilde f=Df$
are polynomials in $z$:
\begin{equation}\label{ff}
f(z)=\sum_{\nu\in P^\vee}f_\nu q^{2(\nu,z)}\,,\qquad \widetilde
f(z)=\sum_{\nu\in P^\vee}\widetilde f_\nu q^{2(\nu,z)}\,.
\end{equation}
Then ${\rm supp}(\widetilde f) \subseteq {\rm supp}(f)$. Further,
let $\lambda$ be a vertex of the polytope $N={\rm supp}(f)$, then
the ratio $c_\lambda=\widetilde f_\lambda/f_\lambda$ of the
corresponding coefficients in \eqref{ff} can be calculated as
follows. First, choose generic $v\in V$ such that $(v,\lambda)\ge
(v,\nu)$ for all $\nu\in N$ and put $$\rho^\vee_\lambda=\frac12
\sum_{\alpha\in R\,:\,(\alpha,v)>0}m_\alpha \alpha^\vee\,.$$ Then
one has $c_\lambda=\sum_{\tau\in W\omega}
q^{2(\tau,\lambda-\rho^\vee_\lambda)}$.
\end{lemma}

\medskip
This can be proven similar to the proof of \eqref{action} and
\eqref{eig} in \cite{M1,M2}. A key point is that inclusion
$A\#C\subseteq B\#C$ for convex domains $A,B,C$ implies
$A\subseteq B$.

\medskip
A similar result is true for quasipolynomials. Recall that for a
quasipolynomial $\phi=q^{2(x,z)}f(z)$ its support, by our
conventions, coincides with the support of $f$. Thus, we have the
following analog of the lemma above.
\medskip
\begin{lemma}\label{leading}
Suppose that both $\phi$ and $\widetilde \phi=D\phi$ are
quasipolynomials in $z$:
\begin{equation}\label{fifi}
\phi(z)=q^{2(x,z)}\sum_{\nu\in P^\vee}\phi_\nu
q^{2(\nu,z)}\,,\qquad \widetilde \phi(z)=q^{2(x,z)}\sum_{\nu\in
P^\vee}\widetilde \phi_\nu q^{2(\nu,z)}\,.
\end{equation}
Then ${\rm supp}(\widetilde \phi) \subseteq {\rm supp}(\phi)$.
Further, let $\lambda$ be a vertex of the polytope $N={\rm
supp}(\phi)$, then the ratio $c_\lambda=\widetilde
\phi_\lambda/\phi_\lambda$ of the corresponding coefficients in
\eqref{fifi} can be calculated as
\begin{equation}\label{lea}
c_\lambda=\sum_{\tau\in W\omega}
q^{2(\tau,x+\lambda-\rho^\vee_\lambda)} \,,
\end{equation}
where $\rho^\vee_\lambda$ is defined in the lemma above.
\end{lemma}

\bigskip
Now everything is ready to construct a BA function. The idea is
very simple. We start from the quasipolynomial
$\phi=q^{2(x,z)}Q(z)$, where
\begin{equation}\label{Q}
Q(z)=q^{2(\rho^\vee,z)}\prod_{\alpha\in
R_+}\prod_{j=1}^{m_\alpha}[(\alpha^\vee,z)+j]\,[(\alpha^\vee,z)-j]\,,
\end{equation}
where $[a]$, as usual, denotes $q^a-q^{-a}$. This polynomial is
especially chosen to guarantee that $\phi$ satisfies the
conditions \eqref{axz} in $z$. Thus, applying $D=D^\omega_z$
successively to $\phi$, we will always get a quasipolynomial in
$z$ which will still satisfy the conditions \eqref{axz} in $z$.
Let us apply at each step an operator $D-c_i$, so
$\phi_{i+1}=(D-c_i)\phi_i$,\,$\phi_0=\phi$. The coefficients $c_i$
will be adjusted to reduce ${\rm supp}(\phi_i)$, see below.
Finally, $\psi$ will be obtained after repeating this sufficiently
many times. Before proceeding with more details, let us make one
more remark. Notice that the formula \eqref{Q} implies that the
only nonzero terms in the initial quasipolynomial
$\phi=q^{2(x,z)}\sum_{\nu\in P^\vee} \phi_\nu q^{2(\nu,z)}$ are
those with $\nu\in \rho^\vee + Q^\vee\subset P^\vee$. To see this,
one should rewrite each factor $[(\alpha^\vee,z)+j]$ in \eqref{Q}
as $q^{(\alpha^\vee,z)}(q^j-q^{-j-2(\alpha^\vee,z)})$. The same
will be true for all successive functions $\psi_i$. One can see
this directly after rewriting in a similar way the coefficients
\eqref{m2z} of the difference operator $D$.

Now let us make everything more concrete. At the beginning we have
$$N_0={\rm supp}(\phi_0)=\{\nu=\rho^\vee+\sum_{\alpha\in
R_+}l_\alpha\alpha^\vee\, |\,-m_\alpha\le l_\alpha \le
m_\alpha\}\,.$$ Further, as we mentioned already, each of $\phi_i$
will satisfy the conditions \eqref{axzpsi}. Hence, due to lemmas
\ref{lemma1} and \ref{lemma2}, its support $N_i\subseteq N_0$ must
be a union of several copies of the polytope $\mathcal N$ given by
\eqref{nu}. Let us fix generic $v$ lying inside the positive Weyl
chamber $C$. The linear functional $(v,\cdot)$ determines the {\it
height function} on $V$. We will call a vertex of a polytope the
{\it highest} (respectively, the lowest) vertex if it has maximal
(respectively, minimal) height among all vertices. At each step we
will choose the highest vertex $\lambda$ of the polytope $N_i={\rm
supp}(\phi_i)$. For brevity, let us call the corresponding
coefficient $\phi_\lambda$ the {\it highest coefficient} of
$\phi$. Note that the initial polytope $N_0$ (as well as all
successive $N_i$) will be composed of the images of $\mathcal N$
under some of the shifts by $\nu=\sum_{\alpha\in
R_+}l_\alpha\alpha^\vee$ with $l_\alpha=0,\dots,m_\alpha$. If we
look at the smaller polytope $\mathcal N$, then its highest (resp.
lowest) vertex will be $\lambda=\rho^\vee$ (resp.
$\lambda=-\rho^\vee$). Hence, the highest vertex of the polytope
$N_i$ must be of the form
\begin{equation}\label{pr}
\lambda=\rho^\vee+\nu\,,\qquad \nu=\sum_{\alpha\in
R_+}l_\alpha\alpha^\vee\,,\quad  l_\alpha=0,\dots, m_\alpha\,.
\end{equation}
Now we can kill the highest coefficient $\phi_\lambda$ by applying
$D-c_i$ where $c_i=c_\lambda$ is given by \eqref{lea}. Note that
the vector $\rho^\vee_\lambda$ in the formula \eqref{lea} will be
simply $\rho^\vee$ (because $\lambda$ is the highest vertex).
Thus, $c_i=c_\lambda$ will be $c_i=\sum_{\tau\in W\omega}
q^{2(\tau,x+\nu)}$, or simply $$c_i=\mathfrak m_\omega(x+\nu)$$ in
terms of the orbitsum $\mathfrak m_\omega$,
\begin{equation}\label{mu}
\mathfrak m_\omega(x)=\sum_{\tau\in W\omega} q^{2(\tau,x)}\,.
\end{equation}
On the other hand,  let us look now what is happening with the
lowest coefficient $\phi_{-\rho^\vee}$. According to Lemma
\ref{leading}, after application of $D-c_i$ it gets the factor
$c_{-\rho^\vee}-c_i$. Note that for $\lambda=-\rho^\vee$ the
vector $\rho^\vee_\lambda$ in the formula \eqref{lea} will be
simply $-\rho^\vee$. So, the formula \eqref{lea} gives us:
$$c_{-\rho^\vee}=\sum_{\tau\in W\omega} q^{2(\tau,x)}=\mathfrak
m_\omega(x) \,.$$ Obviously, $c_{-\rho^\vee}-c_i=\mathfrak
m_\omega(x)-\mathfrak m_\omega(x+\nu)$ is nonzero as soon as
$\nu\ne 0$ (recall that we assume that $q$ is not a root of
unity).

These considerations imply the existence of BA functions for all
root system, which is one of our main results.

\bigskip
\begin{theorem}\label{exist}
Let $D=D^\omega_z$ denote the Macdonald operator \eqref{m1z}.
Define $\psi(x,z)$ as follows:
\begin{equation}\label{ber}
\psi=\prod_{\nu}\left(D-\mathfrak
m_\omega(x+\nu)\right)\left[q^{2(x,z)}Q(z) \right]\,,
\end{equation}
in accordance with the formulas \eqref{Q}, \eqref{mu}, where the
product is taken over all $\nu\ne 0$ having the form
$\nu=\sum_{\alpha\in R_+}l_\alpha\alpha^\vee$ with
$l_\alpha=0,\dots , m_\alpha$. Then

\noindent(i) $\psi$ has the form \eqref{psi}--\eqref{nu};

\noindent(ii) the coefficient $\psi_{-\rho^\vee}$ in its expansion
\eqref{psi} equals
\begin{equation}\label{bernorm}
\psi_{-\rho^\vee}=\prod_\nu (\mathfrak m_\omega(x)-\mathfrak
m_\omega(x+\nu))\ne 0\,;
\end{equation}

\noindent(iii) $\psi$ is a Baker--Akhiezer function for the system
$R$ with multiplicities $m=\{m_\alpha\}$;

\noindent(iv) as a function of $z$,\ $\psi$ is an eigenfunction of
the Macdonald operator $D$: $$D\psi=\mathfrak m_\omega(x)\psi\,.$$
\end{theorem}

\medskip
\begin{proof}
As we explained above, the constructed function \eqref{ber} will
be a quasipolynomial in $z$ satisfying the conditions
\eqref{axzpsi}. The part (ii) follows immediately from the
construction of $\psi$. It implies that the polytope ${\rm
supp}(\psi)$ contains the polytope \eqref{nu}. On the other hand,
the arguments above show that the highest vertex of ${\rm
supp}\psi$ must be $\lambda=\rho^\vee$, because every "higher"
term $\psi_{\rho^\vee+\nu}$ has been killed after applying
$D-\mathfrak m_\omega(x+\nu)$. Altogether, this gives us that
${\rm supp}(\psi)$ coincides with the polytope $\mathcal N$ in
\eqref{nu}. Thus, the part (i) is also proven.

Part (iii) of the theorem follows from the previous two and the
remark above that $\psi$ satisfies conditions \eqref{axzpsi}.
Finally, $\widetilde \psi=(D-\mathfrak m_\omega(x))\psi$ must be
quasipolynomial in $z$ with all the properties
\eqref{psi}--\eqref{axzpsi}. Hence, it must be proportional to
$\psi$ due to Proposition \ref{prop1}. However, application of
$D-\mathfrak m_\omega(x)$ kills the highest coefficient
$\psi_{\rho^\vee}$, so $\widetilde\psi$ must be zero, which proves
the last part.
\end{proof}

\medskip
\begin{remark}
Formula \eqref{ber} is a trigonometric version of a related
formula from \cite{C00}, which, in its turn, is a discrete version
of the Berest's formula \cite{Be}.
\end{remark}

\section{Bispectral duality}\label{duu}

In this section we will explain how one should normalize BA
function to achieve a certain symmetry between $x$ and $z$
variables. We start  by looking closely at the rank one case.

\subsection{$A_1$ case}\label{a1}

For the rank-one case the existence of a BA function is a very
simple fact, since in this case the number of "free" parameters
equals the number of conditions \eqref{axzpsi}. Let us consider
the root system $R=\{\pm 2\}\subset \mathbb R$ with $m_\alpha=m\in
\mathbb Z_+$. It is convenient to fix a scalar product on
$V=\mathbb R$ as $(u,v)=\frac12 uv$. In this case we will have
$R=R^\vee$,\ $Q=Q^\vee=2\mathbb Z$  and $P=P^\vee=\mathbb Z$. We
will fix $R_+=\{2\}$, so $\rho=\rho^\vee=m$. In accordance with
\eqref{psi}, the Baker--Akhiezer function $\psi$ depends on two
scalar variables $x, z$ and  has the following form:
\begin{equation}\label{psi1}
\psi=q^{xz}\sum_{\nu=-m}^m \psi_\nu q^{\nu z}\,.
\end{equation}
By definition, it must satisfy the following conditions:
\begin{equation}\label{axzpsi1}
\left(T_z^j-T_z^{-j}\right)\psi=0 \qquad\text{for each}\
j=1,\dots,m\ \text{\ and}\quad q^{2z}=1\,.
\end{equation}
Similarly to \eqref{lin}, these conditions lead to the following
linear system for the coefficients $\psi_\nu$:
\begin{equation}\label{lin1}
 \sum_{j=0}^m a_{ij}\psi_{-m+2j}=0\,, \qquad i=1\dots m \,,
\end{equation}
 where $a_{ij}=(q_j)^i-(q_j)^{-i}$ with $ q_j := q^{-m+2j+x}$.
Introduce $m\times(m+1)$ matrix $A=(a_{ij})$, then the linear
system above takes the form $Av=0$, where $v$ is the column
$v=(\psi_{-m},\psi_{-m+2}, \dots, \psi_{m})$. We know already that
for generic $x$ the matrix $A$ has the maximal rank (equal to
$m$), hence, its kernel is one-dimensional. Using Cramer's rule,
we find the values of $\psi_\nu$ (up to a common factor):
$$\psi_{-m+2s}=(-1)^j\det A^{(s)}\,,$$ where $A^{(s)}$ is obtained
from $A$ by deleting its $s$-th column. This gives us the ratio
$\psi_{-m+2s}/\psi_{-m}$ as $$ \psi_{-m+2s}/\psi_{-m}=(-1)^s\det
A^{(s)}/\det A^{(0)}\,.$$ To calculate this explicitly, we use the
following lemma.

\medskip
\begin{lemma}\label{det}
For arbitrary $q_1,\dots, q_n$ consider the matrix
$A=(a_{ij})_{i,j=1\dots n}$ whose entries are
$a_{ij}=(q_j)^i-(q_j)^{-i}$. Then $$ \det A= \prod_{i<j}\left(
q_i^{1/2}q_j^{1/2}-q_i^{-1/2}q_j^{-1/2}\right)
\left(q_i^{-1/2}q_j^{1/2}-q_i^{1/2}q_j^{-1/2}\right)\prod_i
\left(q_i-q_i^{-1}\right)\,. $$
\end{lemma}
Expanding the determinant, it is easy to see that this formula is
equivalent to the Weyl denominator formula for the $C_n$-root
system. There is also a simple direct way of proving it, using
that $\frac{q^i-q^{-i}}{q-q^{-1}}$ is polynomial in $q+q^{-1}$.

\medskip

Applying the lemma, we calculate the determinants and find
that $$
\frac{\psi_{-m+2s}}{\psi_{-m}}=\frac{\left(q_0-q_0^{-1}\right)\prod_{j=1}^{m_\alpha}\left(
q_0^{1/2}q_j^{1/2}-q_0^{-1/2}q_j^{-1/2}\right)
\left(q_0^{-1/2}q_j^{1/2}-q_0^{1/2}q_j^{-1/2}\right)
}{\left(q_s-q_s^{-1}\right)\prod_{\genfrac{}{}{0pt}{}{j\ge 0}{j\ne
s} }\left( q_s^{1/2}q_j^{1/2}-q_s^{-1/2}q_j^{-1/2}\right)
\left(q_s^{-1/2}q_j^{1/2}-q_s^{1/2}q_j^{-1/2}\right) }\,. $$
Substituting $q_j= q^{-m+2j+ x}$ we arrive after simple
transformations at the formula
\begin{equation}\label{ratio1}
\frac{\psi_{-m+2s}}{\psi_{-m}}=\prod_{j=1}^{s}\frac{\left(q^{-m+j-1}-q^{m-j+1}
\right)\left(q^{-m+j-1+x}-q^{m-j+1-x}\right)}
{\left(q^j-q^{-j}\right)\left(q^{j+x}-q^{-j-x}\right)}\,.
\end{equation}
In particular, for $s=m$ we have :
\begin{equation}\label{vert1}
\psi_{m}/\psi_{-m}=\prod_{j=1}^{m} \frac{q^{j-x}-q^{-j+x}}
{q^{j+x}-q^{-j-x}}\,.
\end{equation}
Let us fix $\psi_m$ in the following form:
$$\psi_{m}=\prod_{j=1}^m\left(q^{j-x}-q^{-j+x}\right)\,.$$ Then
the relation \eqref{vert1} gives us that
\begin{gather}\label{norm1}
\psi_{-m}=\prod_{j=1}^m\left(q^{j+x}-q^{-j-x}\right)\,.
\end{gather}
Note that our choice of $\psi_m$ implies that all $\psi_\nu$ will
be Laurent polynomials in $q^x$ (see formula \eqref{ratio1}).

This allows us to prove the following proposition.

\begin{prop}\label{r1}
(i) The function $\psi(x,z)$ given by
\eqref{psi1}, \eqref{ratio1} and \eqref{norm1} is a
Baker--Akhiezer function for $R=A_1$;

\noindent (ii) $\psi$ satisfies the difference equation $
L\psi=(q^z+q^{-z})\psi$\,,where $L$ coincides with the operator
\eqref{rank1}:$$L=\frac{q^{x-m}-q^{-x+m}}{q^x-q^{-x}}T_x+
\frac{q^{x+m}-q^{-x-m}}{q^x-q^{-x}}T_x^{-1}\,;$$

\noindent (iii) $\psi$ is symmetric in $x$ and $z$:\
$\psi(x,z)=\psi(z,x)$.

\end{prop}

\medskip
\begin{proof}
Part (i) is proven above. To prove (ii) we apply the standard
argument due to Krichever \cite{Kr}. Namely, let us consider the
function $\phi=L\psi-(q^z+q^{-z})\psi\,$. The first remark is that
$\phi$ still satisfies the conditions \eqref{axzpsi1}:\,$
\left(T_z^j-T_z^{-j}\right)\phi=0$ for each $j=1,\dots,m$ and
$q^{2z}=1$. Indeed, it is obvious for $L\psi$ since the operator
$L$ does not involve $z$. Further, $f(z)=q^z+q^{-z}$ satisfies the
conditions \eqref{axzpsi1}, hence, $f(z)\psi$ will satisfy them,
too.

 Our second remark is that $\phi$ can be presented as follows:
\begin{equation}\label{phi}
\phi=q^{xz}\sum_{j=0}^{m+1} \phi_j q^{(-m-1+2j) z}\,.
\end{equation}
This follows directly from the formula \eqref{psi1} and our
definition of $\phi$.

Let us calculate now the coefficient $\phi_0$ using
\eqref{norm1}and the definition of $\phi$. This gives:
$$\phi_0=\frac{q^{x+m}-q^{-x-m}}{q^x-q^{-x}}\psi_{-m}(x-1)-\psi_{-m}(x)=0\,.$$
In the same way, $\phi_{m+1}=0$. Thus, the expansion \eqref{phi}
contains $m$ terms only, hence, it must be zero due to Lemma
\ref{lemma1}. This proves part (ii).

To prove part (iii), first notice that according to \eqref{ratio1}
and \eqref{norm1} $\psi$ has no singularities in the $x$-variable
and it may be presented as
\begin{equation}\label{psi2}
\psi=q^{xz}\sum_{\nu=-m}^m a_\nu(z) q^{\nu x}\,.
\end{equation}
We know that, as a function of $x$, $\psi$ is an eigenfunction of
the operator \eqref{rank1}. Invoking Lemma \ref{inverse}, we
conclude that $\psi$ satisfies the following conditions in $x$:
\begin{equation}\label{axxpsi1}
\left(T_x^j-T_x^{-j}\right)\psi=0 \qquad\text{for each}\
j=1,\dots,m\ \text{\ and}\quad q^{2x}=1\,.
\end{equation}
Thus, $\psi(x,z)$ must coincide with $\psi(z,x)$ up to a
$z$-depending factor: $\psi(x,z)=F(z)\psi(z,x)$. Switching $x$ and
$z$, we conclude that $\psi(z,x)=G(x)\psi(x,z)$. This implies that
$F=G^{-1}$ is constant. Expanding the coefficients $\psi_\nu(x)$
with the help of the formulas \eqref{ratio1} and \eqref{norm1}, we
see that $\psi(x,z)$ contains the term
$$(-1)^mq^{-m(m+1)}q^{xz}q^{mx+mz}\,.$$ Since this term is
symmetric in $x$ and $z$, we conclude that $F=1$ and
$\psi(x,z)=\psi(z,x)$.

\end{proof}

\subsection{Normalized BA-function}
Now we are going to extend the results of the previous subsection
to the higher rank case. Above we have proved that a
Baker--Akhiezer function $\psi$ is determined uniquely (up to a
$x$-depending factor) by its properties
\eqref{psi}--\eqref{axzpsi}.  Let us impose the following
normalization condition on $\psi$, prescribing its leading
coefficient $\psi_{\rho^\vee}$ to be the following:
\begin{equation}
\label{norm0} \psi_{\rho^\vee}=\prod_{\alpha\in
R_+}\prod_{j=1}^{m_\alpha} [j-(\alpha, x)]\,,\qquad
[a]:=q^a-q^{-a}\,.
\end{equation}

\medskip
\begin{definition}
A {\bf normalized} BA function is a (unique) function $\psi(x,z)$
with the properties \eqref{psi}--\eqref{axzpsi} and normalization
\eqref{norm0}.
\end{definition}

\medskip
Let us discuss briefly the geometry of the polytope $\mathcal
N={\rm supp}(\psi)$ defined by \eqref{nu}.
We mentioned already that $\mathcal N$
is the Minkowski sum of the segments associated with the vectors
$m_\alpha\alpha^\vee$\ ($\alpha^\vee\in R_+^\vee$). It is
convenient to use a more symmetric definition of $\mathcal N$,
which is obviously equivalent to \eqref{nu}:
\begin{equation}
\label{nuu} \mathcal N = \{ \nu= \sum_{\alpha \in R_+} l_\alpha
\alpha^\vee \mid -\frac12 m_\alpha\le l_\alpha \le \frac12
m_\alpha \}\,.
\end{equation}
To better understand its structure, let us choose a generic
direction $v$ in $V$ and consider the height function $(v,\cdot)$
on $V$. Then the highest and the lowest among the points of
$\mathcal N$ will be the points $\frac12\sum_{\alpha \in R_+} \pm
m_\alpha \alpha^\vee$, where the signs in the sum either coincide
with the signs of $(\alpha^\vee, v)$, or are exactly the opposite.
This shows that the vertices of $\mathcal N$ have the form
$\frac12\sum_{\alpha \,:\,(\alpha,v)>0}m_\alpha \alpha^\vee$,
where $v$ is a generic vector. For instance, taking $v$ from the
{\it positive Weyl chamber} $C$, i.e. such that $(v,\alpha)>0$ for
all $\alpha\in R_+$, we obtain $\rho^\vee$ as one of the vertices
of $\mathcal N$.
Other vertices will correspond to other chambers and will be of
the form $w\rho^\vee$ with $w\in W$. If we take now the vertex
$\rho^\vee$, then its adjacent vertices will correspond to the
Weyl chambers, adjacent to the positive one. There are exactly $n$
of them, $n={\rm rank}R$, and they are the images $s_iC$ of $C$
under the simple reflections $s_1,\dots, s_n$.
It is known that a simple reflection $s_i$ leaves
the set $R_+\backslash\alpha_i$ invariant, while sending
$\alpha_i$ to $-\alpha_i$. Thus, $s_i\rho^\vee=\rho^\vee-m_i
\alpha_i^\vee$,\ $m_i:=m_{\alpha_i}$. By the way, this implies
that $(\alpha_i,\rho^\vee)=m_i$, i.e. $\rho^\vee$ can be rewritten
in terms of the fundamental coweights $b_i$ as
\begin{equation}\label{cow}
\rho^\vee=\sum_{i=1}^n m_ib_i\,,\qquad m_i=m_{\alpha_i}\,.
\end{equation}
 So, the edges of $\mathcal N$ coming out from the vertex $\rho^\vee$ are
given by the vectors $-m_i\alpha_i^\vee$ where
$\alpha_1,\dots,\alpha_n$ are the simple roots from $R_+$. In the
same way, the vectors $-m_iw\alpha_i^\vee$ will lead from the
vertex $w\rho^\vee$ to its adjacent vertices. Summarizing, we
arrive at the following result.

\begin{lemma}\label{nuvert}
The vertices of the polytope \eqref{nuu} have the form
$w\rho^\vee$,\ $w\in W$. If $\nu,\,\nu'$ are two adjacent
vertices, then $\nu'$ equals $\nu+m_\alpha\alpha^\vee$ for a
proper $\alpha\in R$. Moreover, in this case we have
$\nu'=s_\alpha \nu$ and $(\alpha, \nu')=-(\alpha, \nu)=m_\alpha$.
\end{lemma}

\bigskip
Our next proposition justifies our choice of normalization
\eqref{norm0}.

\begin{prop}\label{no}
The normalized BA function $\psi$ has the following properties:

\noindent (i) for all $w\in W$ the coefficient $\psi_\nu$ in the
vertex $\nu=w\rho^\vee$ of the polytope $\mathcal N$ has the form
\begin{equation}
\label{norm} \psi_{w\rho^\vee}=\prod_{\alpha\in
wR_+}\prod_{j=1}^{m_\alpha} [j-(\alpha, x)]\,;
\end{equation}

\noindent (ii) for all $\nu$ lying on the (one-dimensional) edges
of $\mathcal N$ the corresponding $\psi_\nu$ are polynomial in
$x$;

\noindent (iii) $\psi$ is quasipolynomial in both $x$ and $z$.
\end{prop}

\begin{proof}
Let $\nu$ and $\nu'$ be two adjacent vertices of $\mathcal N$, so
for a proper $\alpha\in R$ we have $\nu'=s_{\alpha^\vee}(\nu)=
\nu+m_\alpha\alpha^\vee$. Let us introduce
$\psi_j=\psi_{\nu+j\alpha^\vee}$ {and} $q_j= q^{2j-m_\alpha+
(\alpha, x)}$ ($j=0,1,\dots, m_\alpha$). Introduce also a matrix
$A$ with entries $a_{ij}=(q_j)^i-(q_j)^{-i}$\,($1\le i\le
m_\alpha$,\,$0\le j\le m_\alpha$). In these notations, the
conditions \eqref{axzpsi} can be expressed as $Av=0$, where $v$ is
the column $v=(\psi_0, \dots, \psi_{m_\alpha})$. This system is
completely analogous to the system  \eqref{lin1}. Repeating the
same arguments, we arrive at the formula
\begin{equation}\label{ratio}
\psi_s/\psi_0=\prod_{j=1}^{s}\frac{[-m_\alpha+j-1]\,
[(\alpha,x)-m_\alpha+j-1]} {[j][j+(\alpha,x)]}\,,
\end{equation}
where $[a]$ denotes $q^a-q^{-a}$. In particular, for $s=m_\alpha$
we have the following expression for the ratio of the coefficients
$\psi_\nu$ and $\psi_{\nu'}$ at two adjacent vertices of the
polytope $\mathcal N$:
\begin{equation}\label{vert}
\psi_{\nu'}/\psi_\nu=\prod_{j=1}^{m_\alpha} \frac{[j-(\alpha, x)]}
{[j+(\alpha, x)]}\,.
\end{equation}

Recall now that $\nu, \nu'$ can be presented as $w\rho^\vee$,
$w'\rho^\vee$ for proper $w, w' \in W$. Moreover, we have
$w'=s_{\alpha}w$, and $\pm\alpha$ are simple roots in the sets
$w'R_+$ and $wR_+$, respectively. This means that if we denote the
common part of the sets $wR_+, w'R_+$ by $S$, then $wR_+=S\cup
\{-\alpha\}$ and $w'R_+=S\cup \{\alpha\}$. Now the formula
\eqref{norm} in question gives us that $$
\psi_{w'\rho^\vee}/\psi_{w\rho^\vee}= \prod_{j=1}^{m_\alpha} \frac
{[j-(\alpha, x)]} {[j+(\alpha, x)]}\,. $$ Since this expression
coincides with the formula \eqref{vert}, and since \eqref{norm} is
valid at one vertex $\nu=\rho^\vee$, we conclude that it must be
valid at all other vertices of $\mathcal N$. This proves part (i)
of the proposition.

The second part follows directly from the formulas \eqref{norm}
and \eqref{ratio}. To prove part (iii), we have to show that all
coefficients $\psi_\nu$ are polynomial in $x$. The BA function in
formula \eqref{ber} is clearly quasipolynomial in $x$, so we only
have to prove the absence of singularities in the normalized BA
function $\psi$.

Suppose $\psi$ has a pole along a hypersurface $F(x)=0$. Without
loss of generality, we may assume that it is irreducible.
Multiplying $\psi$ by $F(x)$ we will obtain the function
$\widetilde\psi$ which will be quasipolynomial in $x$ sharing with
$\psi$ the properties \eqref{psi}-- \eqref{axzpsi}. Now take
generic point $x_0$ on this hypersurface, so $F(x_0)=0$. Then all
the coefficients $\widetilde\psi_\nu=F(x)\psi_\nu$ for $\nu$ lying
on the edges of the polytope $\mathcal N$ must vanish at $x=x_0$
(because $\psi_\nu$ are polynomial, see part (ii)). However,
$\widetilde\psi$ does not vanish at $x_0$ (otherwise $\psi$ would
be nonsingular at $F=0$). Let us denote by $\widetilde {\mathcal
N}$ the support of $\widetilde\psi(x_0,z)$ (i.e. the convex hull
of all $\nu$ with nonzero $\widetilde\psi_\nu(x_0)$). We have,
therefore, a strict inclusion
\begin{equation}\label{incl}
\widetilde {\mathcal N}\subsetneqq {\mathcal N}\,.
\end{equation}
Moreover, this shows that the polytope $\widetilde {\mathcal N}$
has no vertices lying on the edges of ${\mathcal N}$. This will
lead us to a contradiction in a moment.

Indeed, take any $\alpha\in R$ and consider the function
$f(x)=q^{2(\alpha,x)}$ restricted to the hypersurface $H:\ F=0$.
Since $H$ is irreducible, $f$ either is constant along $H$ or
takes infinitely many different values. Suppose now that for every
$\alpha\in R$ the latter alternative holds. Then we would be able
to choose generic $x_0\in H$ such that our proof of lemma
\ref{lemma1} would work, i.e. such that all $m_\alpha\times
m_\alpha$ determinants arising from considering $\alpha$-strings
inside the polytope $\mathcal N$ would be nonzero. However, this
would imply lemma \ref{lemma2} and this would contradict the
inclusion \eqref{incl}.

The upshot is that for some $\alpha\in R$ the function
$q^{2(\alpha,x)}$ is constant along our hypersurface $H$. But in
this case we can choose generic $x_0\in H$ such that lemma
\ref{lemma1}  is still applicable to all roots except $\alpha$.
Hence, the support $\widetilde {\mathcal N}$ of $\widetilde \psi$
due to lemma \ref{lemma2} must contain a polytope $\mathcal N'$,
which by definition is the Minkowski sum of all $m_\beta
\beta^\vee$\ ($\alpha\ne\beta\in R_+$). However, $ {\mathcal
N}=\mathcal N'\#m_\alpha\alpha^\vee$. This means that each copy of
${\mathcal N}'$ lying inside ${\mathcal N}$, has at least one
vertex on the edges of ${\mathcal N}$, but $\widetilde {\mathcal
N}$ has no vertices on the edges of ${\mathcal N}$. This
contradiction proves that $\psi$ has no singularities, therefore,
it is quasipolynomial in $x$.
\end{proof}

\bigskip
So, our choice of normalization \eqref{norm0} guarantees that
$\psi$ will be quasipolynomial in $x$. In a certain sense, it is
the "minimal" quasipolynomial BA function.

\begin{cor}\label{must}
Any BA function $\psi(x,z)$ which is quasipolynomial in $x$ is the
normalized BA function multiplied by some polynomial in $x$.
\end{cor}

This follows immediately from the formula \eqref{norm} which shows
that all the coefficients $\psi_{w\rho^\vee}$ have no common
divisors (as polynomials in $x$).

\subsection{Duality}

A remarkable and important property of the normalized
Baker--Akhiezer function is a certain duality between $x$ and $z$
variables. Namely, it turns out that switching $x$ and $z$ leads
to the BA function for the {\it dual} root system. In particular,
for $A,D,E$ root systems we have simply $\psi(x,z)=\psi(z,x)$. The
next proposition is the main step in establishing this symmetry.

\medskip
\begin{prop}\label{dual1}
The normalized BA function $\psi$ satisfies the following
conditions in $x$: for each $\alpha\in R_+$ and $s=1,\dots
,m_\alpha$
\begin{equation}\label{axxpsi}
\psi\left(x+\frac 12 s\alpha^\vee, z\right)\equiv
\psi\left(x-\frac 12 s\alpha^\vee, z\right)\quad\text{for}\ \
q^{2(\alpha,x)}=1\,.
\end{equation}
\end{prop}

\medskip
\begin{proof} Take any $\alpha\in R$ and consider two adjacent
vertices $\nu,\,\nu'$ of the polytope $\mathcal N$ such that
$\nu'=\nu+m_\alpha\alpha^\vee$. As we know, $\nu$ and $\nu'$ can
be presented as $w\rho^\vee$, $w'\rho^\vee$ for proper $w, w' \in
W$. Moreover, denoting the set $w'R_+\cup wR_+$ by $S$, we will
have: $$wR_+=S\cup \{-\alpha\}\,,\quad w'R_+=S\cup
\{\alpha\}\,,\quad s_\alpha(S)=S\,.$$

Consider the function $\widetilde\psi=\left(T_x^{\frac12
s\alpha^\vee}-T_x^{-\frac12 s\alpha^\vee}\right)\psi$. Similar to
$\psi$, it will be quasipolynomial in $z$,
\begin{equation}\label{sh}
\widetilde\psi = q^{2(x,z)}\sum_{\tau\in P^\vee} a_\tau
q^{2(\tau,z)}\,,
\end{equation}
and it is clear that its support lies inside the polytope
\begin{equation}\label{min}
\widetilde{\mathcal N}=\mathcal N'\#l\,,
\end{equation}
where $l$ is the segment with the endpoints at $\nu-\frac12
s\alpha^\vee$ and $\nu'+\frac12 s\alpha^\vee$ and $\mathcal N'$,
in its turn, is the Minkowski sum of all $m_\beta\beta^\vee$ with
$\beta\in S$ (we recommend the reader to draw a picture for $A_2$ case).

Now take $x_0$ such that $q^{2(\alpha,x_0)}=1$. We would like to
prove that for such a $x_0$\  $\widetilde\psi$ will be zero.
Suppose that it is not the case. We may assume that $x_0$ is
generic enough, i.e. such that we still can apply Lemma
\ref{lemma1} to all $\beta\in S$. This, together with Lemma
\ref{lemma2}, implies that ${\rm supp}(\widetilde\psi)$ can be
presented as $\mathcal N'\#\Omega$ for some domain $\Omega$.
Comparing with \eqref{min} gives us that $\Omega\subset l$. This
implies that at least one of the coefficients $a_\tau$ in
\eqref{sh} with $\tau\in l$ must be nonzero. So, we will arrive at
a contradiction as soon as we prove that all $a_\tau$ for $\tau\in
l$ vanish. We will see in a moment that this essentially reduces
the problem to the rank-one case.

Indeed, notice that $q^{2(x,z)}\sum_{\tau\in l}a_\tau
q^{2(\tau,z)}$, by its construction, coincides with
$\left(T_x^{\frac12 s\alpha^\vee}-T_x^{-\frac12
s\alpha^\vee}\right)\psi^{(0)}$, where $$\psi^{(0)}:=
q^{2(x,z)}\sum_{\tau\in [\nu, \nu']} \psi_\tau q^{2(\tau,z)}$$ is
the part of $\psi$, corresponding to the edge $[\nu,\nu']$ of the
polytope $\mathcal N$. We recall that we consider $x=x_0$ such
that $q^{2(\alpha,x)}=1$.

So, we need to show that
\begin{equation}\label{ax0}
\left(T_x^{\frac12 s\alpha^\vee}-T_x^{-\frac12
s\alpha^\vee}\right)\psi^{(0)}=0\qquad \text{for}\
q^{2(\alpha,x)}=1\,.
\end{equation}
 If our $\psi^{(0)}$ were the normalized BA
function for the root system $R^{(0)}=\{\pm \alpha\}$, this would
follow directly from Proposition \ref{r1}(ii). This is almost the
case, as the formulas \eqref{ratio}--\eqref{vert} clearly
indicate. The only difference is that, due to \eqref{norm}, the
coefficient $\psi_\nu$ with $\nu=w\rho^\vee$ looks as follows: $$
\psi_\nu=\prod_{\beta\in S}\prod_{j=1}^{m_\beta}[j-(\beta,
x)]\prod_{j=1}^{m_\alpha} [(\alpha, x)+j]$$ (here we used that
$wR_+=S\cup \{-\alpha\}$). So, the difference with the similar
formula \eqref{norm1} for the rank one case comes from the factor
$\prod_{\beta\in S}\prod_{j=1}^{m_\beta} [j-(\beta,x)]$. However,
this factor is invariant under reflection $s_\alpha$, hence, it
does not affect the properties \eqref{ax0}.
\end{proof}

\medskip
Next thing to prove is that the support of $\psi$ in the
$x$-variable is the following polytope $\mathcal N^\vee$:
\begin{equation}
\label{nuvee} \mathcal N^\vee = \{ \nu= \rho-\sum_{\alpha \in R_+}
l_\alpha \alpha \mid 0\le l_\alpha \le m_\alpha \}\,.
\end{equation}
To prove this we look first at the formula \eqref{ber}, which
gives an expression for another BA function with
$\psi_{\rho^\vee}$ normalized as in \eqref{bernorm}. It is clear
from \eqref{ber} that such $\psi$ is also quasipolynomial in $x$
with ${\rm supp}_x\psi={\rm supp}_x\psi_{\rho^\vee}$. If we
renormalize now $\psi$ in order to get $\psi_{\rho^\vee}$ as in
\eqref{norm0}, then we will still have that the support of $\psi$
equals the support of $\psi_{\rho^\vee}$ which now will coincide
with $\mathcal N^\vee$. So, for the normalized BA function
$\psi(x,z)$ one has ${\rm supp}_x\psi=\mathcal N^\vee$.

As a result, we see that the normalized BA function can be
presented in a form $$\psi=q^{2(x,z)}\sum
c_{\mu\nu}q^{2(\mu,x)}q^{2(\nu,z)}$$ with the summation taken over
$\mu\in\mathcal N^\vee$ and $\nu\in\mathcal N$, respectively. The
highest term in this expression will correspond to $\mu=\rho$,
$\nu=\rho^\vee$.

\bigskip
Summarizing, we see that the properties of $\psi(x,z)$ in the
$x$-variable are completely analogous to its properties in $z$.
Thus, we obtain the following important property of the normalized
BA function which reflects the symmetry between $x$ and $z$
variables.

\medskip
\begin{theorem}[Duality]\label{symmetry}
Let $\psi(x,z)$ denote the normalized Baker--Akhiezer function
associated to a root system $R$ with multiplicities
$m=\{m_\alpha\}$, and $\psi^\vee(x,z)$ denote a similar function
associated to the dual system $R^\vee$ with
$m_{\alpha^\vee}=m_\alpha$. Then $\psi(x,z)=\psi^\vee(z,x)$.
\end{theorem}

\bigskip
\begin{cor}\label{bsp}  Let $\psi$ be the normalized BA function
for a given root system $R$ and $m=\{m_\alpha\}$. Let $\omega\in
P$ and $\pi\in P^\vee$ be (quasi)minuscule weight and coweight for
the root system $R$. Then $\psi(x,z)$ solves the following
bispectral system of difference equations:
\begin{equation}
\label{bs} \left\{
\begin{array}{l}
D^\omega_z\psi = \mathfrak m_\omega(x)\psi\,, \\ D^\pi_x\psi =
\mathfrak m_\pi(z)\psi\,.
\end{array}
\right.
\end{equation}
Here $D^\pi_x$ and $D^\omega_z$ are the corresponding Macdonald
operators \eqref{m1}-\eqref{m2} and \eqref{m1z}-\eqref{m2z}, while
$\mathfrak m_\lambda$ stands for the orbitsum \eqref{orb}.
\end{cor}

\section{Algebraic integrability and applications}

\subsection{Algebraic integrability}

Above we have shown that the normalized Baker--Akhiezer function
$\psi(x,z)$ associated to a datum $\{R,m\}$ is an eigenfunction of
the Macdonald operators $D^\pi_x$, where $\pi$ is any
(quasi)minuscule coweight for the root system $R$. In fact, $\psi$
is a common eigenfunction of a much bigger commutative ring of
difference operators. This follows in a standard way from its
analytic properties in the $z$-variable (cf. \cite{Kr}). To
formulate the result, let us recall the definition of the ring
$\mathfrak R^\vee$ which consists of all polynomials $f(z)$ of the
form $ f(z)=\sum_{\nu\in P^\vee} f_\nu q^{2(\nu,z)}$, satisfying
conditions \eqref{axz} for each $\alpha\in R$ and $s=1,\dots
,m_\alpha$. It is easy to see that all $W$-invariant polynomials
belong to the ring $\mathfrak R^\vee$. According to Chevalley
theorem, $W$-invariants form a polynomial ring with $n={\rm
rank}R$ generators $m_{b_1},\dots,m_{b_n}$, which are the
orbitsums for the fundamental coweights $b_1,\dots,b_n$. However,
the ring $\mathfrak R^\vee$ is much bigger, for instance, it
contains the principal ideal generated by the polynomial $Q$ given
by \eqref{Q}.

\medskip
\begin{theorem}[Algebraic integrability]\label{alin}
For each polynomial $f(z)$ from the ring $\mathfrak R^\vee$ there
exists a difference operator $D_f$ in $x$ on the lattice $P^\vee$
such that $ D_f \Psi = f(z)\Psi$.
All these operators commute. For a
(quasi)minuscule coweight $\pi$ of $R$ and $f=\mathfrak m_\pi(z)$
the corresponding operator $D_f$ coincides with the Macdonald
operator $D^\pi$ given in \eqref{m}--\eqref{m2}.
\end{theorem}

\medskip
\begin{proof}
Everything is based on the following result.
\begin{lemma}\label{prin}
Any quasipolynomial in $z$ of the form
$\phi=q^{2(x,z)}\sum_{\nu\in P^\vee}\phi_\nu(x)q^{2(\nu,z)}$ which
satisfies the conditions \eqref{axzpsi} can be obtained  by
applying a proper difference operator in $x$ to the BA function
$\psi(x,z)$.
\end{lemma}

To prove the lemma, we recall that according to lemmas
\ref{lemma1} and \ref{lemma2} the support of $\phi$ can be
presented as $N_0\#\mathcal N$, where $\mathcal N$ is the support
of $\psi$ and $N_0$ is some convex polytope. Now choose any vertex
$\nu$ of the ${\rm supp}\phi$ and let $\tau$ be the corresponding
vertex of $N_0$ such that $\nu\in \tau+\mathcal N \subseteq {\rm
supp} \phi$. We can "kill" the coefficient $\phi_\nu$ of $\phi$ by
subtracting a function $T^\tau_x\psi$ taken with a proper
coefficient $d(x)$. The resulting function
$\widetilde\phi=\phi-d(x)T^\tau_x\psi$ still satisfies the
conditions of the lemma, but it has a smaller support. Repeating
this, we will eventually get zero, and this proves the lemma.

Now, to prove the theorem, we notice that for $f\in\mathfrak
R^\vee$ the function $\phi=f(z)\psi$ will satisfy the conditions
of the lemma, hence, $\phi=D_f\psi$ for a proper difference
operator $D_f=\sum_{\nu\in P^\vee}d_\nu T^\nu_x$.
All these operators commute since they have $\psi$ as their
eigenfunction. Indeed, for $f,g\in \mathfrak R^\vee$ we have
$[L_f,L_g]\psi=(fg-gf)\psi=0$. However, one shows easily that if a
difference operator $M$ in $x$ annihilates a (nonzero)
quasipolynomial $\psi$ in $z$, then $M=0$ (otherwise ${\rm
supp}M\psi$ would be nonempty). Hence, $[L_f,L_g]=0$.

Finally, for $f(z)=\mathfrak m_\pi(z)$ and the corresponding $D_f$
we will have $D_f\psi=\mathfrak m_\pi(z)\psi$. The same is true
for the Macdonald operator $D^\pi_x$:\quad $D^\pi_x\psi=\mathfrak
m_\pi(z)\psi$. Hence, these two operators coincide.
\end{proof}

\bigskip
\begin{remark} The duality between $x$ and $z$ implies that the
normalized BA function will be also a common eigenfunction of a
"dual" commutative ring of difference operators in $z$ variable,
isomorphic to the ring $\mathfrak R$ of polynomials with the
properties \eqref{axx}. Thus, we have a {\it bispectral pair} of
commutative rings, in the spirit of \cite{W2}.
\end{remark}

\medskip
\begin{remark}
In the limit $q\to 1$ this theorem reduces to the result from
\cite{VSC}, where the algebraic integrability was established for
the quantum trigonometric Calogero--Sutherland--Moser problem and,
more generally, for its generalizations \cite{OP} related to the
root systems. Proof in \cite{VSC} was based on results by Heckman
and Opdam \cite{HO,H}, who developed a nice theory of
multivariable hypergeometric functions related to root systems. An
independent direct proof was obtained in \cite{C00}. For $A_{n-1}$
case, i.e. for the Macdonald--Ruijsenaars operators \eqref{ru}
with $t=q^m$\ ($m\in\mathbb Z$) the algebraic integrability was
established by Etingof--Styrkas in \cite{ES}. They constructed
$\psi$ in terms of certain intertwiners from representation theory
for quantum groups. Their formula for $\psi$ was made more explicit
by Felder and Varchenko in \cite{FV}.
\end{remark}

\subsection{Liouville
integrability}

As a part of the previous theorem, we have constructed $n={\rm
rank}R$ commuting difference operators $D_i$ corresponding to the
basic $W$-invariants $f=\mathfrak m_{b_i}(z)$, which are the
orbitsums \eqref{orb} for the fundamental coweights $b_i\in
P^\vee$. Those of $b_i$ which are (quasi)minuscule, will lead to
the Macdonald operators \eqref{m}, \eqref{m1} with
$t_\alpha=q^{-m_\alpha}$. For others, there is no such simple
explicit formula. However, it is easy to evaluate leading terms in
$D_i$ and they look similar to the leading terms in Macdonald
operators.

\medskip
\begin{prop}\label{int}
For any coweight $\pi\in P^\vee$ and its orbitsum $f=\mathfrak
m_\pi(z)$ the operator $D_f$ constructed in Theorem \ref{alin}
will have the form $$
 D_f = \sum_{\tau\in W\pi}a_\tau
T^{\tau}+\dots\,,$$ where the dots stand for a sum of the terms
$a_\nu T^\nu$ with $\nu\in \pi+Q^\vee$ lying inside the convex
hull of the orbit $W\pi$, and the leading coefficients $a_\tau$
are given by \eqref{delta}, \eqref{coef} with
$k_\alpha=-m_\alpha$. The operator $D_f$ is $W$-invariant.
\end{prop}

\medskip
\begin{proof}
Formula for $a_\tau$ follows directly from the construction of
$D_f$ in theorem \eqref{alin} and formula \eqref{norm} for the
leading coefficients in $\psi$. It remains to prove that $D_f$ is
$W$-invariant. To this end we have the following symmetry property
of $\psi$.

\begin{lemma}\label{w}
The normalized BA function is $W$-invariant in the following
sense: $\psi(wx,wz)=\psi(x,z)$ for any $w\in W$.
\end{lemma}
To prove this we notice that for the leading coefficients
\eqref{norm} of $\psi$ one has
$\psi_{w\rho^\vee}(wx)=\psi_{\rho^\vee}(x)$. Since
$\widetilde\psi(x,z):=\psi(wx,wz)$ shares with $\psi$ the
properties \eqref{psi}--\eqref{axzpsi}, they must coincide due to
the uniqueness of $\psi$.

Using the lemma and $W$-invariance of $f(z)$, one gets that
$[D\psi](wx,wz)=f(wz)\psi(wx,wz)=[D\psi](x,z)$. From this it
easily follows the property $a_{w\nu}(wx)=a_\nu(x)$ for the
coefficients of $D$, i.e. its $W$-invariance.
\end{proof}

\bigskip
Let us look what happens if we change the multiplicities
$m_\alpha$. The formula \eqref{coef} shows that the leading
coefficient $a_\pi$ in $D_f$ will be rational in $t=\{t_\alpha\}$.
In fact, it is not difficult to prove that all the coefficients of
the operator $D_f$ with $f=\mathfrak m_\pi(z)$ will be rational in
$t$.

\medskip
\begin{lemma}\label{ad} All the operators $D_f$ depend rationally on $t$.
\end{lemma}

\begin{proof}
Recall that the construction of $D_f$ was given in Lemma
\ref{prin}. Now let $f(z)=\mathfrak m_\pi(z)$ be an orbitsum,
$\pi\in P_+$. We can construct $D_f$ as in Lemma \ref{prin},
starting from $\phi=f(z)\psi$ and "killing" at each step the
highest coefficient of $\phi$. This shows that to calculate $D_f$
in this case it is sufficient to know a fixed number of the
coefficients $\psi_\nu$ in the normalized BA function $\psi$,
namely, those with $\nu-\rho^\vee\in -\pi+\overline{W\pi}$. Thus,
if all these coefficients were rational in $t$, then $D_f$ would
be rational in $t$, too. This is not the case, however, because
already $\psi_{\rho^\vee}$ is not rational in $t$. Let us
renormalize $\psi$ in such a way that $\psi_{\rho^\vee}=1$, and
let us calculate $D_f$ using this $\widetilde\psi$. Of course, the
resulting operator $\widetilde D_f$ will differ from $D_f$, but
the relation is simple:$$\widetilde D_f=\Delta\circ D_f \circ
\Delta^{-1}\,,$$ where $\Delta$ is
\begin{equation}\label{delta00}
\Delta=\prod_{\alpha\in R_+}\prod_{j=1}^{m_\alpha}
\left(q^{-(\alpha,x)}- q^{-2j+(\alpha, x)} \right)^{-1}\,,
\end{equation}
which is a specialization of \eqref{delta} in case $t=q^{-m}$. An
important point is that despite the fact that $\Delta$ is not
rational in $t$, a ratio $$\Delta(x)/\Delta(x+\nu)$$ is rational
in $t$ for any $\nu\in P^\vee$. Therefore, it is sufficient to
prove that all the coefficients of the operator $\widetilde D_f$
depend rationally on $t$.

Summarizing, we showed that it is enough to check that for any
fixed $\tau\in Q^\vee_+$ the coefficient $\widetilde\psi_\nu$ with
$\nu={\rho^\vee-\tau}$ will be rational in $t$. Recall now that
the renormalized $\widetilde\psi$, as well as $\psi$ itself,
satisfies a difference equation in the $z$-variable: $$D^\omega_z
\widetilde\psi= \mathfrak m_\omega(x)\widetilde\psi\,,$$ derived
in Theorem \ref{exist}. Using this equation, one can calculate the
coefficients of $\widetilde\psi$ recursively, starting from
$\widetilde\psi_{\rho^\vee}$ (one gets this recursion similar to
\cite{HO}, expanding the coefficients of the operator $D^\omega_z$
in a series in $z$). Since the Macdonald operator is polynomial in
$t$, each $\widetilde\psi_\nu$ with $\nu={\rho^\vee-\tau}$,
calculated from this recursion, will be polynomial in $t$. This
completes the proof.
\end{proof}

\medskip
 As a
corollary, we obtain that all the operators $D_i$ related to the
fundamental coweights, must commute for {\it all values} of $t$
(since they commute for $t=q^{-m}$ with any integer $m$). Let us
consider one of the Macdonald operators as a Hamiltonian of the
corresponding quantum problem, which can be viewed as a
generalization of the trigonometric Ruijsenaars problem to all
root systems. Then what we just proved is the {\it complete
(Liouville) integrability} of this quantum problem. So, in this
way we recover the result obtained by Cherednik \cite{Ch1}:

\begin{cor} A quantum problem related to an arbitrary
Macdonald operator is completely integrable.
\end{cor}

\bigskip
\subsection{Shift operators}
Another important result by Cherednik is the construction of the
so-called {\it shift operators}, which are $q$-versions of the
operators constructed by Opdam for the case $q=1$. Let us explain
how they appear in our approach.

Let us consider a root system $R$ and two sets of {\it integer}
multiplicities $m, \widehat m$. Suppose that $ m\ge \widehat m$,
i.e. $ m_\alpha\ge \widehat m_\alpha$ for all $\alpha\in R$.
Denote by $\psi, \widehat\psi$ corresponding normalized BA
functions. Comparing their properties \eqref{axzpsi}, we then
apply Lemma \ref{prin} to conclude that there exists a difference
operator $S^+$ in $x$ such that $\psi=S^+\widehat\psi$. Let
$D,\widehat D$ be the Macdonald operators \eqref{m1} with
$t=q^{-m}$,\ $\hat t=q^{-\widehat m}$, so $D^\pi\psi=m_\pi\psi$
and $\widehat D^\pi\widehat\psi=m_\pi\widehat\psi$. Hence,
$$(D\circ S^+-S^+\circ \widehat D)\widehat \psi= D
\psi-S^+(m_\pi\widehat \psi)=0\,.$$ Therefore, we arrive at the
identity $$ D\circ S^+=S^+\circ \widehat D\,.$$ It means that the
operator $S^+$ intertwines two Macdonald operators with different
$t$. In the same way, $S^+$ intertwines each pair $D_i,\ \widehat
D_i$ of the corresponding operators from Proposition \ref{int}.

Consider now the case $\widehat t=qt$ (i.e. $\widehat
m_\alpha=m_\alpha-1$ for all $\alpha$). We will use subscripts
denoting by $D_m$ the Macdonald operator with $t=q^{-m}$. So, we
have shown that for any integer $m=\{m_\alpha\}$ we have two
Macdonald operators $D_{m}$,\,$D_{m-1}$ and their intertwiner
$S^+_{m}$ with the relation:
\begin{equation}\label{shift}
D_{m}\circ S^+_m=S^+_m\circ D_{m-1}\,.
\end{equation}
Similar to Lemma \ref{ad}, one proves that the intertwiner $S^+_m$
will depend algebraically on $t=q^{-m}$, so the relation
\eqref{shift} will make sense for {\it any value} of $m$.
Corresponding $S^+_m$ is called a {\it shift operator}, since it
sends eigenfunctions of $D_{m-1}$ to eigenfunctions of $D_{m}$.
{From} its construction we can calculate easily the leading terms
in $S_m$. Introduce $\varrho^\vee\in P^\vee$ as
\begin{equation}\label{varrho}
\varrho^\vee=\frac12\sum_{\alpha\in R_+}\alpha^\vee\,.
\end{equation}

\begin{prop}\label{+}
The shift operator $S^+_m$ above has the form
$$S^+_m=\sum_{\tau\in W\varrho^\vee}b_\tau T^\tau+\dots\,,$$ where
the dots stand for a sum of lower terms $b_\nu T^\nu$ with
$\nu\in\varrho^\vee+Q^\vee$ lying inside the convex hull of the
orbit $W\varrho^\vee$. The leading coefficients $b_\tau$ are given
by the formulas
$$b_{\varrho^\vee}(x)=q^{\sum_{\alpha>0}m_\alpha}\frac
{\Delta_{1-m}(x+\varrho^\vee)}{\Delta_{-m}(x)}\,,\qquad
b_{w\varrho^\vee}(x)=b_{\varrho^\vee}(w^{-1}x)\,,$$ where
$\varrho$ is given by \eqref{varrho} and $\Delta_{-m}$ stands for
the function \eqref{delta00}. The operator $S^+_m$ is
$W$-invariant and it sends the normalized BA function $\psi_{m-1}$
to $\psi_m$.
\end{prop}

 For a proof one should compare the supports and the leading
coefficients in two BA functions, $\psi_m$ and $\psi_{m-1}$.

\medskip Another shift operator will appear if we consider, as
above, two normalized BA functions $\psi_m, \psi_{m-1}$
corresponding to $m=\{m_\alpha\}$ and $\widehat m=m-1$, and
multiply $\psi_{m-1}$ by a polynomial $c_m(z)c_m(-z)$ where
$c_m(z)$ is defined by the formula:
\begin{equation}\label{*}
c_m(z)=\prod_{\alpha\in R_+}[m_\alpha-(\alpha^\vee,z)] \,,\qquad
[a]:=q^a-q^{-a}\,.
\end{equation}
It is easy to see that the resulting function
$\phi=c_m(z)c_m(-z)\psi_{m-1}$ will satisfy the conditions
\eqref{axzpsi}. Hence, due to Lemma \ref{prin}, it has the form
$S^-\psi_m$ for a proper difference operator $S^-=S^-_m$ (in the
$x$-variable). Its leading terms, again, can be calculated easily.
Similarly, we get the intertwining relation $D_{m-1}\circ
S^-_m=S^-_m\circ D_{m}$ which, again, extends analytically to all
values of $m$. Next proposition summarizes the properties of the
constructed $S^-$.

\medskip
\begin{prop}\label{-}For integer $m=\{m_\alpha\}$
there exists a difference operator $S^-=S^-_m$ which intertwines
the Macdonald operators $D_{m},\,D_{m-1}$: $$D_{m-1}\circ
S^-_m=S^-_m\circ D_{m}\,.$$ It has the form $S^-=\sum_{\tau\in
W\varrho^\vee}b_\tau T^\tau+\dots\,,$ where the dots stand for a
sum of lower terms $b_\nu T^\nu$ with $\nu\in\varrho^\vee+Q^\vee$
lying inside the convex hull of the orbit $W\varrho^\vee$. The
leading coefficients $b_\tau$ are given by the formulas
$$b_{\varrho^\vee}(x)=(-1)^{|R_+|}q^{-\sum_{\alpha>0}m_\alpha}
\frac {\Delta_{-m}(x+\varrho^\vee)}{\Delta_{1-m}(x)}\,,\qquad
b_{w\varrho^\vee}(x)=b_{\varrho^\vee}(w^{-1}x)\,,$$ in accordance
with \eqref{varrho}, \eqref{delta00} (here $|R_+|$ denotes the
total number of positive roots). Application of $S^-_m$ to the
normalized BA function $\psi_m$ gives
$S^-_m\psi_m=c_m(z)c_m(-z)\psi_{m-1}$, where $c_m$ is defined in
\eqref{*}.
\end{prop}

\bigskip
\begin{remark}
In case when $R$ consists of two orbits, $R=R_1\cup R_2$, one can
apply similar arguments for the case when multiplicities increase
for one of the orbits, i.e. $m=\widehat m+ 1_i$ with $1_i$ being
the characteristic function of $R_i\subset R$. This proves the
existence of the intertwiners, shifting from $m$ to $m \pm 1_i$.
Corresponding versions of propositions \ref{+}, \ref{-} are
straightforward.
\end{remark}
\medskip
\begin{remark}
 More generally, we can consider a pair of data $\{\widehat R,
\widehat m\}$ and $\{R, m\}$ with $\widehat R\subset R$ and with
$\widehat m_\alpha\le m_\alpha$ for all $\alpha\in \widehat R$.
This leads to the intertwiners between Macdonald operators for
{\it different} root systems. Unlike the shift operators $S_m$
above, they exist for {\it integer} values of $m$ only. As an
example, see \cite{CV} where such an intertwiner was constructed
explicitly for the case $q=1$, and $\widehat R=A_n\subset
R=A_{n+1}$.
\end{remark}
\bigskip
We conclude the section by a description of all possible shift
operators. First we need a result showing that the commutative
ring from Theorem \ref{alin} is maximal in a certain sense.

\begin{prop}\label{uno}
Let $D=D^\pi_x$ be a Macdonald operator related to a root system
$R$ and $t=q^{-m}$ with $m_\alpha\in\mathbb Z_+$ and let
$\mathfrak R^\vee$ be the ring of polynomials in $z$ with the
properties \eqref{axz}. Suppose that a difference operator $L$
commutes with $D$ and has rational coefficients. Then $L$ is one
of the operators constructed in Theorem \ref{alin}, i.e. $L=D_f$
for some $f\in \mathfrak R^\vee$.
\end{prop}

\begin{proof} Suppose that $[L,D]=0$ and consider a function
$\Phi=L\psi$ where $\psi$ is the normalized BA function. Under
assumptions of the theorem, we have that $\Phi(x,z)$ is of the
form $\Phi=q^{2(x,z)}\varphi$ where $\varphi$ is polynomial in $z$
and a rational function in $x$. From the commutativity of $D$ and
$L$ we get that $\Phi$ is also an eigenfunction for the Macdonald
operator $D$: $D\Phi=m_\pi(z)\Phi$. It is not difficult to deduce
that in case $t_\alpha=q^{-m_\alpha}$ with $m_\alpha\in\mathbb
Z_+$ any eigenfunction of $D=D^\pi$ either has infinitely many
poles along hyperplanes of the form $q^{(\alpha,x)}={\rm const}$,
or has no poles at all. By our assumptions, $\Phi$ cannot have an
infinite number of poles, hence, it must be quasipolynomial in
$x$. Then from the same equation, using lemma \ref{leading}, we
derive that ${\rm supp}_x \Phi$ must coincide with the polytope
\eqref{nu}. Further, with the help of lemmas \ref{inverse},
\ref{inverse1} we conclude that $\Phi$ satisfies the conditions
\eqref{axxpsi}. Hence, $\Phi(z,x)$ is a Baker--Akhiezer function
for the dual root system $R^\vee$. Now corollary \ref{must}
implies that $\Phi(z,x)=f(x)\psi^\vee(x,z)$, where
$\psi^\vee(x,z)$ is the normalized BA function for $R^\vee$ and
$f$ is some polynomial. Switching $x, z$ and using duality, we
obtain that $\Phi$ must be proportional to the normalized BA
function: $\Phi(x,z)=f(z)\psi(x,z)$. Notice that since $\Phi$ was
obtained from $\psi$ applying a difference operator in $x$, it
shares with $\psi$ the properties \eqref{axzpsi} in the
$z$-variable. Since this must be valid for any $x$, we obtain that
$f$ itself must be from the ring \eqref{axz}. Thus, $f$ belongs to
$\mathfrak R^\vee$ and $L=D_f$.
\end{proof}

\medskip
\begin{prop}\label{shone}
For positive integer $m=\{m_\alpha\}$ let $D_m$ and $D_{m-1}$ be
two Macdonald operators with $t=q^{-m}$ and $t=q^{-m+1}$,
respectively. Let $S$ be any difference operator in $x$ with
rational coefficients satisfying the intertwining relation
$D_m\circ S=S\circ D_{m-1}$. Then $S= S^+\circ L$ where $S^+$ is
the shift operator from Proposition \ref{+} and $L$ commutes with
$D_{m-1}$.
\end{prop}

\begin{proof}
Let $\psi_m$ and $\psi_{m-1}$ be the normalized BA functions for
multiplicities $m$ and $\widehat m = m-1$, respectively. From the
intertwining relation we obtain that $\widetilde\psi:=S\psi_{m-1}$
will satisfy the same difference equation as does $\psi_m$:\ $D_m
\widetilde\psi=m_\pi(z)\widetilde\psi$. Due to Proposition
\ref{uno}, $\widetilde\psi$ must have the form $f(z)\psi_m$ for
some polynomial $f(z)$. The same arguments as before give that
$f(z)$ must belong to the ring $\mathfrak R^\vee$ related to
$\widehat m=m-1$. Hence, $S\psi_{m-1}=(S^+\circ L)\psi_{m-1}$
where $L=D_f$ is such that $L\psi_{m-1}=f(z)\psi_{m-1}$. This
implies that $S=S^+\circ L$.

\end{proof}

\subsection{Relation to Macdonald polynomials}

Let $\psi(x,z)$ be the normalized BA function constructed as above
starting from the data $\{R, m \}$. Let us consider two functions
$\Phi_\pm$ obtained from $\psi$ by (anti)symmetrization in $x$: $$
\Phi_+=\sum_{w\in W}\psi(wx,z)\,,\qquad \Phi_-=\sum_{w\in
W}(-1)^w\psi(wx,z)\,. $$ Notice that the same will be the result
of (anti)symmetrization in $z$, due to Lemma \ref{w}. So, the
resulting functions $\Phi_\pm$ are (anti)symmetric in $z$, too. In
fact, the functions $\Phi_\pm(x,z)$ are closely related (for
special $z\in P_+$) to the Macdonald polynomials
$P_{\lambda}(x;q,t)$ . Let $\delta=\delta_m(x)$ be the function
defined by \eqref{deltax}, that is, $$ \delta(x)=\prod_{\alpha \in
R_+} \prod_{j=-m_\alpha}^{m_\alpha} [j+(\alpha,x)]\,. $$ Recall
also the definition \eqref{rho1} of $\rho=\rho_m$:
\begin{equation}\label{ro}
\rho_m = \frac 12 \sum_{\alpha \in R_+}m_{\alpha} \alpha\,.
\end{equation}
In particular, in case $m \equiv 1$ the vector $\rho_1$ coincides
with $\varrho$ from \eqref{varrho}.

\bigskip
\begin{theorem}\label{weyl} Let us substitute  $z=\lambda\in P_{+}$ into $\Phi_\pm$. Then
$$\Phi_+(x,\lambda)=c(\lambda,m)P_{\lambda+\rho_m}(x;q,q^{-m})$$
and $$\Phi_-(x,\lambda)=c(\lambda,m) \delta(x)
P_{\lambda-\rho_{m+1}}(x;q,q^{m+1}) \,,$$ where $\delta$ is
defined above and the factor $c(\lambda,m)$ is given by the
formula
\begin{align}\label{cc}
c(\lambda,m)=\prod_{\alpha\in
R_+}\prod_{j=1}^{m_\alpha}[j-(\alpha^\vee,\lambda)]\,.
\end{align}
We suppose that $(\lambda,\alpha^\vee)>m_\alpha$ for all
$\alpha\in R_+$ to ensure that $c(\lambda,m)\ne 0$. In other
words, $\lambda\in \rho_{m+1}+P_+$.
\end{theorem}

\begin{proof}
Let $\pi$ be any (quasi)minuscule coweight of $R$ and $D=D^\pi_x$
denote the corresponding Macdonald operator \eqref{m1} with
$t=q^{-m}$. Each of the functions
$\psi(wx,\lambda)=\psi(x,w^{-1}\lambda)$ satisfies the same
equation $D\psi=\mathfrak m_\pi(\lambda)\psi$, hence we have:
\begin{equation}\label{e}
D\Phi_\pm=\mathfrak m_\pi(\lambda)\Phi_\pm\,\qquad \mathfrak
m_\pi(\lambda)=\sum_{\tau\in W\pi}q^{2(\tau,\lambda)}\,.
\end{equation}
Using Lemma \ref{conj}, we conclude that $p(x)=\Phi_+$ and
$\widehat p(x)=\delta^{-1}\Phi_-$ will satisfy the equations
$$Dp=\mathfrak m_\pi(\lambda)p\,,\quad \widehat D \widehat
p=\mathfrak m_\pi(\lambda)\widehat p\,,$$ where $D$ and $\widehat
D$ are the Macdonald operators \eqref{m1} with $t=q^{-m}$ and
$t=q^{m+1}$, respectively. These are exactly the defining
equations \eqref{eig} for the Macdonald polynomials
$P_{\lambda+\rho_m}(x;q,q^{-m})$ and
$P_{\lambda-\rho_{m+1}}(x;q,q^{m+1})$. So, to prove the theorem,
we only have to check that $p$ and $\widehat p$ are symmetric
polynomials and to calculate their leading terms.

The $W$-invariance of $p, \widehat p$ is obvious since $\Phi_\pm$
is (anti)invariant. Calculating the leading term is quite
straightforward since we know the leading terms in
$\psi(x,\lambda)$. So the only non-trivial thing to prove is that
$\Phi_-$ is divisible by $\delta(x)$. Here we can use the
properties \eqref{axxpsi} of the BA function. Let us rewrite them
after shifting in $x$ by $-\frac12 s\alpha^\vee$ : $$\psi\left(x,
z\right)- \psi\left(x- j\alpha^\vee, z\right)=0\quad\text{for}\ \
q^{2(\alpha,x)}=q^{2j}\,,\quad j=1,\dots, m_\alpha\,.$$ Choose now
any $\alpha\in R$ and $x$ such that $(\alpha,x)=j$, then  its
image $x'=s_{\alpha}x$ under reflection $s_\alpha$ will be $x'=x-
j\alpha^\vee$, so according to the property above, we will have:
$$\psi(x,z)-\psi(s_\alpha x,z)=\psi(x,z)-\psi(x,s_\alpha
z)=0\quad\text{for}\quad(\alpha,x)=j\,.$$ Now let us split the sum
$\Phi_-(x,z)=\sum_{w\in W}(-1)^w\psi(x,wz)$ into pairs of terms
with $w,\,w'=s_\alpha w$. As a result, we obtain that
$$\Phi_-(x,z)=0\qquad\text{ for }\quad (\alpha,x)=1,\dots
m_\alpha\,.$$ $\Phi_-$ also vanishes for $(\alpha,x)=0$ due to its
antiinvariance. Invoking all $\alpha\in R$, we arrive at the
following result.

\medskip
\begin{lemma} For any $\alpha\in R_+$ the function $\Phi_-(x,z)=\sum_{w\in
W}(-1)^w\psi(x,wz)$ vanishes along the hyperplanes
$(\alpha,x)=0,\pm 1,\dots, \pm m_\alpha$.
\end{lemma}

\medskip
This is true for any $z$. However, for $z=\lambda\in P$ the
functions $\Phi_\pm$ are quasiperiodic with respect to the lattice
$\mathcal L=\omega P^\vee$ where $\omega=\pi i(\log q)^{-1}$.
Indeed, under a shift $x\to x+ l$\ ($l\in\mathcal L$) each of
$\psi(x,w\lambda)$\ ($w\in W$) gets the same factor as does the
function $q^{(\lambda+\rho,x)}$. These translation properties
imply that as soon as $\Phi_-$ vanishes for $2(\alpha,x)=2j$, it
will also vanish for $q^{2(\alpha,x)}=q^{2j}$. This proves that
$\Phi_-$ for $z=\lambda\in P$ is divisible by the polynomial
$\delta(x)$, thus, completing the proof of the theorem.

\end{proof}

\medskip
\begin{remark} The proof shows that $\Phi_-(x,\lambda)$ will be
zero for all $\lambda\in P_+$ which do not belong to
$\rho_{m+1}+P_+$. Indeed, in this case $\Phi_-$ must be divisible
by $\delta(x)$ but has a smaller support.
\end{remark}

\bigskip
Our expression for $P_{\lambda+\rho_m}(x;q,q^{-m})$ via the
function $\Phi_+$ is well-defined for sufficiently large
$\lambda$, namely for $\lambda\in \rho_{m+1}+P_+$. For smaller
$\lambda$ it doesn't work, which reflects the known fact that in
case $t=q^{-m}$ with $m_\alpha\in \mathbb Z_+$ some of $P_\mu$ are
not well-defined. Let us call the Macdonald polynomials
$P_\mu(x;q,q^{-m})$ with $(\mu,\alpha)>2m_\alpha$ for all
$\alpha>0$ the {\it stable Macdonald polynomials}. They are always
well-defined (if $q$ is not a root of unity). The following
"localization" property is a direct corollary of Theorem
\ref{weyl}.

\bigskip
\begin{cor}
For $t=q^{-m}$ with $m_\alpha\in \mathbb Z_+$ the stable Macdonald
polynomials $P_\lambda (x;q,q^{-m})$ are localized, i.e. in their
expression \eqref{macp} through the orbitsums
$$P_\lambda=\mathfrak m_\lambda+\sum_{\nu \prec
\lambda}a_{\lambda\nu} \mathfrak m_\nu\,,$$ only the terms with
$\nu=\lambda-\sum_{\alpha\in R_+}l_\alpha \alpha\,$ and
$l_\alpha=0,1,\dots, m_\alpha$ can appear, all other coefficients
$a_{\lambda\nu}$ vanish. As a result, the total number of nonzero
terms remains bounded as $\lambda$ increases.
\end{cor}

Below is a simple example which illustrates Theorem \ref{weyl}.

\medskip
\begin{example}\label{a1m1} Let us consider the case $R=A_1$ with $m=1$.
Using the results of section \ref{a1}, we obtain (switching $x,z$)
the following formula for the normalized BA function: $$\psi(x,z)=
\left(q^{1+z}-q^{-1-z}\right)q^{(z-1)x}+
\left(q^{1-z}-q^{-1+z}\right)q^{(z+1)x}\,.$$ As a result, we get
the following expressions for
$\Phi_\pm=\psi(x,\lambda)\pm\psi(-x,\lambda)$:
\begin{multline*}
\Phi_\pm(x,\lambda)=\left(q^{1+\lambda}-q^{-1-\lambda}\right)
\left(q^{(\lambda-1)x}\pm q^{(-\lambda+1)x}\right) \\
+\left(q^{1-\lambda}-q^{-1+\lambda}\right)
\left(q^{(\lambda+1)x}\pm q^{(-\lambda-1)x}\right)\,.
\end{multline*}
Taking integer $\lambda\ge 2$, we obtain according to Theorem
\ref{weyl} the following formulas:
\begin{gather*}
P_{\lambda+1}(x;q,q^{-1})=
\left(q^{1-\lambda}-q^{-1+\lambda}\right)^{-1}\Phi_+=\mathfrak
m_{\lambda+1}+ \frac{[1+\lambda]}{[1-\lambda]} \mathfrak
m_{\lambda-1}\,,\\ P_{\lambda-2}(x;q,q^2)=\frac{\Phi_-(x,\lambda)}
{[1-\lambda][x-1][x][x+1]}\,.
\end{gather*}

\end{example}

\bigskip
\begin{remark}
In the case $m=0$ the BA function is a pure exponential,
$\psi=q^{2(x,z)}$. Theorem \ref{weyl} reduces in this case to
Weyl's character formula. Thus, one can think of Theorem
\ref{weyl} as a generalized Weyl formula with $\psi$ being a
"perturbed" exponent. Note that there is a similar result
involving "non-symmetric" Macdonald polynomials (see \cite{M3}).
Our $\psi(x,\lambda)$ for $\lambda\in P$ is polynomial, too, but
it differs from non-symmetric Macdonald polynomials. In fact,
$\psi$ is related to a natural $q$-analogue of the hypergeometric
function by Heckman and Opdam \cite{HO}. However, we will not
discuss this relation here.
\end{remark}

\medskip
\begin{remark} In case $R=A_n$ our formula for $\Phi_-$ is the
formula conjectured by Felder and Varchenko \cite{FV} and proved
by Etingof and Styrkas \cite{ES} in the context of the
representation theory for quantum groups.
\end{remark}

\bigskip
The following proposition is another direct corollary of Theorem
\ref{weyl}.

\medskip
\begin{prop}\label{shact} Let $P_\mu^{(m)}=P_\mu(x;q,q^m)$ denote the Macdonald
polynomial for a root system $R$ and $t=q^m$, where
$m=\{m_\alpha\}$ are positive integers. Consider the shift
operators ${S^\pm_m}$ constructed in Propositions \ref{+} and
\ref{-}. Define $\widehat{S}^\pm_m$ as
$$\widehat{S}^+_m=\delta_m^{-1}\circ{S^+_m}\circ\delta_{m-1}\,,\qquad
\widehat{S}^-_m=\delta_{m-1}^{-1}\circ{S^-_m}\circ\delta_{m}\,.$$
Then these operators act in the following simple way onto
Macdonald polynomials:
$$\widehat{S}^+_mP_{\lambda-\rho_m}^{(m)}=c_m(\lambda)
P_{\lambda-\rho_{m+1}}^{(m+1)}$$ for all
$\lambda\in\rho_{m+1}+P_+$ and
$\widehat{S}^+_mP_{\lambda-\rho_m}^{(m)}=0$ otherwise.  Further,
$$\widehat{S}^-_mP_{\lambda-\rho_{m+1}}^{(m+1)}=
c_m(-\lambda)P_{\lambda-\rho_m}^{(m)}$$ for all
$\lambda\in\rho_{m+1}+P_+$. In these formulas $c_m$ denotes the
function \eqref{*}.
\end{prop}

\bigskip
\begin{proof}
Let $\psi_m(x,z)$ denote the normalized BA function for the system
$R$ with multiplicities $m$. According to Theorem \ref{weyl} we
have: $$\Phi_-^{(m)}(x,\lambda)=\sum_{w\in
W}(-1)^w\psi_m(x,w\lambda)=c(\lambda,m)\delta_m
P_{\lambda-\rho_{m+1}}^{(m+1)}(x)\,.$$ As we know, the shift
operator $S_m^-$ sends $\psi_m(x,w\lambda)$ to
$c_m(\lambda)c_m(-\lambda)\psi_{m-1}(x,w\lambda)$ (here we use
$W$-invariance of $c_m(z)c_m(-z)$). Applying it to $\Phi_-$, we
obtain:
$$S_m^-\Phi_-^{(m)}(x,\lambda)=c_m(\lambda)c_m(-\lambda)\Phi_-^{(m)}\,.$$
After rewriting it in terms of Macdonald polynomials, we get:
$$S_m^-\left[c(\lambda,m)\delta_m
P_{\lambda-\rho_{m+1}}^{(m+1)}\right]=c_m(\lambda)c_m(-\lambda)
c(\lambda,m-1)\delta_{m-1} P_{\lambda-\rho_m}^{(m)}\,,$$ which
leads directly to the last formula from the proposition. Another
part, involving $S_m^+$, can obtained in the same way.
\end{proof}

\bigskip
Later we will apply these results to prove the so-called norm
identity for Macdonald polynomials. In order to do this, we need
one more result.

\medskip
\begin{prop}\label{shconj}
Let $\langle , \rangle_k$ denote the scalar product \eqref{scpr},
depending on the parameters $k_\alpha$ which are supposed to be
{\bf positive} integers. Consider two shift operators $\widehat
{S}^\pm_k$, constructed in the previous proposition. Then
$\widehat S^-_k$ is "adjoint" to $\widehat S^+_k$ in the following
sense: $$\langle \widehat S^+_k f, g
\rangle_{k+1}=(-1)^{|R_+|}q^{2\sum_{\alpha>0}k_\alpha}\langle
f,\widehat S_k^- g \rangle_{k}$$ for any two polynomials $f,\,g$.
\end{prop}

\medskip
\begin{proof}
First, recall that according to \cite{M1} each of the Macdonald
operators $D_k$, being restricted to the space of $W$-invariant
polynomials, is self-adjoint with respect to the scalar product
$\langle , \rangle_k$. We should warn the reader that this is not
true for $k_\alpha\in \mathbb Z_-$. Nevertheless, for all $k$ the
Macdonald operator $D_k$ is {\it formally self-adjoint}. This is a
purely algebraic statement, related to the following rule of
calculating the adjoint:
\begin{equation}\label{adjo1}
(a_\tau T^\tau)^*=\left(\Delta_k
\overline{\mathstrut\Delta}_k\right)^{-1}\circ T^\tau\circ
\left(\Delta_k\overline{\mathstrut\Delta}_k \overline{\mathstrut
a}_\tau\right)\,,
\end{equation}
where $\overline{\mathstrut f}=f(-x)$.
 The point is that this definition leads to the relation $$\langle
a_\tau T^\tau f,g \rangle _k= \langle f, (a_\tau T^\tau)^* g
\rangle_k$$ as soon as we are sure that
$a_\tau\Delta_k\overline{\mathstrut\Delta}_k$ is polynomial.

In a similar way, let us define another type of adjoint by the
rule
\begin{equation}\label{adjo2}
(b_\nu
T^\nu)^*=\left(\Delta_k\overline{\mathstrut\Delta}_k\right)^{-1}\circ
T^\nu\circ
\left(\Delta_{k+1}\overline{\mathstrut\Delta}_{k+1}\overline{\mathstrut
b}_\nu\right)\,.
\end{equation}
 This
definition leads to the relation $$\langle b_\nu T^\nu f,g \rangle
_{k+1}= \langle f, (b_\nu T^\nu)^* g \rangle_k$$ as soon as we
know that $b_\nu\Delta_{k+1}\overline{\mathstrut\Delta}_{k+1}$ is
polynomial.

Now let us consider the shift operator $\widehat S^+_k$ from
Proposition \ref{shact}. As we know, it satisfies the intertwining
relation $$D_{k+1} \widehat S^+_k= \widehat S^+_k D_k\,,$$ where
$D_k$ now stands for the Macdonald operator with $t=q^k$. Taking
the adjoints according to the definition \eqref{adjo2}, we obtain
that $$\left(\widehat S^+_k\right)^* D_{k+1}= D_k \left(\widehat
S^+_k\right)^*$$ (here we used that $D_k$ and $D_{k+1}$ are
formally self-adjoint in the sense of \eqref{adjo1}).

Notice that this coincides with the intertwining relation for the
shift operator $\widehat S^-_k$. One easily compares their leading
terms, using the formulas from Propositions \ref{+} and \ref{-}
and the definition \eqref{adjo2} of the adjoint (in doing this it
is useful to present the function \eqref{deltax} in terms of
\eqref{delta} as $\delta_m=(-1)^{|R_+|}\Delta_{m+1}/\Delta_{-m}$).
A straightforward calculation gives that the leading terms
coincide up to a factor $(-1)^{|R_+|}q^{2\sum k_\alpha}$. Thus,
the intertwiners should also coincide, due to Proposition
\ref{shone}. This proves that $$(\widehat
S^+_k)^*=(-1)^{|R_+|}q^{2\sum k_\alpha}\widehat S^-_k\,,$$ where
the adjoint is understood in the formal sense \eqref{adjo2}.

Now, in order to derive the relation between the scalar products,
it remains to check that if we present the operator $\widehat
S^+_k$ as a sum $\sum b_\nu T^\nu$ then each coefficient $b_\nu$
after multiplying by
$\Delta_{k+1}\overline{\mathstrut\Delta}_{k+1}$ becomes
polynomial. Clearly, it is enough to prove that the operator
$A=S_m^+\circ\delta_{m-1}$ has polynomial coefficients. To prove
this, we consider a quasipolynomial
$\phi(x,z)=q^{2(x,z)}\delta_{m-1}(x)$ which obviously satisfies
the conditions \eqref{axxpsi} for $j=1,\dots, m_\alpha-1$. By
Lemma \ref{prin}, $\phi$ can be presented in the form
$\phi=L_z(\psi_{m-1})$ for some difference operator $L$ in the
$z$-variable, where $\psi_{m-1}$ is the normalized BA function.
Using that $S_m^+\psi_{m-1}=\psi_m$, we see that
$S_m^+\phi=L_z\psi_m$, i.e. it is polynomial in $x$. This proves
that the operator $A=S_m^+\circ\delta_{m-1}$ maps polynomials into
polynomials, thus, it must have polynomial coefficients. This
completes the proof of the proposition.
\end{proof}

\subsection{Macdonald--Cherednik identities}

Our results imply immediately three Macdonald's conjectures about
polynomials $P_\lambda$, which were first proved (for all reduced
root systems) by Cherednik \cite{Ch1,Ch2}.

First, let us derive the norm identity in a traditional way,
following the idea first used by Opdam \cite{O} in case $q=1$. Let
$P_{\mu}^{(k)}$, as above, denote the Macdonald polynomial
$P_{\mu}(x;q,q^k)$ related to the system $R$ and parameters
$k_\alpha\in \mathbb Z_+$, and $\langle , \rangle_k$ stands for
the corresponding scalar product \eqref{scpr}. Now we take an
arbitrary weight $\lambda\in\rho_{m+1}+P_+$ and rewrite the scalar
product $\langle P_{\lambda-\rho_{m+1}}^{(m+1)} ,
P_{\lambda-\rho_{m+1}}^{(m+1)}\rangle_{m+1}$ with the help of
propositions \ref{shact} and \ref{shconj}:
\begin{multline*}
\langle P_{\lambda-\rho_{m+1}}^{(m+1)} ,
P_{\lambda-\rho_{m+1}}^{(m+1)}\rangle_{m+1}=[c_m(\lambda)]^{-2}
\langle \widehat S^+_m P_{\lambda-\rho_{m}}^{(m)} , \widehat S^+_m
P_{\lambda-\rho_{m}}^{(m)}\rangle_{m+1}=\\(-1)^{|R_+|}q^{\sum
m_\alpha}[c_m(\lambda)]^{-2} \langle P_{\lambda-\rho_{m}}^{(m)} ,
\widehat S^-_m\widehat S^+_m
P_{\lambda-\rho_{m}}^{(m)}\rangle_{m}=\\ (-1)^{|R_+|}q^{\sum
m_\alpha}\frac{c_m(-\lambda)}{c_m(\lambda)}\langle
P_{\lambda-\rho_{m}}^{(m)} ,
P_{\lambda-\rho_{m}}^{(m)}\rangle_{m}\,.
\end{multline*}
This reduces the problem of calculating the norms at $t=q^{m+1}$
to the same problem but for $t=q^m$. For instance, in the simplest
case when all $m_\alpha$ are equal, we may descend to $m=0$.
Macdonald polynomials $P_{\lambda}^{(0)}$ are simply the orbitsums
$m_\lambda(x)$, so in this case $\langle P_{\lambda}^{(0)},
P_{\lambda}^{(0)}\rangle_0=|W|$. This gives the norm identity:
$$\langle P_{\lambda-\rho_{k}}^{(k)} ,
P_{\lambda-\rho_{k}}^{(k)}\rangle_{k}=
|W|\,q^{\sum_{\alpha>0}k_\alpha(k_\alpha-1)}\prod_{\alpha\in
R_+}\prod_{j=1}^{k_\alpha-1}\frac
{[(\alpha^\vee,\lambda)+j]}{[(\alpha^\vee,\lambda)-j]}\,.$$ This
result extends easily to the two-orbit case by using the shift
operators which lower the multiplicity $m_\alpha$ at one of the
orbits only.

\bigskip
Now let us return to Theorem \ref{weyl} and apply the duality.
Recall that starting from a root system $R$ and multiplicities
$m_\alpha\in\mathbb Z_+$ we have constructed a function
$$\Phi_-(x,z)=\sum_{w\in W}(-1)^w\psi_m(x,wz)=\sum_{w\in
W}(-1)^w\psi_m(wx,z)\,.$$ We can consider also a similar function
$\Phi^\vee_-$ for the dual root system $R^\vee$. Then Theorem
\ref{symmetry} implies that these two functions are related
through interchanging the arguments:
$$\Phi_-(x,z)=\Phi_-^\vee(z,x)\,.$$ Now take $x=\lambda\in
\rho_{m+1}+P_+$ and $z=\mu\in\rho^\vee_{m+1}+P_+^\vee$ and use
Theorem \ref{weyl}. This gives that
$$c(\lambda,m)\delta_m(\mu)P_{\lambda-\rho_{m+1}}(\mu;q,q^{m+1})=
c^\vee(\mu,m)\delta^\vee_m(\lambda)
P^\vee_{\mu-\rho^\vee_{m+1}}(\lambda;q,q^{m+1})\,,$$ where we used
$\vee$ to distinguish the objects related to the dual root system
$R^\vee$. Taking into account formulas \eqref{cc}, \eqref{deltax}
and their natural counterparts for the system $R^\vee$, we arrive
after simple transformations to the relation
\begin{multline}\label{s1}
P_{\lambda-\rho_{m+1}}(\mu;q,q^{m+1})\prod_{\alpha\in
R_+}\prod_{j=0}^{m_\alpha} [(\alpha^\vee,\lambda)+j]^{-1}
=\\P_{\mu-\rho^\vee_{m+1}}(\lambda;q,q^{m+1})\prod_{\alpha\in
R_+}\prod_{j=0}^{m_\alpha} [(\alpha,\mu)+j]^{-1} \,.
\end{multline}
Put now $\mu=\rho^\vee_{m+1}$, then the right-hand side will not
depend on $\lambda$ (since $P_{0}(x)\equiv 1$). This gives us that
$$P_{\lambda-\rho_{m+1}}(\rho^\vee_{m+1};q,q^{m+1})= {\rm
const}\prod_{\alpha\in R_+}\prod_{j=0}^{m_\alpha}
[(\alpha^\vee,\lambda)+j]\,.$$ To determine the constant, we
substitute $\lambda=\rho_{m+1}$ which leads to the evaluation
identity:
 $$P_{\lambda-\rho_{m+1}}(\rho^\vee_{m+1};q,q^{m+1})=
\prod_{\alpha\in R_+}\prod_{j=0}^{m_\alpha}
\frac{[(\alpha^\vee,\lambda)+j]}
{[(\alpha^\vee,\rho_{m+1})+j]}\,.$$ Denoting $\lambda-\rho_{m+1}$
by $\nu$ and $m_\alpha+1$ by $k_\alpha$, we can rewrite it as
follows:$$P_{\nu}(\rho^\vee_{k};q,q^{k})=q^{-2(\rho_k^\vee,\nu)}
\prod_{\alpha\in R_+}\prod_{j=0}^{k_\alpha-1}
\frac{1-q^{2j+2(\alpha^\vee,\nu+\rho_k)}}{1-q^{2j+2(\alpha^\vee,\rho_{k})}
}\,.$$ This can be presented as $$P_{\nu}(\rho^\vee_{k};q,q^{k})=
q^{-2(\rho_k^\vee,\nu)} \prod_{\alpha\in
R_+}\prod_{j=0}^{(\alpha^\vee,\nu)-1}
\frac{1-q^{2j+2k_\alpha+2(\alpha^\vee,\rho_k)}}
{1-q^{2j+2(\alpha^\vee,\rho_{k})} }\,.$$ An advantage of this form
is that it works for all (e.g. non-integer) $k_\alpha$. Indeed, it
is clearly rational function of $t_\alpha=q^{k_\alpha}$. On the
other hand, it is easy to see that for all $k_\alpha$ the value
$P_\nu(\rho^\vee_k;q,q^k)$ must be rational in $t=q^k$.  Since
this identity is valid for integer $k$, it remains valid for all
$k$.

Using this identity (and its counterpart for the system $R^\vee$)
one can rewrite the relation \eqref{s1} in a more compact form,
known as the symmetry identity:
$$\frac{P_{\lambda-\rho_{m+1}}(\mu;q,q^{m+1})}
{P_{\lambda-\rho_{m+1}}(\rho^\vee_{m+1};q,q^{m+1})}=
\frac{P^\vee_{\mu-\rho^\vee_{m+1}}(\lambda;q,q^{m+1})}
{P^\vee_{\mu-\rho^\vee_{m+1}}(\rho_{m+1};q,q^{m+1})}\,.$$ Changing
notations, we can rewrite it as follows:
$$\frac{P_{\lambda}(\mu+\rho^\vee_k;q,q^{k})}
{P_{\lambda}(\rho^\vee_{k};q,q^{k})}=
\frac{P^\vee_{\mu}(\lambda+\rho_k;q,q^{k})}
{P^\vee_{\mu}(\rho_{k};q,q^{k})}\,\qquad \forall \lambda\in
P_+\,,\quad\forall \mu\in P^\vee_+\,.$$ Again, it is easy to see
that both sides are rational functions of $t=q^k$. Thus, in this
form this identity is valid for all (e.g. non-integer) $k_\alpha$.

 This last
identity is very important, since it leads directly (see \cite
{Ch2}) to the recurrence relations between Macdonald polynomials
$P_\lambda$ with different $\lambda$, which is a higher analogue
of the three-term relation for classical orthogonal polynomials.

\subsection{Limiting case $q=1$}
\label{q=1}

 Let us describe briefly what happens in the limit $q,t\to
1$ if we keep $t_\alpha$ to be of the form $t_\alpha =
q^{-m_\alpha}$ with fixed $m_\alpha$. It is convenient to present
$q$ as $q=e^{\epsilon}$ and rescale the $x$-variable: $x\to
\epsilon^{-1}x$. The normalized BA function constructed in
previous sections is a common eigenfunction of the Macdonald
operators, acting in $x$ and $z$. Let us rewrite these operators
in new notations.

First, in the $x$-variable we have the operators $D^\pi$ related
to (quasi)minuscule co\-weights of the root system $R$ and given
by formulas \eqref{m}--\eqref{m2} with $k_\alpha=-m_\alpha$. To
avoid non-essential details, let us consider minuscule coweights
only. Then in new notations the formula \eqref{m} takes the form:
\begin{equation}
\label{1m} D_{x}^\pi = \sum_{\tau\in W\pi}a_\tau T_{x}
^{\epsilon\tau}\,, \qquad a_\tau (x) =
\prod_{\genfrac{}{}{0pt}{}{\alpha \in R :}{(\alpha, \tau)>0}}\,
\frac {\sinh (-\epsilon m_\alpha+(\alpha, x))} {\sinh (\alpha,
x)}\,.
\end{equation}
Here $T_{x}^{\epsilon\tau}$ denotes the shift in the $x$-variable
by the vector $\epsilon\tau$.

Similarly, in the $z$-variable we have the operators $D_z^\omega$
related to minuscule weights $\omega$ of the system $R$:
\begin{equation}
\label{2m} D_{z}^\omega = \sum_{\tau\in W\omega}a_\tau T_{z}
^{\tau}\,, \qquad a_\tau (z) = \prod_{\genfrac{}{}{0pt}{}{\alpha
\in R :}{(\alpha, \tau)>0}}\, \frac {\sinh
\epsilon(-m_\alpha+(\alpha^\vee, z))} {\sinh \epsilon(\alpha^\vee,
z)}\,.
\end{equation}

The normalized BA function $\psi(x,z)$ related to the root system
$R$ and multiplicities $m_\alpha\in\mathbb Z_+$, is a common
eigenfunction of all these operators. It can be characterized
uniquely in terms of its analytic properties in the $z$-variable,
as it was done above or, alternatively, in terms of its properties
in $x$ (due to duality). Similarly, one can construct it using the
operator $D_z^\omega$ as in Theorem \ref{exist} or, alternatively,
by using in a similar way the operator $D_x^\pi$.

Now let us look what happens in the limit $\epsilon\to 0$. It is
quite clear that the operator \eqref{2m} in this limit takes the
form:
\begin{equation}
\label{02m} D_{z}^\omega = \sum_{\tau\in W\omega}a_\tau T_{z}
^{\tau}\,, \qquad a_\tau (z) = \prod_{\genfrac{}{}{0pt}{}{\alpha
\in R :}{(\alpha, \tau)>0}}\, \frac {(\alpha^\vee, z)-m_\alpha}
{(\alpha^\vee, z)}\,.
\end{equation}
On the other hand, it is slightly more difficult to see what is a
proper limit of the operator \eqref{1m} as $\epsilon$ goes to
zero. In order to do this one should expand $D_x^\pi$ in a series
in $\epsilon$ and pick up the first nontrivial term of this
expansion. It turns out to be the following second-order
differential operator:
\begin{equation}\label{01m}
L_m=\Delta-\sum_{\alpha\in R_+}2
m_\alpha\coth(\alpha,x)\partial_\alpha+4(\rho_m,\rho_m) \,,
\end{equation}
where $\Delta$ is the Laplacian in $V=\mathbb R^n$,\,
$\partial_\alpha$ stands for the derivative in $\alpha$-direction,
and $\rho_m$ is given by \eqref{rho1}. This operator plays the
central role in Heckman--Opdam theory of multivariable
hypergeometric functions \cite{HO}, and it is gauge-equivalent to
the generalized Calogero--Sutherland operator from \cite{OP}.

{From} this it is natural to expect that a proper limit of the BA
function $\psi(x,z)$ must give a common eigenfunction both for the
rational Macdonald operators \eqref{02m} and the operator $L_m$.
However, it is quite difficult to see this directly, looking at
our formula for $\psi$. The best way is to repeat the main
constructions independently for this degenerate case. Since now we
have no symmetry between $x$ and $z$ variables, we may choose two
different ways to describe $\psi$: either in terms of
$x$-properties, or in terms of its properties in the $z$-variable.
The $x$-part of the story is very similar to what we had before:
$\psi$ has the form $$\psi(x,z)=e^{2(x,z)}\sum_{\nu\in P}\psi_\nu
e^{2(\nu,x)}$$ with the summation taken over the polytope
\eqref{nuvee}. The main difference is that the conditions
\eqref{axxpsi} are replaced by the following: for each $\alpha\in
R_+$ and $s=1,\dots, m_\alpha$
\begin{equation}\label{0axxpsi}
(\partial_\alpha)^{2s-1}\psi\equiv 0\quad\text{for}\ \
q^{2(\alpha,x)}=1\,.
\end{equation}
Here $\partial_\alpha$ denotes the derivative in
$\alpha$-direction in the $x$-variable. One can check that the
operator \eqref{01m} with $k_\alpha=-m_\alpha$ preserves these
properties, which gives an analog of Proposition \ref{inv}. After
that everything becomes more or less a straightforward
modification of the previous constructions.

The $z$-properties of $\psi$ in this case change more
significantly: it becomes qua\-si\-poly\-nomial in $z$ in a
standard sense: $$\psi(x,z)=P(x,z)\exp 2(x,z)\,,\qquad
P=C(x)\prod_{\alpha\in R_+}(\alpha^\vee,z)^{m_\alpha}+\dots\,, $$
where $C(x)=\prod_{\alpha\in R_+}(\sinh(\alpha,x))^{m_\alpha}$ and
the dots stand for lower terms in $z$. A detailed exposition of
this latter approach is given in \cite{C00}, where the relation of
the $\psi$-function to the hypergeometric function and Jacobi
polynomials by Heckman and Opdam is also discussed in detail.

To illustrate the difference between these two approaches, let us
compare how the formula \eqref{ber} will look for each of them.
First, one can construct $\psi$ using $L=L_m$. Let us introduce a
trigonometric polynomial $Q(x)$ as
$$Q(x)=e^{2(\rho,x)}\prod_{\alpha\in R_+}
(2\sinh(\alpha,x))^{2m_\alpha}\,, $$ where $\rho=\rho_m$ is given
by \eqref{rho1}. Then $\psi$, up to a certain $z$-depending
factor, is given by the formula $$\psi\sim\prod_\nu
(L-4(z+\nu)^2)[Q(x)e^{2(x,z)}]\,,$$ where the product is taken
over all $\nu\ne 0$ having the form
$\nu=\sum_{\alpha>0}l_\alpha\alpha$ with $l_\alpha=0,1,\dots,
m_\alpha$.

Alternatively, to construct $\psi$ one can use $D=D_z^\omega$
given by \eqref{02m}. Then, up to a certain $x$-depending factor,
$\psi$ is given by the expression (see \cite{C00})
$$\psi\sim(D-m_\omega(x))^M[q(z)e^{2(x,z)}]\,,\qquad
M=\sum_{\alpha\in R_+}m_\alpha\,,$$ where
$m_\omega(x)=\sum_{\tau\in W\omega}e^{2(\tau,x)}$ and the
polynomial $q(z)$ looks as follows: $$q(z)=\prod_{\alpha\in
R_+}\prod_{j=1}^{m_\alpha}((\alpha^\vee,z)^2-j^2)\,.$$

The first expression for $\psi$ can be viewed as a trigonometric
generalization of the original formula by Berest \cite{Be}, while
the second one is its difference version. As we see, they both are
specializations of the same formula \eqref{ber}.

\section{$BC_n$ case and Koornwinder polynomials}

In this section we consider the case of the non-reduced root
system $R=BC_n$. A proper generalization of the Macdonald theory
in this case was proposed by Koornwinder \cite{Ko}. It depends on
five parameters (apart from $q$) and generalizes Askey--Wilson
polynomials \cite{AW} to higher dimensions. The one-dimensional
case will be essential for us, so we start considering $R=BC_1$.

\subsection{Rank-one case}
In this case we have $4$ parameters $a,b,c,d$ (apart from $q$),
and the corresponding one-dimensional difference operator $D$,
suggested by Askey and Wilson, looks as follows \cite{AW}:
\begin{equation}\label{aw}
D=v^+(z)(T-1)+v^-(z)(T^{-1}-1)+q(abcd)^{-1}+q^{-1}abcd\,,
\end{equation}
where $T^{\pm 1}$ denotes the shift by $\pm 1$ in $z$ and the
coefficients $v^{\pm}$ are given by the following formula:
\begin{align}\label{a+}
v^+(z)=&\frac{(aq^{z}-a^{-1}q^{-z})(bq^{z}-b^{-1}q^{-z})(cq^{z}+c^{-1}q^{-z})
(dq^{z}+d^{-1}q^{-z})}{(q^z-q^{-z})
(q^{\frac12+z}-q^{-\frac12-z})(q^z+q^{-z})
(q^{\frac12+z}+q^{-\frac12-z})}\,,\\\label{a-}v^-(z)=&v^+(-z)\,.
\end{align}
Our notations differ from those of Askey and Wilson \cite{AW}:
what they denote by $(q,a,b,c,d)$ is $(q^2,a^2,b^2,-c^2,-d^2)$ in
our notations.

The following function $\Delta(z)=\Delta(z;a,b,c,d,q)$ plays an
important role in Askey--Wilson's theory:
\begin{equation}\label{awdelta}
\Delta(z)=(abcd/q)^{-z}\frac{(q^{4z},q^{2+4z};q^4)_\infty}
{(a^2q^{2z},b^2q^{2z},-c^2q^{2z},-d^2q^{2z};q^2)_\infty}\,.
\end{equation}
Here we used the standard notations: $$(t;q)_\infty:=\prod_{i\ge
0}(1-tq^i)\,,\qquad (t_1,\dots,t_n;q)_\infty:=\prod_{s=1}^n
(t_s;q)_\infty\,.$$
 Using it, one can present the coefficients $v^\pm$ of the
Askey--Wilson operator $D$ as
\begin{equation}\label{awcoef}
v^+(z)=\Delta(z+1)/\Delta(z)\,,\qquad
v^-(z)=\Delta(-z-1)/\Delta(-z)\,.
\end{equation}
Introduce the {\it dual parameters}  $\widetilde a,\widetilde b,
\widetilde c, \widetilde d$ as follows:
\begin{equation}\label{duabcd}
\widetilde a=(abcd/q)^{\frac12}\,,\quad \widetilde
b=(abq/cd)^{\frac12}\,,\quad \widetilde
c=(acq/bd)^{\frac12}\,,\quad \widetilde d=(adq/bc)^{\frac12}\,.
\end{equation}
We will denote by $\widetilde \Delta$ the function \eqref{awdelta}
with the dual parameters:
\begin{equation}\label{awdudelta}
\widetilde \Delta(z)=\Delta(z;\widetilde a,\widetilde b,\widetilde
c,\widetilde d,q)\,.
\end{equation}

 Now let us make some special choice of parameters $a,b,c,d$ in
\eqref{a+}. Namely, we put
\begin{equation}\label{integers}
a=q^{-l},\quad b=q^{-l'},\quad c=q^{-m},\quad d=q^{-m'}\,,
\end{equation}
where $l,l',m,m'\in \frac12\mathbb Z$ are some (half)integers with
the requirement that
\begin{equation}\label{req}
\frac12+l+l'\in \mathbb
Z\qquad\text{and}\qquad\frac12+m+m'\in\mathbb Z\,.
\end{equation}
We will assume that $l,l',m,m'$ are positive. Introduce $N$ as
\begin{equation}\label{M}
N=1+l+l'+m+m'\in\mathbb Z\,.
\end{equation}
Using it, one can present the constant term in \eqref{aw} as $$
q(abcd)^{-1}+q^{-1}abcd=q^N+q^{-N}\,.$$ Below we will denote by
$(\widetilde l, \widetilde l', \widetilde m, \widetilde m')$ the
dual set of multiplicities, determined in accordance with
\eqref{duabcd}:
\begin{equation}\label{dupa}
\begin{pmatrix}
\widetilde l\\ \widetilde l'\\ \widetilde m \\\widetilde m'
\end{pmatrix}
=
\begin{pmatrix}
1/2\\-1/2\\-1/2 \\-1/2
\end{pmatrix}
+\frac12
\begin{pmatrix}1&\phantom{-}1&\phantom{-}1&\phantom{-}1\\
1&\phantom{-}1&-1&-1\\1 &-1&\phantom{-}1&-1\\1&-1&-1&\phantom{-}1
\end{pmatrix}
\begin{pmatrix}
l\\ l'\\ m \\ m'
\end{pmatrix}
\,.
\end{equation}
Notice that this transformation, as well as \eqref{duabcd}, is
{\it involutive}. Geometrically, it reduces to the orthogonal
reflection with respect to the hyperplane $l-l'-m-m'=1$ in
$\mathbb R^4$.

 Now we are going to formulate an analogue of Proposition
\ref{inv} for the Askey--Wilson operator. First, introduce the
following shorthand notation $u\preccurlyeq (v,w)$ in the
situation when for $3$ given numbers $u,v,w$ at least one of the
differences $v-u, w-u$ belongs to $\mathbb Z_{\ge 0}$. Now let us
consider a ring $\mathfrak R$ which consists of all polynomials
$$f(z)=\sum_{j\in\mathbb Z}f_jq^{jz}\,,$$ satisfying the following
$N=1+l+l'+m+m'$ conditions:
\begin{align}\label{awaxz1}
&f(z-s)=f(z+s)&\qquad\text{ for}&&\   q^{2z}=1&\quad\text{ and}\
&0<s\preccurlyeq(l,l')\,,\\\label{awaxz4}&
f(z-s)=f(z+s)&\qquad\text{ for}&&\ q^{2z}=-1&\quad\text{ and}\
&0<s\preccurlyeq(m,m')\,.
\end{align}

\bigskip
\begin{prop}\label{awinv}
The Askey--Wilson operator \eqref{aw}--\eqref{a-} with the
parameters of the form \eqref{integers}--\eqref{req} preserves the
ring $\mathfrak R$ above:\ $D(\mathfrak R)\subseteq \mathfrak R$.
\end{prop}

\medskip
Proof can be found along the lines of our proof of Proposition
\ref{per}. We need also the following inverse result.

\medskip
\begin{lemma}\label{awinverse}Let $D$ be the Askey--Wilson operator
\eqref{aw} with $a,b,c,d$ as in \eqref{integers}--\eqref{req}. Let
$f$ be an eigenfunction  of $D$ of the form $f=r(x,z)q^{2xz}$
where $r$ is rational in $z$. Then $f$ is, in fact, non-singular
and satisfes conditions \eqref{awaxz1}--\eqref{awaxz4}.
\end{lemma}

\medskip
\begin{proof}
First, it is easy to see from equation $Df=\lambda f$ that in the
case \eqref{integers}--\eqref{req} any eigenfunction either has
infinite number of poles, or has no poles at all. This proves that
$f$ is nonsigular in $z$. Then the conditions
\eqref{awaxz1}--\eqref{awaxz4} can be derived similar to lemma
\ref{inverse1}.
\end{proof}

\bigskip
All this suggests the following definition of a Baker--Akhiezer
function for the Askey--Wilson operator with multiplicities
$(l,l',m,m')$.

\begin{definition} A function $\psi(x,z)$ of the form
\begin{equation}\label{awpsi}
\psi=q^{2xz}\sum_{|j|\le N} \psi_{j}q^{jz}\,,\qquad
\psi_j=\psi_j(x)\,,
\end{equation}
is called a Baker--Akhiezer function for the Askey--Wilson
operator if it satisfies $N=1+l+l'+m+m'$ conditions
\eqref{awaxz1}--\eqref{awaxz4} in $z$ for all $x$.
\end{definition}

\bigskip In the same way as in case $R=A_1$, one proves that such a
$\psi$ does exist and it is unique up to an $x$-depending factor.
Also it is quite clear that the only nonzero coefficients $\psi_j$
in \eqref{awpsi} will be those with $j+N\in 2\mathbb Z$. However,
we were not able to calculate $\psi_j$ explicitly. This makes the
following result somewhat less trivial.

\medskip
\begin{prop}\label{nontr} Suppose that the multiplicities
$l,l',m,m'\in\frac12\mathbb Z_+$ satisfy \eqref{req} and are such
that the dual multiplicities \eqref{dupa} are positive, too. Then
if one puts $\psi_{N}=\widetilde \Delta^{-1}(x)$ in accordance
with the formulas \eqref{awdelta}--\eqref{integers}, then all the
coefficients $\psi_j$ will be polynomial in $x$ with
$\psi_{-N}(x)=\psi_N(-x)$.
\end{prop}

\medskip
As we already said, we cannot prove this proposition directly. So
we use the following strategy: we will construct an eigenfunction
$\psi$ of the Askey--Wilson operator which has the form
\eqref{awpsi} and satisfies the requirements of the proposition.
The constructed $\psi$ will be automatically a BA function due to
lemma \ref{awinverse}.

To construct eigenfunctions of the operator \eqref{aw} for the
case \eqref{integers}--\eqref{req}, we will use shift operators,
which in rank-one case can be computed directly. We need the
operators which shift the parameters $a,b,c,d$ in  \eqref{aw}. We
say that an operator $S$ shifts from $(a,b,c,d)$ to $(\hat
{\mathstrut a},\hat {\mathstrut b},\hat {\mathstrut c},\hat
{\mathstrut d})$ if the following intertwining relation holds:
\begin{equation}\label{intrtw}
\widehat D\circ S=S\circ D\,,
\end{equation}
where $D$ and $\widehat D$ are the Askey--Wilson operators
\eqref{aw} which correspond to $(a,b,c,d)$ and $(\hat {\mathstrut
a},\hat {\mathstrut b},\hat {\mathstrut c},\hat {\mathstrut d})$,
respectively.

\medskip
\begin{prop} Define $4$ operators $S_1,S_2,S_3,S_4$ as
$$S_i=A^+_i(z)T^{1/2}+A^-_i(z)T^{-1/2}\,,\quad i=1,\dots,4\,,$$
where $A^-_i(z)=A^+_i(-z)$ and $A^+_i$ are given by the formulas:
\begin{align*}
A^+_1=&\frac{(aq^{z-\frac12}-a^{-1}q^{\frac12-z})
(bq^{z-\frac12}-b^{-1}q^{\frac12-z})
(cq^{z-\frac12}+c^{-1}q^{\frac12-z})
(dq^{z-\frac12}+d^{-1}q^{\frac12-z}) }{(q^z-q^{-z})
(q^z+q^{-z})}\\ A^+_2=&\frac{(aq^{z-\frac12}-a^{-1}q^{\frac12-z})
(bq^{z-\frac12}-b^{-1}q^{\frac12-z})}{(q^z-q^{-z})
(q^z+q^{-z})}\,,\\
A^+_3=&\frac{(aq^{z-\frac12}-a^{-1}q^{\frac12-z})
(cq^{z-\frac12}+c^{-1}q^{\frac12-z})}{(q^z-q^{-z})
(q^z+q^{-z})}\,,\\
A^+_4=&\frac{(aq^{z-\frac12}-a^{-1}q^{\frac12-z})
(dq^{z-\frac12}+d^{-1}q^{\frac12-z}) }{(q^z-q^{-z})
(q^z+q^{-z})}\,.
\end{align*}
Then $S_i$ are the shift operators for the Askey--Wilson operator
and the corresponding shifts of the parameters are as follows:

$S_1$:\quad $(a,b,c,d) \to
(q^{-\frac12}a,q^{-\frac12}b,q^{-\frac12}c,q^{-\frac12}d)$\,,

$S_2$:\quad $(a,b,c,d) \to
(q^{\frac12}a,q^{\frac12}b,q^{-\frac12}c,q^{-\frac12}d)$\,,

$S_3$:\quad $(a,b,c,d) \to
(q^{\frac12}a,q^{-\frac12}b,q^{\frac12}c,q^{-\frac12}d)$\,,

$S_4$:\quad $(a,b,c,d) \to
(q^{\frac12}a,q^{-\frac12}b,q^{-\frac12}c,q^{\frac12}d)$\,.
\end{prop}

\medskip
\begin{proof}
It reduces to a straightforward though tedious calculation.
Because of a certain symmetry between $a,b,c,d$ it suffices to
check the intertwining relation for $S_1$ and $S_2$ only.
\end{proof}

\medskip
\begin{proof}[Proof of Proposition \ref{nontr}]
Recall that we assume that both $l,l',m,m'$ and their dual
\eqref{dupa} are positive. Consider the case when $\widetilde l$
and $\widetilde m$ are both integers. Then there is a proper
composition of shifts $S_i$ which shifts from
$(1,q^{1/2},1,q^{1/2})$ to $(q^{-l},q^{-l'},q^{-m},q^{-m'})$.
Namely, one should apply $S_1$\ $\widetilde l$ times, $S_2$ \
$(\widetilde l'+1/2)$ times, $S_3$ \ $\widetilde m$ times and
$S_4$ \ $(\widetilde m'+1/2)$ times. The Askey--Wilson operator
with the parameters $(1,q^{1/2},1,q^{1/2})$ is a trivial one,
$D_0=T+T^{-1}$, with $\psi_0=q^{2xz}$ being its eigenfunction.
Applying the composition of the shifts above to $\psi_0$, we
obtain an eigenfunction $\psi$ of the operator $D$. It is
obviously quasipolynomial in $x$. Moreover, it will be
quasipolynomial in $z$, too (this is not completely obvious, since
the shifts have singularities, but follows, for instance, from
lemma \ref{awinverse}). Now in order to find the leading
coefficients, it is sufficient to look for each application of
$S_i$ at the asymptotics of $\psi$ in $z$ at $\pm\infty$, using
the explicit formulas for $S_i$. After some simple inductive
calculations, we arrive directly at the formula for $\psi_{\pm N}$
from the proposition.

The three other possible cases, namely

\noindent (1) $l,m'\in\mathbb Z$\,,\qquad (2) $l',m\in\mathbb
Z$\,,\qquad (3) $l',m'\in\mathbb Z$\,,

\noindent can be considered in a similar manner. The only
difference is that in cases 2 and 3 one applies shifts starting
from $D$ with $(a,b,c,d)$ being $(1,q^{1/2},q^{1/2},q)$ or
$(1,q^{1/2},q,q^{1/2})$, respectively. These operators are also
almost trivial: they are obtained from $D_0=T+T^{-1}$ by a simple
gauge. For instance, the Askey--Wilson operator with parameters
$(1,q^{1/2},q^{1/2},q)$ has the form
\begin{multline*}
D=\frac{q^{1+x}+q^{-1-x}}{q^x+q^{-x}}(T-1)+
\frac{q^{1-x}+q^{-1+x}}{q^{-x}+q^{x}}(T-1)+(q+q^{-1})\\=
(q^x+q^{-x})^{-1}\circ (T+T^{-1})\circ (q^x+q^{-x})\,.
\end{multline*}
Then we apply $S_1$\  $(\widetilde l+1/2)$ times, $S_2$ \
$\widetilde l'$ times, $S_3$ \  $\widetilde m$ times and $S_4$ \
$(\widetilde m'+1/2)$ times, arriving at $D$ with the parameters
$(q^l,q^{l'},q^m,q^{m'})$. Other arguments remain the same. The
case 3 is analogous.
\end{proof}

\bigskip
This leads us directly to the main result of this section.

\begin{theorem}\label{awmain}
Let $\Psi(x,z)$ denote a BA function with the parameters
$l,l',m,m'\in\mathbb Z_+$ normalized as in Proposition
\ref{nontr}. Then $\Psi$ solves the following bispectral system:
\begin{equation}
\label{awbs} \left\{
\begin{array}{l}
D_z\Psi = (q^{2x}+q^{-2x})\Psi\,, \\ \widetilde D_x\psi =
(q^{2z}+q^{-2z})\Psi\,.
\end{array}
\right.
\end{equation}
Here $D_z$ is the Askey--Wilson operator \eqref{aw} related to
$l,l',m,m'$, and $\widetilde D_x$ acts in $x$ and is related to
the dual multiplicities \eqref{dupa}. Moreover, if the dual
parameters are positive, we will have the duality as follows:
$\Psi(z,x)=\widetilde \Psi(x,z)$ where $\widetilde \Psi(x,z)$
denotes the normalized BA function related to the dual parameters.
\end{theorem}

\medskip
\begin{proof}
The first equation (in $z$) follows from Proposition \ref{awinv}
and the uniqueness of $\psi$. It does not depend on our particular
way of normalizing $\psi$. To obtain the second equation we apply
the standard argument: consider the function
$\phi=(q^{2z}+q^{-2z})\Psi$. It is quasipolynomial in $z$,
satisfying conditions \eqref{awaxz1}--\eqref{awaxz4}. Then it must
be obtained from $\Psi$ by applying a proper difference operator:
$$L\Psi=(q^{2z}+q^{-2z})\Psi\,,\quad
L=b_+(x)T_x+b_-(x)T_x^{-1}+b_0(x)\,.$$ The coefficients $b_\pm$
can be found easily as soon as we know the leading coefficients
$\Psi_{\pm N}$ in $\Psi$. To calculate $b_0$, however, we need one
more term, say, $\Psi_{-N+2}$ in $\Psi$. Unfortunately, we have
not found anything better than to calculate it directly from the
difference equation $D_z\Psi=(q^{2x}+q^{-2x})\Psi$. This is pretty
straightforward and we shall not reproduce this calculation here.
As a result, one finds that the the operator $L$ is nothing but
the dual Askey--Wilson operator $\widetilde D_x$.

Finally, the duality between $x$ and $z$ follows similar to the
case $R=A_1$, see Proposition \ref{r1}.
\end{proof}

\bigskip

\subsection{BA function for Koornwinder operator}

Now let us consider the difference operator by Koornwinder which
generalizes the Askey--Wilson operator to higher dimensions. This
operator $D$ depends on five parameters $a,b,c,d,t$ apart from $q$
and it looks as follows \cite{Ko1}:
\begin{equation}\label{ko}
D=\sum_{i=1}^n
v_i^+(T_i-1)+v_i^-(T_i^{-1}-1)+abcdq^{-1}\frac{1-t^{2n}}{1-t^2}+
(abcd)^{-1}q\frac{1-t^{-2n}}{1-t^{-2}}\,,
\end{equation}
where $T_i^s$ stands for a shift by $s$ in $z_i$ and the
coefficients $v_i^\pm(z_1,\dots,z_n)$ are given by the formulas:
\begin{equation}\label{kocoef}
v_i^\pm(z)=v^\pm(z_i)\prod_{j\ne i}\frac{\left(tq^{\pm
z_i-z_j}-t^{-1}q^{\mp z_i+z_j}\right) \left(tq^{\pm
z_i+z_j}-t^{-1}q^{\mp z_i-z_j}\right)} {\left(q^{\pm
z_i-z_j}-q^{\mp z_i+z_j}\right) \left(q^{\pm z_i+z_j}-q^{\mp
z_i-z_j}\right)}\,,
\end{equation}
with the functions $v^\pm(z_i)$ obtained by substituting $z=z_i$
into the formulas \eqref{a+}--\eqref{a-}.

The underlying geometrical structure here is an affine root system
$C^\vee C_n$ in notations of \cite{M5}. For our purposes, however,
it will be enough to consider a usual root system $R$ of
$B_n$-type: $$R=\{\pm e_i\}\cup \{\pm e_i \pm e_j \mid i\ne
j\}\,.$$ We fix its positive half $R_+$ as
$$R_+=\{e_1,\dots,e_n\}\cup\{e_i\pm e_j\mid i<j\}\,.$$ The root
lattice is $Q=\mathbb Z^n$ and its positive part $Q_+$ is defined
as: $$Q_+=\{\nu\in\mathbb Z^n\mid \sum_{i=1}^j \nu_i\ge 0\,,\quad
j=1,\dots,n\}\,.$$ Weight lattice $P$ in this case is also the
standard lattice $\mathbb Z^n$:
$$P=\{\nu=\nu_1e_1+\dots+\nu_ne_n\mid \nu_i\in \mathbb Z\}\,,$$
while the cone of the dominant weights $P_+$ looks as follows:
$$P_+=\{\nu\in P\mid \nu_1\ge \nu_2\ge \dots \ge \nu_n\ge 0\}\,.$$
The Weyl group $W$ acts by permuting the variables and flipping
their signs arbitrarily.

As before, by a polynomial we mean a finite sum $f(z)=\sum_{\nu\in
P}q^{2(\nu,z)}$, keeping calling functions like $f(z)q^{2(x,z)}$
{\it quasipolynomial} in $z$. Algebra of $W$-invariant polynomials
in this case is a linear span of orbitsums
\begin{equation}\label{bcorb}
\mathfrak m_\lambda(z)=\sum_{\tau\in W\lambda}q^{2(\tau,z)}\,,
\end{equation}
and it is a symmetric polynomial algebra of the generators
$y_i=q^{2z_i}+q^{-2z_i}$.

Now let us specialize the parameters $a,b,c,d,t$ as follows:
\begin{equation}\label{bcintegers}
t=q^{-k},\quad a=q^{-l},\quad b=q^{-l'},\quad c=q^{-m},\quad
d=q^{-m'}\,,
\end{equation}
where $k\in\mathbb Z_+$ and $l,l',m,m'\in \frac12\mathbb Z_+$ are
some (half)integers with the requirement as in rank-one case that
\begin{equation}\label{bcreq}
\frac12+l+l'\in \mathbb
Z\qquad\text{and}\qquad\frac12+m+m'\in\mathbb Z\,.
\end{equation}
We will denote by $M$ the whole set $$M=(k,l,l',m,m')\,.$$Below we
will also use the multiplicities $m_\alpha$ defined in the
following way:
\begin{align}\label{multi}
&m_\alpha=k\quad\text{for}\quad\alpha=\pm e_i\pm
e_j\qquad\text{and}\\\label{multip}
&m_\alpha=1+l+l'+m+m'\quad\text{for}\quad \alpha=\pm e_i\,.
\end{align}
Introduce a vector $\rho$ depending on $M=(k,l,l',m,m')$ as
follows:
\begin{equation}\label{bcrho}
\rho=\rho_M=\frac12\sum_{\alpha\in R_+}m_\alpha \alpha\,,
\end{equation}
Using it, one can present the constant term in \eqref{ko} in the
following way: $$abcdq^{-1}\frac{1-t^{2n}}{1-t^2}+
(abcd)^{-1}q\frac{1-t^{-2n}}{1-t^{-2}}=\sum_i q^{2(\rho,\pm
e_i)}\,,$$ which makes it similar to the constant term in
\eqref{m1}.

Let $\mathfrak R$ denote a ring which consists of all polynomials
$f(z)=\sum_{\nu\in\mathbb Z^n}q^{2(\nu,z)}$ with the following
properties:

 (1) for all $i=1,\dots,n$ and $0<s\preccurlyeq(l,l')$
\begin{equation}\label{bcaxz1}
(T^s_i-T^{-s}_i)f=0\qquad\text{ for}\ q^{2z_i}=1\,;
\end{equation}

 (2) for all $i=1,\dots,n$ and $0<s\preccurlyeq(m,m')$
\begin{equation}\label{bcaxz2}
(T^s_i-T^{-s}_i)f=0\qquad\text{ for}\ q^{2z_i}=-1\,;
\end{equation}

(3) for all $1\le i<j\le n$ and $s=1,\dots,k$
\begin{equation}
\label{bcaxz3} (T^s_i-T^{s}_j)f=0\qquad\text{ for}\
q^{2z_i}=q^{2z_j}\,;
\end{equation}

(4) for all $1\le i<j\le n$ and $s=1,\dots,k$
\begin{equation}
\label{bcaxz4} (T^s_i-T^{-s}_j)f=0\qquad\text{ for}\
q^{2z_i}=q^{-2z_j}\,.
\end{equation}
We used the same notation $s\preccurlyeq(l,l')$ as above, denoting
that at least one of the differences $l-s,l'-s$ is a nonnegative
integer.

\bigskip
\begin{prop}\label{bcinv}
The Koornwinder operator \eqref{ko}--\eqref{kocoef} with the
parameters  \eqref{bcintegers}--\eqref{bcreq} preserves the ring
$\mathfrak R$ as above:\ $D(\mathfrak R)\subseteq \mathfrak R$.
\end{prop}

Moreover, in a similar way one can check that all $n$ commuting
difference operators $D_1,\dots, D_n$ constructed in \cite{vD1},
will preserve the ring $\mathfrak R$ in case
\eqref{bcintegers}--\eqref{bcreq}.

\bigskip
Now a Baker--Akhiezer function $\psi(x,z)$ is defined similarly to
the case of a reduced root system:

\noindent (1) $\psi$ has the form
\begin{equation}\label{bcpsi}
\psi=q^{2(x,z)}\sum_{\nu\in \mathcal N}\psi_\nu(x)q^{2(\nu,z)}\,,
\end{equation}
where the summation is taken over $\nu\in P$ lying inside the
polytope
\begin{equation}\label{bcnu}
\mathcal N = \{ \nu= \sum_{\alpha \in R_+} l_\alpha \alpha^\vee
\mid -\frac12 m_\alpha\le l_\alpha \le \frac12 m_\alpha \}\,;
\end{equation}

\noindent (2) $\psi(x,z)$ satisfies the conditions
\eqref{bcaxz1}--\eqref{bcaxz4} in $z$ for all $x$\,.

One proves, similar to Section \ref{baf}, that such $\psi$ does
exist and is unique up to an $x$-depending factor. It can be
expressed by a formula, similar to \eqref{ber}. Namely, introduce
first a function $c(\lambda,M)$ depending on $\lambda\in\mathbb
R^n$ and the multiplicities $M=(k,l,l',m,m')$ in the following
way:
\begin{multline}\label{bcc}
c(\lambda,M)=\prod_{\genfrac{}{}{0pt}{}{\alpha\in
R_+^0}{0<j\preccurlyeq(l,l')}}
[j-(\alpha,\lambda)]
\prod_{\genfrac{}{}{0pt}{}{\alpha\in
R_+^0}{0<j\preccurlyeq(m,m')}}
[j-(\alpha,\lambda)]
\prod_{\genfrac{}{}{0pt}{}{\alpha\in
R_+^1}{j=1,\dots, k}}
[j-(\alpha,\lambda)]
\,,
\end{multline}
where $R^0$ and $R^1$ denote the sets of short and long roots,
respectively, and $[a]$ stands, as usual, for $q^a-q^{-a}$.

Now define the polynomial $Q(z)$ as follows:
\begin{equation}\label{bcQ}
Q(z)=(-1)^{n(l+l'+\frac12)}q^{2(\rho,z)}c(z,M)c(-z,M)\,,
\end{equation}
where $c(z,M)$ is defined above and the vector
 $\rho$ is given by \eqref{bcrho}.
 Introduce also the
notation $\mathfrak m$ for the orbitsum
\begin{equation}\label{bcmu}
\mathfrak m(x)=\sum_{i=1}^n(q^{2x_i}+q^{-2x_i})\,.
\end{equation}

\bigskip
\begin{theorem}
Let $D$ be the Koornwinder operator \eqref{ko} with the parameters
as in \eqref{bcintegers}-\eqref{bcreq}. Define $\psi(x,z)$ as
follows:
\begin{equation}\label{bcber}
\psi=\prod_{\nu}\left(D-\mathfrak
m(x+\nu)\right)\left[q^{2(x,z)}Q(z) \right]\,,
\end{equation}
in accordance with the formulas \eqref{bcQ}, \eqref{bcmu}, where
the product is taken over all $\nu\ne 0$ having the form
$\nu=\sum_{\alpha\in R_+}l_\alpha\alpha^\vee$ with
$l_\alpha=0,\dots , m_\alpha$. Then

\noindent(i) $\psi$ is a BA function for the Koornwinder operator;

\noindent(ii) the coefficient $\psi_{-\rho}$ in its expansion
\eqref{bcpsi} equals $$\prod_\nu (\mathfrak m(x)-\mathfrak
m(x+\nu))\ne 0\,;$$

\noindent(iii) as a function of $z$,\ $\psi$ is an eigenfunction
of the Koornwinder operator $D$: $D\psi=\mathfrak m(x)\psi$.
\end{theorem}

\bigskip
To normalize $\psi$, consider the {\bf dual parameters}
$$\widetilde M=(\widetilde k,\widetilde l,\widetilde l',\widetilde
m,\widetilde m')\,,$$ where $\widetilde k=k$ while other $4$
parameters transform according to \eqref{dupa}. Let us normalize
$\psi$ in the following way:
\begin{equation}
\label{bcnorm0} \psi_{\rho}(x)=c(x,\widetilde M)
\end{equation}
in accordance with the formula \eqref{bcc}.

\bigskip
\begin{theorem}[Duality]
Let $\psi(x,z)$ be a BA function related to the parameters
$l,l',m,m',k$ and normalized as above. Suppose that the dual
parameters $\widetilde l,\widetilde l',\widetilde m,\widetilde
m',\widetilde k$ are positive. Then $\psi(x,z)$ is quasipolynomial
in both $x$ and $z$ and has the following duality property:
$$\psi(z,x)=\widetilde\psi(x,z)\,,$$ where $\widetilde\psi$
denotes the normalized BA function related to the dual parameters.
In particular, $\psi(x,z)$, as a function of $z$, is an
eigenfunction of the Koornwinder operator \eqref{ko}, while in $x$
it satisfies a similar difference equation related to the dual
parameters.
\end{theorem}

\bigskip
The algebraic integrability of the Koornwinder operator in case
\eqref{bcintegers} and the existence of the shift operators is a
straightforward generalization of the similar results for reduced
root systems.

\subsection{Koornwinder polynomials}

The Koornwinder polynomials
$$P_\lambda(z)=P_\lambda(z;q,a,b,c,d,t)$$ can be defined similar
to Macdonald ones, as polynomial eigenfunctions of the Koornwinder
operator \eqref{ko}, see \cite{Ko,vD1}. Since our notations are
slightly different, we reproduce here their definition for the
reader's convenience. Namely, $P_\lambda$ has the form
\begin{equation}\label{kop}
P_\lambda=\mathfrak m_\lambda+\sum_{\nu \prec \lambda}a_{\lambda\nu}
\mathfrak m_\nu\,,\quad \lambda\in P_+\,,
\end{equation}
in notations of the previous section, with $\nu \prec \lambda$
meaning that $\lambda-\nu\in Q_+$. For generic values of the
parameters $a,b,c,d,t$ the polynomial $P_\lambda$ is uniquely
determined from the equation:
\begin{equation}\label{koeigen}
D P_\lambda=c_{\lambda \lambda} P_\lambda\,,
\end{equation}
where $D$ is the Koornwinder operator \eqref{ko} and the
eigenvalue $c_{\lambda\lambda}$ has the form
\begin{equation}\label{koeig}
c_{\lambda\lambda}=\sum_{\tau=\pm e_i} q^{2(\tau,\lambda+\rho)}\,,
\end{equation}
with $\rho=\rho_M$ given by formulas \eqref{bcintegers} and
\eqref{multi}--\eqref{bcrho}.

To compare with the notations in \cite{vD1}, one should put
$q^2=e^{i\omega}$, then our orbitsums $\mathfrak m_\lambda(z)$
correspond to $m_\lambda(\omega z)$ in notations of \cite{vD1}.
Our defining equation \eqref{koeigen} corresponds to the equation
$$\widehat D_1 p_\lambda=E_{1,n}(\lambda+\rho)p_\lambda$$ which is
the case $r=1$ of eq. (3.72) from \cite{vD1}.

\bigskip
 Now the relation of $\psi$ to Koornwinder polynomials
(generalized Weyl formula), the norm formula, evaluation identity
and duality, they all can be derived similar to the case of a
reduced root system. To avoid repetitions, we skip the details
(see \cite{vD2} for the formulation of
all these identities).

\section{Integrable deformation of the Macdonald--Ruijsenaars operators}

In this section we discuss a version of the Macdonald--Ruijsenaars
operators related to a certain "deformed" $A_n$ system. This
system $R=A_n(m)\subset \mathbb R^{n+1}$ was introduced in
\cite{VFC} in the following way: $R=R^0\cup R^1$ where
\begin{align} R^0=&\{\pm(e_i-e_j)\mid 1\le i<j\le
n\}\cong A_{n-1}\subset \mathbb R^n\,,\notag\\\label{dsy}
R^1=&\{\pm(e_i-\sqrt me_{n+1})\mid i=1,\dots ,
n\}\qquad\text{with}\\m_\alpha\equiv &m \quad \text{for}\
\alpha\in R^0\,,\qquad m_\alpha\equiv 1 \quad \text{for}\
\alpha\in R^1\,.\notag
\end{align}
Here $m$ is a parameter, which at first will be an integer.

\subsection{Deformed Macdonald--Ruijsenaars operator}
Let us consider the following difference operator related to the
deformed $A_n$ system:
\begin{equation}\label{dan}
\widetilde D=\widetilde a_1(z)T^{e_1}+\cdots +\widetilde
a_n(z)T^{e_n}+ \widetilde a_{n+1}(z)T^{\sqrt m e_{n+1}}\,,
\end{equation}
where $e_1,\dots,e_{n+1}$ is the standard basis in $V=\mathbb
R^{n+1}$ and the coefficients $a_i$ look as follows:
\begin{align}
&\widetilde a_i(z)= \frac {[z_i-\sqrt mz_{n+1}-\frac{m+1}{2}]}
{[z_i-\sqrt mz_{n+1}]}\prod_{j\neq i}^n \frac {[z_i-z_j-m]}
{[z_i-z_j]}\,,\notag\\\label{dan1}&\\ & \widetilde
a_{n+1}(z)=\frac {[1]} {[m]} \prod_{j=1}^n \frac {[\sqrt
mz_{n+1}-z_j-\frac{m+1}{2}]} {[\sqrt
mz_{n+1}-z_j+\frac{m-1}{2}]}\,.\notag
\end{align}
Here, as before, the square brackets are used to denote
$[a]:=q^a-q^{-a}$.

This operator is a discretization of the deformed Calogero--Moser
operator proposed in \cite{VFC}, with $m$ being the deformation
parameter. Its rational version was considered in \cite{C00,C97}.
For $m=1$ it reduces to a special case of the
Macdonald--Ruijsenaars operator $D_1$ in (\ref{ru}).

 We fix a
"positive half" of the system \eqref{dsy} as $$R_+=\{e_i-e_j\mid
1\le i<j\le n\}\cup \{e_i-\sqrt me_{n+1}\mid i=1,\dots , n\}\,.$$
 Next, we introduce
"weight lattice": $$P=\mathbb Ze_1\oplus\dots\oplus\mathbb
Ze_n\oplus\sqrt m\mathbb Ze_{n+1}\,.$$ Respectively, we will call
any finite sum $f(x)=\sum_{\nu\in P}f_\nu q^{2(\nu,x)}$ a
polynomial in $x$. We don't need coroots and coweights in this
case, and there will be no substantial difference between
variables $x,z\in \mathbb R^{n+1}$ below, so, for instance,
polynomials in $z$ are defined in the same way. As before, we will
apply the term {\it quasipolynomial} to a function
$\psi(x,z)=q^{2(x,z)}f$ where $f$ is polynomial either in $x$ or
in $z$.

Now let us consider a ring $\widetilde{\mathfrak R}$ of
polynomials $f(z)$ with the properties \eqref{axz} (with
$\alpha^\vee:=\alpha$). In our case they can be rewritten as
follows:
\begin{align}\label{daxz1}
f(z+se_i)=f(z+se_j)\quad\text{for}\quad q^{2z_i}=q^{2z_j}\  (1\le
s\le m,\quad 1\le i<j\le n)\,,\\\label{daxz2} f(z+e_i)=f(z+\sqrt
me_{n+1})\quad\text{for}\quad q^{1-m+2z_i}=q^{2\sqrt mz_{n+1}}\
(i=1,\dots, n)\,.
\end{align}

\begin{prop}\label{dinv} Deformed Macdonald--Ruijsenaars operator
\eqref{dan}--\eqref{dan1} with $m\in\mathbb Z_+$ preserves the
ring $\widetilde{\mathfrak R}$ of the polynomials with the
properties \eqref{daxz1}--\eqref{daxz2}:\  $\widetilde
D(\widetilde{\mathfrak R}) \subseteq \widetilde{\mathfrak R}$.
\end{prop}

\begin{proof} For the conditions \eqref{daxz1} the arguments
repeat those from the proof of Proposition \ref{inv}, because the
operator $D$ is symmetric in $z_1,\dots ,z_n$. For the remaining
conditions \eqref{daxz2} everything reduces to the proof that the
operators $D_0=\widetilde a_j(z)T^{e_j}$ \ (for $j\ne i$) and
$D_1=\widetilde a_i(z)T^{e_i}+\widetilde a_{n+1}(z)T^{\sqrt
me_{n+1}}$ preserve the property \eqref{daxz2}. This can be
checked straightforwardly, in the spirit of Lemma \ref{pres1}.
\end{proof}

\subsection{BA function}
Introduce  $\rho=\rho(m)$ similar to \eqref{rho1}:
\begin{equation}
\label{drho1} \rho = \frac 12 \sum_{\alpha \in R_+}m_{\alpha}
\alpha \,.
\end{equation}
In our case we have explicitly: $$\rho=\frac m2
(n-1,n-3,\dots,-n+1,0)+\frac 12 (1,\dots,1,-\sqrt mn)\,.$$ We
define a polytope $\mathcal N$ similar to \eqref{nu}:
\begin{equation}
\label{dnu} \mathcal N = \{ \nu= -\rho+\sum_{\alpha \in R_+}
l_\alpha \alpha \mid 0\le l_\alpha \le m_\alpha \}\,.
\end{equation}
Now the definition of a BA function repeats our definition in case
of a root system.

\begin{definition} A function $\psi(x,z)$ which is
quasipolynomial in $z$, $$\psi = q^{2(x,z)}\sum_{\nu\in P}
\psi_\nu q^{2(\nu,z)}\,,$$ with ${\rm supp}(\psi)\subseteq
\mathcal N$ and which satisfies the conditions
\eqref{daxz1}--\eqref{daxz2} in $z$ is called a BA function for
the system \eqref{dsy}.
\end{definition}

In exactly the same way as in Section \ref{baf}, one proves that
$\psi$ is defined uniquely up to an $x$-depending factor.
Moreover, analyzing the corresponding linear conditions for the
coefficients $\psi_\nu$, we come to the following choice of
normalizing $\psi$:
\begin{equation}
\label{dnorm0} \psi_{\rho}=\prod_{\alpha\in
R_+}\prod_{j=1}^{m_\alpha} [(\alpha,\frac12 j\alpha -x)]\,,
\end{equation}
with $[a]$ denoting the $q$-number $[a]:=q^a-q^{-a}$. We will call
this $\psi$ the {\it normalized BA function}. One can calculate
the leading coefficients $\psi_\nu$ at other vertices of the
polytope $\mathcal N$. It has $(n+1)!$ vertices which are in
one-to-one correspondence with all permutations $\sigma\in
S_{n+1}$. Namely, for $\sigma\in S_{n+1}$ let us introduce a
vector $t_\sigma\in\mathbb R^{n+1}$ as $$t_\sigma=(\sigma_1,\dots,
\sigma_n, \sigma_{n+1}/\sqrt m)\,.$$ Now denote by $\sigma R_+$
the following subset in $R$: $$\sigma R_+=\{\alpha\in R\mid
(\alpha, t_\sigma)<0\}\,,$$ and introduce $\sigma\rho$ as
$$\sigma\rho = \frac 12 \sum_{\alpha \in \sigma R_+}m_{\alpha}
\alpha \,.$$ These are exactly the vertices of the polytope
\eqref{dnu}. For instance, taking $\sigma=e$ we get
$\sigma\rho=\rho$, and taking $\sigma=(n+1,n,\dots,1)$ we get
$\sigma\rho=-\rho$.

\bigskip
Similarly to Proposition \ref{no}, one gets the following result.
\begin{prop}
The leading coefficients $\psi_{\sigma\rho}$ of the normalized BA
function have the form:
\begin{equation}
\label{dnorm} \psi_{\sigma\rho}=\prod_{\alpha\in \sigma
R_+}\prod_{j=1}^{m_\alpha} [(\alpha,\frac12 j\alpha-x)]\,.
\end{equation}
 The
normalized BA function is quasipolynomial in both $x$ and $z$.
\end{prop}

\bigskip
The existence of $\psi$ can be proven similar to Theorem
\ref{exist}. To formulate the result, introduce a polynomial
$Q(z)$ as
\begin{equation}\label{dQ}
Q(z)=q^{2(\rho,z)}\prod_{\alpha\in
R_+}\prod_{j=1}^{m_\alpha}[(\alpha,z+\frac12 j\alpha)]\,
[(\alpha,z-\frac12 j\alpha)]\,,
\end{equation}
with square brackets denoting $q$-number as before. Define also a
"deformed orbitsum" $\mathfrak m(x)$ as follows:
\begin{equation}\label{dorb}
\mathfrak m(x)=q^{2x_1}+\dots
+q^{2x_n}+\frac{q-q^{-1}}{q^m-q^{-m}}q^{2\sqrt mx_{n+1}}\,.
\end{equation}

\begin{theorem}\label{dexist}
Let $\widetilde D$ be the deformed Macdonald--Ruijsenaars operator
\eqref{dan}--\eqref{dan1}. Define $\psi(x,z)$ as follows:
\begin{equation}\label{dber}
\psi=\prod_{\nu}\left(\widetilde D-\mathfrak
m(x+\nu)\right)\left[q^{2(x,z)}Q(z) \right]\,,
\end{equation}
in accordance with the formulas \eqref{dQ}, \eqref{dorb}, where
the product is taken over all $\nu\ne 0$ having the form
$\nu=\sum_{\alpha\in R_+}l_\alpha\alpha$ with $l_\alpha=0,\dots ,
m_\alpha$. We have the following:

\noindent(i) $\psi$ is a BA function for the system \eqref{dsy};

\noindent(ii) its coefficient $\psi_{-\rho}$ equals $\prod_\nu
(\mathfrak m(x)-\mathfrak m(x+\nu))\ne 0$;

\noindent(iii) as a function of $z$,\ $\psi$ is an eigenfunction
of the operator $\widetilde D$: $\widetilde D\psi=\mathfrak
m(x)\psi$.
\end{theorem}

Thus,  renormalizing the constructed $\psi$ one gets the
normalized Baker--Akhiezer function $\Psi$ for the system
\eqref{dsy}.

Now one can derive the duality similar to Theorem \ref{symmetry}.
In this case it is simply the symmetry between $x,z$.

\begin{theorem}
The normalized BA function $\Psi$ constructed above, is symmetric
under permutation of its arguments: $\Psi(x,z)=\Psi(z,x)$.
\end{theorem}

\subsection{Quantum integrability}
We start from discussing the {\it algebraic integrability} of the
deformed Macdonald--Ruijsenaars operator \eqref{dan}--\eqref{dan1}
with $m\in \mathbb Z_+$ , which is a direct corollary of the
existence of a BA function for the system \eqref{dsy}.

\medskip
\begin{theorem}[Algebraic integrability]\label{dalin} Let
$m\in \mathbb Z_+$ and $\Psi(x,z)$ be the normalized BA function
for the deformed $A_n$ system \eqref{dsy}. Then for each
polynomial $f(z)$ from the ring $\widetilde{\mathfrak R}$ there
exists a difference operator $D_f$ in $x$ on the "weight" lattice
$P$ such that $ D_f \Psi = f(z)\Psi$. All these operators commute.
For $f=\mathfrak m(z)\in\widetilde{\mathfrak R}$ given by
\eqref{dorb}, the corresponding operator $D_f$ is the deformed
Macdonald--Ruijsenaars operator $\widetilde D$ given
\eqref{m}--\eqref{m2}.
\end{theorem}

\medskip
Proof is the same as in Theorem \ref{alin}. In this way we obtain
a commutative ring of difference operators, isomorphic to the ring
\eqref{daxz1}--\eqref{daxz2}. Moreover, due to the symmetry
between $x$ and $z$, we obtain a {\it bispectral pair} of
commutative rings with $\Psi(x,z)$ being their common
eigenfunction.

The ring $\widetilde{\mathfrak R}$ is big enough, it contains, for
instance, a principal ideal generated by the polynomial
\eqref{dQ}. Now let us choose the following special elements
$\mathfrak m_1, \dots, \mathfrak m_{n+1}$ of the ring
$\widetilde{\mathfrak R}$:
\begin{equation}\label{dorbi}
\mathfrak m_s(z)=q^{2sz_1}+\dots
+q^{2sz_n}+\frac{q^s-q^{-s}}{q^{sm}-q^{-sm}}q^{2s\sqrt
mz_{n+1}}\,.
\end{equation}
In particular, for $s=1$ we obtain the deformed orbitsum from
\eqref{dorb}. Notice that for $m=1$ these polynomials turn into
the Newton basis in the ring of symmetric functions.

In accordance with Theorem \ref{dalin}, to each $\mathfrak m_s$
corresponds a certain difference operator $D_s=D_{\mathfrak m_s}$,
and they all commute. Similar to Lemma \ref{ad}, one can show that
these operators $D_s$ admit analytic continuation in $m$ and,
thus, they give rise to a commutative family for any value of the
parameter $m$. Thus, we arrive at quantum integrability of the
deformed Macdonald--Ruijsenaars operator.

\bigskip
\begin{theorem} The operator $\widetilde D$ is completely
integrable, i.e. it can be included into a commutative family
$D_1=\widetilde D, D_2,\ldots ,D_{n+1}$ of difference operators,
which in case $m=1$ coincide with the Macdonald--Ruijsenaars
family \eqref{ru}.
\end{theorem}

\bigskip
One obtains a natural elliptic version of the operator \eqref{dan}
replacing all expressions like $[a]=q^a-q^{-a}$ in its
coefficients by their elliptic analogues $\sigma(a)$ with
$\sigma(z)=\sigma(z|\omega,\omega')$ being the Weierstrass
$\sigma$-function. For such an operator a proper version of
Proposition \ref{dinv} holds. This indicates that it is
(algebraically) integrable, too. As a concluding remark, we
mention that in a similar manner one can construct generalized
Macdonald operators for other "deformed" root systems, some of
which were presented in \cite{CFV1}. Details will appear
elsewhere.

\end{document}